\documentclass[a4paper,twoside,11pt]{article} 

\usepackage{amsmath,amsfonts,amsthm,amssymb,bbm,caption,cite,graphicx,fancyhdr,hyperref,mathtools,tikz}
\usetikzlibrary{calc}
\graphicspath{{figures/}}

\makeatletter

\@addtoreset{equation}{section}
\makeatother

\newtheorem{lem}{Lemma}[section]
\newtheorem{prop}[lem]{Proposition}
\newtheorem{thm}[lem]{Theorem}
\newtheorem{cor}[lem]{Corollary}

\theoremstyle{definition}
\newtheorem{defi}[lem]{Definition}
\newtheorem{ex}[lem]{Example}
\newtheorem{asm}[lem]{Assumption}
\theoremstyle{remark}
\newtheorem{rem}[lem]{Remark}

\newcommand{\ol}[1]{\overline{#1}}

\newcommand{\Mc}[1]{\mathcal{#1}}
\newcommand{\Mb}[1]{\mathbb{#1}}
\newcommand{\Mf}[1]{\mathfrak{#1}}
\newcommand{\Mr}[1]{\mathrm{#1}}
\renewcommand{\Cup}{\bigcup}
\renewcommand{\Cap}{\bigcap}

\newcommand{\indi}{\mathbbm{1}}
\newcommand{\zero}{\boldsymbol{0}}
\newcommand{\der}{\mathrm{d}}

\DeclareMathOperator{\supp}{supp}

\DeclareMathOperator{\Lip}{Lip}
\DeclareMathOperator{\lCla}{lCla}
\DeclareMathOperator{\uBCL}{uBCL}

\makeatletter

\newcommand{\DeclareDescCond}[3]{%
	\newcounter{#1}%
	\expandafter\renewcommand\csname the#1\endcsname{(#3)}%
	\expandafter\newcommand\csname #2\endcsname[2][]{%
		\ifhmode
		\unskip
		\par
		\fi
		\if\relax\detokenize{##1}\relax
		\stepcounter{#1}%
		\protected@edef\desccond@currentlabel{%
			\csname the#1\endcsname
		}%
		\expandafter\xdef
		\csname desccond@label@##2\endcsname
		{\desccond@currentlabel}%
		\expandafter\gdef
		\csname desccond@primes@##2\endcsname{}%
		\item[\desccond@currentlabel]%
		\@ifundefined{phantomsection}{}{\phantomsection}%
		\protected@edef\@currentlabel{\desccond@currentlabel}%
		\def\@currentcounter{#1}%
		\label{##2}%
		\else
		\@ifundefined{desccond@label@##1}{%
			\PackageError{DescCond}%
			{Undefined parent condition `##1'}%
			{The parent condition must be declared before its variants.}%
		}{%
			\expandafter\g@addto@macro
			\csname desccond@primes@##1\endcsname
			{'}%
			\protected@edef\desccond@currentlabel{%
				\csname desccond@label@##1\endcsname
				\csname desccond@primes@##1\endcsname
			}%
			\item[\desccond@currentlabel]%
			\@ifundefined{phantomsection}{}{\phantomsection}%
			\protected@edef\@currentlabel{\desccond@currentlabel}%
			\def\@currentcounter{#1}%
			\label{##2}%
		}%
		\fi
		\ignorespaces
	}%
}

\makeatother

\makeatother

\DeclareDescCond{measurecond}{Mcond}{M\arabic{measurecond}}
\DeclareDescCond{measurecond2}{mcond}{m\arabic{measurecond2}}
\DeclareDescCond{functcond}{Fcond}{F\arabic{functcond}}
\DeclareDescCond{modularcond}{rcond}{R\arabic{modularcond}}
\DeclareDescCond{medmcond}{DMcond}{D\arabic{medmcond}}
\DeclareDescCond{rescond}{RFcond}{RF\arabic{rescond}}

\title{The strongly local part of nonlinear Dirichlet structures and its energy measures}
\author{K\^ohei Sasaya
	\thanks{Graduate School of Mathematical Sciences, the University of Tokyo, Tokyo 153-8914, Japan}}
\date{\today}
\begin{document}

\maketitle 
\begin{abstract}
This paper substantially extends the results of~\cite{Sas2026}:
we propose a framework for nonlinear generalizations of regular Dirichlet forms and prove that each energy in this framework has a largest strongly local part.
We also prove the existence and uniqueness of canonical energy measures of this part, which satisfy convexity and locality properties and, under additional assumptions, certain chain rules in one variable.
In particular, our framework requires neither $p$-homogeneity nor strong locality.
Moreover, under a Ball--Carlen--Lieb-type upper bound~\cite{BCL1994}, we prove more detailed chain rules for the directional derivatives of these energy measures.
Our framework is applicable both to $p$-energy forms on fractals and to a class of symmetric nonlinear Dirichlet forms in the sense of Cipriani and Grillo~\cite{CG2003}, including those of Orlicz--Sobolev and variable-exponent Sobolev type.
\end{abstract}
\tableofcontents
\section{Introduction}
In this paper, we consider possibly nonlinear generalizations of (regular) Dirichlet forms and the associated energy measures of their strongly local parts. Namely, for a locally compact separable metrizable space $X$ and a full-support Radon measure $\Mf{m}$, a regular Dirichlet form on the measure space $(X,\Mf{m})$ is a pair consisting of a functional $\Mc{E}:\Mc{F}\to[0,\infty)$ and its domain $\Mc{F}\subset L^2(X,\Mf{m})$ with certain properties. Dirichlet forms were introduced by Beurling and Deny~\cite{BD1958,BD1959} as a generalization of the pair consisting of the Dirichlet energy $\int_\Omega |\nabla f|^2\der x$ and its domain, a Hilbert space equipped with the Sobolev norm, in order to develop a potential theory independent of specific kernels. As Beurling and Deny anticipated, Dirichlet forms have been extensively used in the study of harmonic analysis, differential equations, and probability theory. Their nonlinear generalizations, the main subject of this paper, have also been studied in this century as natural extensions. More specifically, we consider generalizations of the $p$-Dirichlet energy
\[\Mc{E}(f)=\int_{\Omega}|\nabla f|^p\der x\quad (\Omega\subset\Mb{R}^n\text{: domain},\ p\in[1,\infty)),\]
its domain, and the associated localizations
\[\Gamma^c\langle f\rangle(A)=\int_A|\nabla f|^p\der x \quad (A\subset\Omega),\]
which we call energy measures, with particular emphasis on the case $p\ne2$.

\begin{figure}
\centering
\includegraphics[width=5cm]{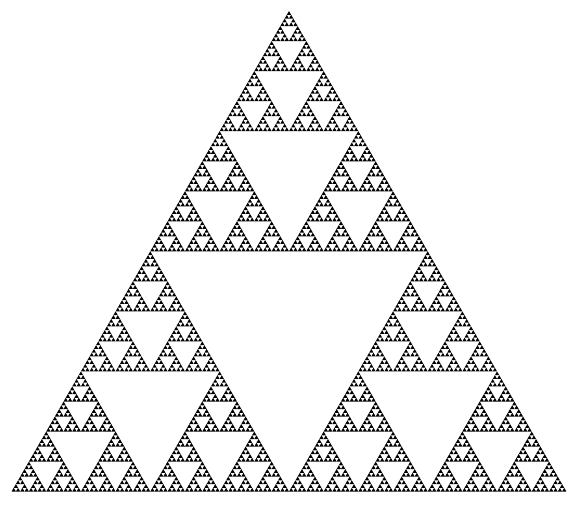}
\caption{Standard Sierpi\'nski gasket} \label{fig:SG}
\end{figure}

In contrast to the case of Dirichlet energies in $\Mb{R}^n$, the existence of energy measures and expressions of this kind is not obvious, even in the bilinear case, in nonsmooth settings such as fractals. For example, Dirichlet forms on the standard Sierpi\'nski gasket~(Figure \ref{fig:SG}) were constructed via the corresponding stochastic processes studied in \cite{BP1988,Gol1987,Kus1987}, and the mutual singularity of any energy measure and the standard Hausdorff measure was shown in \cite{Kus1989} (cf. \cite{Hin2005,HN2006,KM2020}). For $p\ne2$, counterparts of the $p$-Dirichlet energy on fractals were first considered in \cite{HPS2004}. Both their construction (e.g.~\cite{CGQ2022,Kig2023,MS2025}) and their structure (e.g.~\cite{KS2024+,Sas2026}) have been extensively studied in recent years, especially in strongly local settings, that is, when $\Mc{E}(f+g)=\Mc{E}(f)+\Mc{E}(g)$ holds whenever $f,g\in C_c(X)$ belong to the domain of $\Mc{E}$ and the supports of $f$ and $g+a$ do not intersect for some $a\in\Mb{R}$: here, $C_c(X)$ denotes the space of compactly supported real-valued continuous functions on $X$. In particular, our prior work \cite{Sas2026} shows the existence and uniqueness of energy measures associated with $p$-Dirichlet energies under axiomatic assumptions, extending the results in \cite{Sil1974,LeJ1978} for Dirichlet forms, that is, for $p=2$.

In settings without strong locality, we do not even have a standard expression for nonlinear energies on which to base local energy measures. For classical regular Dirichlet forms, the Beurling--Deny decomposition
\begin{equation}\label{eq:BDdecomp}
\Mc{E}(f)=\Mc{E}^{(c)}(f)+\int_{X\times X\setminus \Mr{diag}}(f(x)-f(y))^2 J(\der x, \der y)+\int_X f(x)^2\Mf{k}(\der x)
\end{equation}
holds for some strongly local energy $\Mc{E}^{(c)}$ and positive Radon measures $J$ on $X\times X\setminus \Mr{diag}$ and $\Mf{k}$ on $X$, and for any $f\in\Mc{F}\cap C_c(X)$, where $\Mr{diag}=\{(x,x)\mid x\in X\}$ (see \cite[Theorem 3.2.1]{FOT2011}). This kind of decomposition is also assumed in recent studies of nonlinear energies (e.g.~\cite{CY2025+,Yan2026+}), but has not been established and may even fail in general axiomatic settings. Possibly nonlocal, nonlinear energies are also defined without explicit kernel representations, as in the case of traces of strongly local energies in \cite{KS2024+}. Thus, whether a counterpart of \eqref{eq:BDdecomp} holds and what the precise decomposition formula for nonlinear energies should be are significant questions in the study of nonlinear generalizations of Dirichlet forms.

In this paper, as a first step toward resolving the problem above, we intrinsically identify the counterpart of $\Mc{E}^{(c)}$ and the associated canonical energy measures in nonlinear settings without assuming $p$-homogeneity of the energies. We first work without a reference measure on the core $\Mc{M}(\Mc{E})\cap C_c(X)$, where $\Mc{M}(\Mc{E})$ is the modular space defined in Definition~\ref{defi:modref}, and then extend the construction to an $L^q$ setting under additional assumptions. We assume that the functional $\Mc{E}$ satisfies properties \ref{cond:rconv}--\ref{cond:prereg}, which are counterparts of those of regular Dirichlet forms and include convexity, symmetry, certain contraction properties, and the vanishing of energies of certain truncations in the limit; see Section \ref{sec:results} for details. On the core, we prove the existence of a largest strongly local part $\Mc{E}^{(c)}$ in the sense that $\Mc{E}^{(c)}$ is strongly local and satisfies \ref{cond:rconv}--\ref{cond:prereg} and $\Mc{E}^{(c)}\le\Mc{E}$, and that if $\tilde{\Mc{E}}$ also satisfies the preceding conditions, then $\tilde{\Mc{E}}\le\Mc{E}^{(c)}$. In particular, if $\Mc{E}$ is $p$-homogeneous, then so is $\Mc{E}^{(c)}$. Moreover, we prove the existence and uniqueness of canonical energy measures $\{\Gamma^c\langle f\rangle\}$ such that $\Gamma^c\langle f\rangle(X)=\Mc{E}^{(c)}(f)$, $\Gamma^c\langle f\rangle(A)=\Gamma^c\langle g\rangle(A)$ whenever $f,g\in\Mc{M}(\Mc{E})\cap C_c(X)$ and $f-g$ is constant on $A$, and $f\mapsto \Gamma^c\langle f\rangle(A)$ satisfies \ref{cond:rconv}--\ref{cond:prereg} for any Borel set $A$.

In the $L^q$ setting of Theorem~\ref{thm:main2}, the structural conditions are imposed on a dense linear core $\Mc{C}\subset\Mc{M}(\Mc{E})\cap C_c(X)$. Reflexivity and suitable bounds under scalar multiplication yield chain rules for energy measures on $\Mc{C}$. In this setting, $p$-homogeneity gives the usual exact chain rule. A Ball--Carlen--Lieb-type upper bound~\cite{BCL1994} further allows us to define directional derivatives of energy measures as signed Borel measures. Under additional assumptions, we establish chain rules for these derivatives in both arguments; see Theorem~\ref{thm:chain2} (cf.~\cite{KS2024+,Sas2026}).

Our results essentially recover those in \cite[Theorem 2.4]{Sas2026} and relax the assumptions of strong locality, the $p$-Clarkson inequalities,\footnote{See Remark \ref{rem:clk}~(4) for the definition.} and $p$-homogeneity for the existence and uniqueness results. In particular, to the best of the author's knowledge, we give the first constructions of energy measures of the strongly local parts of possibly nonlocal $p$-energies with $p\ne2$ in cases where these parts are not apparent from the defining constructions; see Remark \ref{lem:pemapply}~(2),(3). Relaxing the $p$-Clarkson inequalities is also useful when $\Mc{E}$ is $p$-homogeneous (cf. \cite[Example 2.9]{EM2025+}).

Our results are also applicable to a certain class of nonlinear Dirichlet forms in the sense of Cipriani and Grillo~\cite{CG2003}, which provide another framework for nonlinear generalizations of Dirichlet forms; see Definition~\ref{defi:NLDF} and, for related developments, \cite{BH2024,Bri2026,BDS2025+}. For example, $\Mc{E}:L^2(0,1)\to[0,\infty]$ defined by
\[\Mc{E}(f)=
\begin{cases}
\int_0^1 (|f'|^2+|f'|^3)\der x&f\in W^{1,3}_0(0,1)\\
\infty&\text{otherwise},
\end{cases}
\]
which is inhomogeneous, is a nonlinear Dirichlet form in this sense.
This framework also includes examples from Orlicz--Sobolev spaces and Sobolev spaces with variable exponents, as well as gradient energies on sub-Rie\-mannian and reversible Finsler manifolds. Localization and finer calculus rules are discussed in the section titled ``Perspectives and desiderata'' in the latest survey~\cite{Bri2026} in this area. Our results therefore contribute to the development of these studies, although the scope of our theory is subject to some limitations.

Here we briefly explain the central idea for obtaining $\Mc{E}^{(c)}$. Let $(\Mc{E},\Mc{F})$ be a regular Dirichlet form, and let $T_n:\Mb{R}\to\Mb{R}$ be the triangle wave defined by $T_n(t)=\min_{k\in\Mb{Z}}|t-2^{1-n}k|$. Then, for $f\in\Mc{F}\cap C_c(X)$, by the chain rule (see \cite[Theorem 3.2.2 and footnote 8 on p.~129]{FOT2011}) we have $\Mc{E}^{(c)}(T_n\circ f)=\Mc{E}^{(c)}(f)$ and
\begin{multline*}
\lim_{n\to\infty}\int_{X\times X\setminus \Mr{diag}}|T_n\circ f(x)-T_n\circ f(y)|^2J(\der x, \der y)\\
=\lim_{n\to\infty}\int_X|T_n\circ f(x)|^2\Mf{k}(\der x)=0
\end{multline*}
by the dominated convergence theorem, where $J,\Mf{k}$ are the measures in \eqref{eq:BDdecomp}. Hence, we obtain
\begin{equation}\label{eq:limenergy}
	\Mc{E}^{(c)}(f)=\lim_{n\to\infty}\Mc{E}^{(c)}(T_n\circ f)=\lim_{n\to\infty}\Mc{E}(T_n\circ f).
\end{equation}
For nonlinear energies, we seek to use the equality \eqref{eq:limenergy} to define $\Mc{E}^{(c)}$.

We first define Radon measures $\{\Pi\langle f\rangle\}$ on $X$, which we call pre-energy measures, using a construction similar to that in \cite{Sas2026}. More precisely, for fixed $f,g\in\Mc{M}(\Mc{E})\cap C_c(X)$ and $a\in\Mb{R}$, we define a sequence of functions $\{f_{n}^{g,a}\}_{n\in\Mb{N}}$, where $f_n^{g,a}$ is the minimum of $T_n\circ f$ and a suitable truncation of $g$ (see Section \ref{sec:pushforward} for the precise definition), and define the Radon measure $\Pi\langle f\rangle$ so that $\Pi\langle f\rangle(\{g\le a\})=\lim_{n\to\infty}\Mc{E}(f_n^{g,a})$.\footnote{For an intuitive understanding, see Lemma \ref{lem:fundestm}.} In particular, we obtain $\Pi\langle f\rangle(X)=\lim_{n\to\infty}\Mc{E}(T_n\circ f)$. However, we do not know whether $f\mapsto \Pi\langle f\rangle(A)$ is convex. Hence, we define the measures $\Gamma^c\langle f\rangle$ as ``convex envelopes'' of the pre-energy measures $\{\Pi\langle f\rangle\}$, set $\Mc{E}^{(c)}(f)=\Gamma^c\langle f\rangle (X)$, and establish their properties and maximality (Theorem \ref{thm:main1}). In general, whether $\Pi=\Gamma^c$ holds remains open; see Section \ref{sec:arg}.

Admittedly, we do not resolve the initial question of a nonlinear Beurling--Deny decomposition. However, we believe that our structural results for nonlinear generalizations of Dirichlet forms contribute to further studies in this direction.

The organization of this paper is as follows. Section~\ref{sec:results} is devoted to introducing our framework and results. Section~\ref{sec:modular} presents preliminary results on generalized semimodulars, a key notion for handling inhomogeneous energies in this paper. In Section~\ref{sec:pushforward}, we study the properties of Radon measures on $\Mb{R}$ that can be viewed as the pushforward measures $g_*\Pi\langle f\rangle$. In Section~\ref{sec:construction}, we construct $\Pi\langle f\rangle$ and $\Gamma^c\langle f\rangle$ for $f\in\Mc{M}(\Mc{E})\cap C_c(X)$ and prove the first form of our main theorem. Under lower semicontinuity in $L^q$, Section~\ref{sec:lsc} extends $\Mc{E}^{(c)},\Gamma^c$ beyond the core and, with additional assumptions, establishes chain rules in one variable, thereby completing the proof of the second form of our main theorem. We discuss reflexivity of the associated function spaces and more detailed chain rules for the resulting energy measures under Clarkson-type and Ball--Carlen--Lieb-type inequalities in Section~\ref{sec:CI}. Sections~\ref{sec:ex} and~\ref{sec:arg} are devoted to examples to which our results apply and remaining problems, respectively. Appendix~\ref{ap:further} discusses notable relations among the conditions that appear in our framework and results, and Appendix~\ref{sec:medm} establishes the existence of a minimal energy-dominant measure (cf.~\cite{Hin2010}).



\subsection*{Notation}
\begin{itemize}
	\item For a set $X$, $\Mb{R}^X$ denotes the real-valued functions on $X$; $\indi_A$ and $\sqcup$ denote the indicator of $A$ and disjoint union. We use $\vee$ and $\wedge$ for maximum and minimum, pointwise for functions, with $a\vee f\wedge b=(a\vee f)\wedge b$ for $a\le b$.

	\item On a topological space $X$, $\ol{U}$ and $\Mc{B}(X)$ denote closure and the Borel $\sigma$-algebra; $U\Subset V$ means that $\ol{U}$ is compact and contained in $V$, and $\supp(f)=\ol{\{x\in X:f(x)\ne0\}}$ for continuous $f$ on $X$. We write $C(X)$ for real-valued continuous functions, with $C_c(X)$ and $C_b(X)$ the subspaces of compactly supported and bounded functions, respectively. We set $\|f\|_\infty=\sup_{x\in X}|f(x)|$ for $f\in C_b(X)$.

	\item For a metric $d$ on $X$, $B_d(x,r)=\{y\in X:d(x,y)<r\}$ for $r>0$, $d(A,B)=\inf\{d(x,y):x\in A,\ y\in B\}$ for nonempty subsets $A,B$ of $X$, and $C(X,d)$ denotes $C(X)$ when $X$ is equipped with the topology induced by $d$. For a map $f:(X,d)\to(Y,d_Y)$ between metric spaces, we set
	\[
		\Lip(f)=\inf\{L\ge0:d_Y(f(x),f(y))\le Ld(x,y)\text{ for all }x,y\in X\}.
	\]

	\item For (possibly signed) measures $\mu,\nu$ on the same measurable space $X$, $|\mu|$ denotes the total variation measure. We write $\mu\ll\nu$ for absolute continuity of $|\mu|$ with respect to $|\nu|$, and $\mu\le\nu$ if $\mu(A)\le\nu(A)$ for every measurable set $A$.
	For a measurable map $f$, $f_*\mu$ denotes the pushforward of $\mu$; $\mu|_A$ and $\mu|_{\Mc{A}}$ denote its restrictions to a measurable set $A$ and a sub-$\sigma$-algebra $\Mc{A}$, respectively.
	For $h\in L^1(|\mu|)$, $h\der\mu$ denotes the signed measure $A\mapsto\int_A h\der\mu$; this notation also applies to nonnegative measurable $h$ when $\mu\ge0$.
	If $X$ is topological and $\mu$ is a (signed) Borel measure, its support is
	\[
		\supp[\mu]=\{x\in X:|\mu|(U)>0\text{ for every open }U\ni x\}.
	\]
	\item $\Mc{L}^1$ denotes one-dimensional Lebesgue measure. We always assume piecewise affine functions and piecewise $C^1$ functions on $\Mb{R}$ are continuous.

	\item We write $\partial_iP=\partial P/\partial x_i$ for the $i$th partial derivative and $\nabla f$ for the weak gradient of a Sobolev function $f$ (the weak derivative in one dimension). The space $W_0^{1,p}(\Omega)$ is the closure of $C_c^\infty(\Omega)$ in $W^{1,p}(\Omega)$.
\end{itemize}
We define some basic functions on $\Mb{R}$ by
\begin{align*}
	\eta_{[a,b]}(t)&=(a\vee t\wedge b)-(a\vee 0\wedge b) & S_n^a(t)&=0\vee(a+2^{1-n}-t)\wedge 2^{-n} \\
	T_{n,\alpha}(t)&=\min_{k\in\Mb{Z}}|t-k\alpha2^{1-n}|,  &T_n(t)&=T_{n,1}(t)
\end{align*}
for $a,b\in\Mb{R}$ with $a\le b$, $\alpha\in(0,\infty)$ and $n\in\Mb{N}$.

\begin{rem}
Let $X$ be Hausdorff and $\mu$ a (signed) Borel measure with $\supp[\mu]=X$. Then $\supp(f)=\supp[f\der\mu]$ for every $f\in C(X)$ for which $f\der\mu$ is defined.
\end{rem}

\section{Framework and main results}\label{sec:results}
For the treatment of $p$-energy forms and nonlinear Dirichlet forms in a common framework, we use the following generic notion.
\begin{defi}[{Generalized semimodular (on an $\Mb{R}$-vector space),\cite{SZ2025+}}]\label{defi:gsm}
Let $V$ be an $\Mb{R}$-vector space. We say that $\varrho:V\to[0,\infty]$ is a \emph{generalized semimodular} if $\varrho$ satisfies the following conditions.
\begin{description}
	\rcond{cond:rconv} $\varrho$ is convex.
	\rcond{cond:rsym} $\varrho$ is symmetric, i.e., $\varrho(v)=\varrho(-v)$ and $\varrho(0)=0$.
	\rcond{cond:rlc} $\lim_{\lambda\uparrow 1}\varrho(\lambda v)=\varrho(v)$ holds for any $v\in V$.
\end{description}
\end{defi}
The notion of a modular in this sense was first introduced by Nakano~\cite{Nak1950,Nak1951}, and its generalizations have been used in the studies of Orlicz spaces and other generalizations of Lebesgue spaces. For standard references, see \cite{Mus1983} or \cite{DHHR2011}. We also note that there are subtle differences among definitions of modular-like functions. In this paper, we adopt the one in \cite{SZ2025+}, for the convenience of citations.\footnote{For the equivalence between our formulation and that in \cite{SZ2025+}, see Section \ref{sec:modular}.} Moreover, we abbreviate generalized semimodular as \emph{gsm}, as in \cite{SZ2025+}.
\begin{defi}\label{defi:modref}
Let $\varrho$ be a gsm on $V$. We write
\[\Mc{D}(\varrho)=\{v\in V\mid \varrho(v)<\infty\},\quad \Mc{M}(\varrho)=\{v\in V\mid \lambda v\in\Mc{D}(\varrho)\text{ for some }\lambda>0   \}\]
and 
\[\|v\|_{LN(\varrho)}=\inf\{\lambda>0\mid \varrho (\lambda^{-1}v)\le 1 \}\quad \text{for }v\in\Mc{M}(\varrho).\]
Then, it is easily shown that $(\Mc{M}(\varrho),\|\cdot\|_{LN(\varrho)})$ is a seminormed vector space (cf.~\cite[Theorem 2.1.7]{DHHR2011}). $\Mc{D}(\varrho)$ (resp. $\Mc{M}(\varrho)$) is called the \emph{effective domain} (resp. \emph{modular space}) of $\varrho$. The seminorm $\|v\|_{LN(\varrho)}$ is called the  \emph{Luxemburg--Nakano norm}, or simply \emph{Luxemburg norm}, of $\varrho$.

For $p\ge1$, $\varrho$ is called $p$-homogeneous if $\varrho(a v)=|a|^p\varrho(v)$ for any $a\in\Mb{R}$ and $v\in V$. $\varrho$ is called \emph{reflexive} (as a gsm) if the completion of $(\Mc{M}(\varrho)/\{\|\cdot\|_{LN(\varrho)}=0\}, \|\cdot\|_{LN(\varrho)})$ is reflexive.
\end{defi}
\begin{defi}\label{defi:condgsm}
Let $X$ be a (nonempty) topological space, $V$ be a linear subspace of $\Mb{R}^X$ and $\varrho$ be a gsm on $V$. We define the properties of $\varrho$ as follows.
\begin{description}
		\rcond{cond:rwNC}(Weak normal contraction property)  Let $\varphi:\Mb{R}\to\Mb{R}$ be a piecewise affine function with $\varphi(0)=0$, and $|\varphi'(t)|\le 1\ \Mc{L}^1$-a.e.~$t$. Then $\varphi\circ f\in\Mc{M}(\varrho)$ and $\varrho(\varphi\circ f)\le\varrho(f)$ whenever $f\in\Mc{M}(\varrho)$.
	\rcond{cond:rvanish}(Vanishing truncation property) $0\vee f\wedge t\in\Mc{M}(\varrho)$ for any $t>0$ and $\lim_{t\downarrow0}\varrho(0\vee f\wedge t )=0$, whenever $f\in\Mc{M}(\varrho)\cap C_c(X)$.
		\rcond{cond:rstrsubad}(Strong subadditivity) If $f, g\in\Mc{M}(\varrho)$, then {\rm(}$f\vee g, f\wedge g\in V$ and{\rm)} $\varrho(f\vee g)+\varrho(f\wedge g)\le \varrho(f)+\varrho(g)$.
	\rcond{cond:prereg}(Regularity) $\Mc{M}(\varrho)\cap C_c(X)$ is dense in $(C_c(X),\|\cdot\|_{\infty})$.
		\rcond{cond:rstrloc}(Strong local property) Assume that $f,g\in \Mc{M}(\varrho)\cap C_c(X)$, and that $\supp(f)\cap\supp(g+a)=\emptyset$ for some $a\in\Mb{R}$. Then, $\varrho(f+g)=\varrho(f)+\varrho(g)$.
\end{description}
\end{defi}
\ref{cond:rwNC}, \ref{cond:rstrsubad}--\ref{cond:rstrloc} are counterparts of corresponding properties of Dirichlet forms. The name of \ref{cond:rstrsubad} comes from \cite{KS2024+}, and this condition is also called \emph{submodularity} in \cite{GM2015,Bri2026}. The terminology for \ref{cond:prereg} follows the definition of regular resistance forms in \cite[Definition 6.2]{Kig2012}. \ref{cond:rvanish} is somewhat technical and serves as an alternative to lower semicontinuity without underlying $L^q$ spaces.

Our first form of the main result does not require any reference measures. 
\begin{thm}\label{thm:main1}
Let $X$ be a (nonempty) locally compact metrizable space and $\Mc{E}$ be a gsm on $V\subset C_c(X)$, satisfying \ref{cond:rwNC}--\ref{cond:prereg}.
	\begin{enumerate}
		\item There exists a maximum gsm $\Mc{E}^{(c)}$ on $V$ satisfying \ref{cond:rwNC}--\ref{cond:rstrloc} and $\Mc{E}^{(c)}\le \Mc{E}$. In other words, if a gsm $E$ on $V$ satisfies \ref{cond:rwNC}--\ref{cond:rstrloc} and $E\le \Mc{E}$, then $E\le\Mc{E}^{(c)}$.
		\item There exists a unique family $\{\Gamma^c\langle f\rangle\}_{f\in\Mc{M}(\Mc{E})}$ of finite Radon measures on $X$ with the following properties.
		\begin{description}
			\Mcond{cond:mene}$\Mc{E}^{(c)}(f)=\Gamma^c\langle f\rangle (X)$ for any $f\in\Mc{M}(\Mc{E})$.
			\Mcond{cond:mgsm} For each $A\in\Mc{B}(X)$, the functional $f\mapsto \Gamma^c\langle f\rangle(A)$ is a gsm on $\Mc{M}(\Mc{E})$.
			\Mcond{cond:mlocal} $\Gamma^c\langle f\rangle(A)=\Gamma^c\langle g\rangle (A)$ whenever $f,g\in\Mc{M}(\Mc{E})\cap C_c(X)$ and $f-g$ is constant on $A$.
		\end{description}
		\item If $\Mc{E}$ is $p$-homogeneous for some $p\ge1$, then $\Mc{E}^{(c)}$ and $f\mapsto \Gamma^c\langle f\rangle(A)$ are so, for any $A\in\Mc{B}(X)$. Moreover, the triangle inequality 
		\begin{equation}\label{eq:tri}
		\Gamma^c\langle f\rangle (A)^\frac{1}{p}+\Gamma^c\langle g\rangle (A)^\frac{1}{p}\ge\Gamma^c\langle f+g\rangle (A)^\frac{1}{p}
		\end{equation}
		holds for any $f,g\in\Mc{M}(\Mc{E})$ and $A\in\Mc{B}(X)$.
	\end{enumerate}
	
\end{thm}
We call $\Mc{E}^{(c)}$ the \emph{strongly local part} of $\Mc{E}$, and also call $\{\Gamma^c\langle f\rangle\}_{f\in\Mc{M}(\Mc{E})}$ the \emph{canonical energy measures} of $\Mc{E}^{(c)}$.

For practical use, we also give another form of the main result with a reference measure.
\begin{rem}\label{rem:abbrinj}
If $\Mf{m}$ is a full-support Radon measure on a locally compact metrizable space $X$, then there exists a natural injection $\iota:C_c(X)\to L^q(\Mf{m})$ for $q\ge 1$. Identifying $\iota(f)$ with $f$, the properties \ref{cond:rwNC}--\ref{cond:rstrloc} are also well-defined for a gsm on (a subspace of) $L^q(\Mf{m})$. We use this extended notion in Theorem \ref{thm:main2}.
\end{rem}
\begin{thm}\label{thm:main2}
Let $q\in[1,\infty)$, $X$ be a (nonempty) locally compact metrizable space, $\Mf{m}$ be a full-support Radon measure on $X$ and $\Mc{E}$ be a lower semicontinuous gsm on $L^q(\Mf{m})$. Assume that there exists an $\Mb{R}$-linear subspace $\Mc{C}\subset\Mc{M}(\Mc{E})\cap C_c(X)$ such that $\Mc{E}|_{\Mc{C}}$ satisfies \ref{cond:rwNC}--\ref{cond:prereg} and $\Mc{C}$ is dense in $(\Mc{M}(\Mc{E}), \|\cdot\|_{LN(\Mc{E})}+\|\cdot\|_{L^q(\Mf{m})})$.
\begin{enumerate}
	\item There exists a maximum gsm $\Mc{E}^{(c)}$ on $L^q(\Mf{m})$ such that $\Mc{E}^{(c)}\le\Mc{E}$ and $\Mc{E}^{(c)}|_{\Mc{C}}$ satisfies \ref{cond:rwNC}--\ref{cond:rstrloc}. If $\Mc{E}$ is reflexive, then $\Mc{E}^{(c)}$ also satisfies \ref{cond:rwNC} and \ref{cond:rstrsubad}.
	\item There exists a unique family $\{\Gamma^c\langle f\rangle\}_{f\in\Mc{M}(\Mc{E})}$ of finite Radon measures on $X$ satisfying \ref{cond:mene}--\ref{cond:mlocal}, with $f,g\in\Mc{C}$ in \ref{cond:mlocal}, and the following additional property.
	\begin{description}
		\Mcond{cond:mconti} For any $A\in\Mc{B}(X)$, $f\mapsto \Gamma^c\langle f\rangle(A)$ is continuous on $\Mc{M}(\Mc{E})$ with respect to the seminorm $\|\cdot\|_{LN(\Mc{E})}$.
	\end{description}
	\item If $\Mc{E}$ is $p$-homogeneous for some $p\ge1$, then $\Mc{E}^{(c)}$ and $f\mapsto \Gamma^c\langle f\rangle(A)$ are so, for any $A\in\Mc{B}(X)$. Moreover, the triangle inequality \eqref{eq:tri} holds for any $f,g\in\Mc{M}(\Mc{E})$ and $A\in\Mc{B}(X)$.
	\item(Chain rule). If $\Mc{E}$ is reflexive and there exist continuous functions $\Phi_1, \Phi_2:\Mb{R}\to[0,\infty)$ such that
	\[ \Phi_1(t)\Mc{E}(f)\le\Mc{E}(tf)\le\Phi_2(t)\Mc{E}(f) \]
	for any $t\in\Mb{R}$ and $f\in\Mc{C}$, then
	\begin{equation}\label{eq:chaininthm2}
		\bigl(\Phi_1\circ(\varphi')\circ f\bigr)\der\Gamma^c\langle f\rangle
		\le \Gamma^c\langle\varphi\circ f\rangle
		\le \bigl(\Phi_2\circ(\varphi')\circ f\bigr)\der\Gamma^c\langle f\rangle
	\end{equation}
	for any $f\in\Mc{C}$ and locally Lipschitz function $\varphi$ with $\varphi(0)=0$.
\end{enumerate}
\end{thm}

For examples to which our framework applies, see Section \ref{sec:ex}. For more detailed chain rules for specific cases (and the well-definedness of \eqref{eq:chaininthm2}), see Sections \ref{sec:lsc} and \ref{sec:CI}. In particular, our results essentially recover \cite[Theorem 2.4]{Sas2026}.

\section{Basic properties of generalized semimodulars}\label{sec:modular}
As a preliminary step, we introduce the basic properties of generalized semimodulars. We only focus on those which we use in this paper. For the details of proofs and other properties, see \cite[Appendix A]{SZ2025+} or \cite[Section 2]{DHHR2011} (with the following Corollary \ref{cor:semimodular}), for example. Throughout this section, we assume that $\varrho:V\to[0,\infty]$ is defined on an $\Mb{R}$-vector space $V$ and satisfies \ref{cond:rconv} and \ref{cond:rsym}; we start by discussing a condition equivalent to \ref{cond:rlc}.
\begin{lem}\label{lem:rlcalter}
	The following condition \ref{cond:rlcalter} is equivalent to \ref{cond:rlc} (under \ref{cond:rconv} and \ref{cond:rsym}).
	\begin{description}
		\rcond[cond:rlc]{cond:rlcalter} $\lim_{\lambda\uparrow 1}\varrho(\lambda v)=\varrho(v)$ whenever $\sup\{\lambda\ge0\mid \varrho(\lambda v)<\infty\}=1$.
	\end{description}
\end{lem}
\begin{proof}
	Since \ref{cond:rconv} and \ref{cond:rsym} hold, for any fixed $v\in V$, $I_v\coloneq \{\lambda \ge0\mid \varrho(\lambda v)<\infty\}$
	is $\{0\}$ or an interval including $0$. Moreover, $\varrho(\lambda v)$ is continuous on the interior of $I_v$ by \ref{cond:rconv} (see \cite[Theorem 6.2]{HL2001} for example). Hence, \ref{cond:rlcalter} implies \ref{cond:rlc}. The other implication is obvious.
\end{proof}
Hereafter, we assume $\varrho$ to be a gsm. We first note that our Definition \ref{defi:gsm} is equivalent to \cite[Definition A.1]{SZ2025+}, because $0\in\Mc{D}(\varrho)\ne\emptyset$ by \ref{cond:rsym}. 

A special case of a gsm is defined as follows.
\begin{defi}[semimodular]
$\varrho$ is called a semimodular if it satisfies the following condition.
\begin{description}
	\rcond{cond:rsemimod} $\varrho(\lambda v)=0$ for any $\lambda>0$ if and only if $v=0$.
\end{description}
\end{defi}
Since $\varrho (\lambda v)\le \lambda\varrho(v)+(1-\lambda)\varrho(0)=\lambda\varrho(v)$ for $\lambda\in(0,1)$, it is evident that $\|\cdot\|_{LN(\varrho)}$ is a norm if and only if $\varrho$ is a semimodular.
\begin{lem}[{\cite[Corollary A.4]{SZ2025+}}]\label{lem:lsc}
	$\varrho$ is lower semicontinuous on $(\Mc{M}(\varrho),\|\cdot\|_{LN(\varrho)})$.
\end{lem}
\begin{cor}\label{cor:semimodular}
$\varrho: V/\{\|\cdot\|_{LN(\varrho)}=0\}\to [0,+\infty]$ is well-defined, and is a semimodular on $V/ \{\|\cdot\|_{LN(\varrho)}=0\}$.
\end{cor}
Thanks to Corollary \ref{cor:semimodular}, we can derive the properties of gsms from those of semimodulars.

We next discuss the upper semicontinuity of $\varrho$.
\begin{lem}[{\cite[Lemma A.3]{SZ2025+}}]\label{lem:continuity}
Let $\{v_n\}_{n\ge 1}$ be a sequence of vectors in $\Mc{M}(\varrho)$. Then $\lim_{n\to\infty}\varrho(\lambda v_n)=0$ for any $\lambda>0$ if and only if
$v_n\to 0\ (n\to\infty)$ in $\|\cdot\|_{LN(\varrho)}$.
\end{lem}
\begin{defi}
We define
\[\Mc{D}^i(\varrho)=\{v\in V\mid \varrho(\lambda v)<\infty \text{ for some }\lambda>1 \}.\]
\end{defi}
\begin{lem}\label{lem:limeval}
Let $\{u_n\}_{n\ge 1},\{v_n\}_{n\ge 1}$ be sequences of vectors in $V$ such that there exist $\lambda_*>1$ and $u\in \Mc{D}^i(\varrho)$ satisfying $\varrho(\lambda u_n)\le \varrho(\lambda u)$ for any $n$ and $\lambda\in(1,\lambda_*)$, and $v_n\to 0\ (n\to\infty)$ in $\|\cdot\|_{LN(\varrho)}$. Then, $u_n+v_n\in\Mc{D}^i(\varrho)$ for sufficiently large $n$ and $\limsup_{n\to\infty}\varrho(u_n+v_n)\le\varrho(u)$. In particular, $\varrho$ is continuous at each $u\in\Mc{D}^i(\varrho)$.
\end{lem}
\begin{proof}
Choose $\lambda_0\in(1,\lambda_*)$ such that $\varrho(\lambda_0u)<\infty$. Since $t\mapsto\varrho(tu)$ is continuous on $(0,\lambda_0)$ (see \cite[Theorem 6.2]{HL2001} for example), for any $\epsilon>0$, there exist $\delta\in(0,\lambda_0-1)$ and $n_0>0$ such that
\[\varrho((1+\delta)u_n)\le \varrho((1+\delta)u)+\epsilon\le \varrho(u)+2\epsilon,\ \varrho\left(\frac{1+\delta}{\delta}v_n\right)\le\epsilon\]
for any $n\ge n_0$, by Lemma \ref{lem:continuity}. Then, the convexity of $\varrho$ yields
\begin{equation*}
	\varrho(u_n+v_n)\le\frac{1}{1+\delta}\varrho((1+\delta)u_n)+\frac{\delta}{1+\delta}\varrho\left(\frac{1+\delta}{\delta}v_n\right)
	\le\varrho(u)+3\epsilon.
\end{equation*}
Letting $n\to\infty$ and $\epsilon\downarrow 0$, we obtain the desired inequality. Moreover, this inequality yields
$\varrho(\lambda(u_n+v_n))\le\varrho(\lambda u)+1<\infty$ for fixed $\lambda\in(1,\lambda_0)$ and sufficiently large $n$, and hence $u_n+v_n\in\Mc{D}^i(\varrho)$. The last part of the assertion follows from the former one with $u_n=u$, and Lemma \ref{lem:lsc}.
\end{proof}
\begin{cor}\label{cor:Dint}
$\Mc{D}^i(\varrho)$ is the interior of $\Mc{D}(\varrho)$, with respect to $(\Mc{M}(\varrho),\|\cdot\|_{LN(\varrho)})$.
\end{cor}
\begin{proof}
Since Lemma \ref{lem:limeval} holds, $v$ is in the interior of $\Mc{D}(\varrho)$ whenever $v\in\Mc{D}^i(\varrho)$. Conversely, if $v\in \Mc{D}(\varrho)\setminus\Mc{D}^i(\varrho)$, then $v$ is on the boundary of $\Mc{D}(\varrho)$ because $\lambda v\to v\ (\lambda\downarrow 1)$ in $\|\cdot\|_{LN(\varrho)}$.
\end{proof}

\begin{lem}\label{lem:LNcontraction}
Let $X,V$ and $\varrho$ be as in Definition \ref{defi:condgsm}, and suppose that $\varrho$ satisfies \ref{cond:rwNC}. Let $\varphi:\Mb{R}\to\Mb{R}$ be a Lipschitz piecewise affine function with $\varphi(0)=0$. Then, for any $f\in\Mc{M}(\varrho)$, we have $\varphi\circ f\in\Mc{M}(\varrho)$ and
\[\|\varphi\circ f\|_{LN(\varrho)}\le\Lip(\varphi)\|f\|_{LN(\varrho)}.\]
\end{lem}
\begin{proof}
If $\Lip(\varphi)=0$, then $\varphi\equiv0$ and the assertion is immediate. Otherwise, the function $t\mapsto a\varphi(t/(a\Lip(\varphi)))$ satisfies the conditions in \ref{cond:rwNC}. Hence, $\varphi\circ f\in\Mc{M}(\varrho)$ and
\begin{equation}\label{eq:scaledNC}
\varrho(a\varphi\circ f)\le\varrho(a\Lip(\varphi)f)\qquad(a>0).
\end{equation}
The norm estimate follows from the definition of $\|\cdot\|_{LN(\varrho)}.$
\end{proof}

\section{Pushforward measures of pre-energy measures}\label{sec:pushforward}
In this section, we prepare some basic lemmas on $\Mc{E}$ and construct measures on $\Mb{R}$, corresponding to $g_*\Pi\langle f\rangle$ for $f\in\Mc{D}^i(\Mc{E})$ and $g\in\Mc{M}(\Mc{E})$, where $\Pi\langle f\rangle$ is the pre-energy measure mentioned in the Introduction. This ``$g_*\Pi\langle f\rangle$'' is effective even under assumptions weaker than those of Theorem \ref{thm:main1}. In particular, we derive a strong form of \ref{cond:rstrsubad} for strongly local energies from the ``pushforward'' measures.

Throughout this section, we assume the following.
\begin{asm}\label{asm:secpf}
$X$ is a (nonempty) topological space, $V$ is a linear subspace of $C_b(X)$ and $\Mc{E}$ is a gsm on $V$ satisfying \ref{cond:rwNC} and the following \ref{cond:rbvanish} (with $\Mc{E}=\varrho$).
\begin{description}
	\rcond[cond:rvanish]{cond:rbvanish}$0\vee f\wedge t\in\Mc{M}(\varrho)$ and $\lim_{t\downarrow0}\varrho(0\vee f\wedge t )=0$ for any $t>0$ and $f\in\Mc{M}(\varrho)\cap C_b(X)$.
\end{description}
\end{asm}
Obviously, \ref{cond:rvanish} and \ref{cond:rbvanish} are equivalent if $V\subset C_c(X)$. However, for future applications in which $X$ is not locally compact (for example, when $X$ is itself a space of functions or curves), we adopt \ref{cond:rbvanish} in this section.
\begin{lem}\label{lem:domcont}
\begin{enumerate}
	\item If $f,g\in\Mc{M}(\Mc{E})$, then $f\vee(g-a), f\wedge(g+a)\in\Mc{M}(\Mc{E})$ for any $a\ge0$.
	\item Let $\varphi:\Mb{R}\to\Mb{R}$ be a function as in \ref{cond:rwNC}. Then, $\varphi\circ f\in\Mc{D}^i(\Mc{E})$ whenever $f\in\Mc{D}^i(\Mc{E})$.
\end{enumerate}
\end{lem}
\begin{proof}
\begin{enumerate}
	\item The identities
	\begin{align}
		f\vee(g-a)&=\frac{f+g}{2}+\frac{|f-g+a|-a}{2},\label{eq:temp1}\\ f\wedge(g+a)&=\frac{f+g}{2}-\frac{|f-g-a|-a}{2},\label{eq:temp2}
	\end{align}
	ensure (1) by \ref{cond:rwNC}.
	\item By \eqref{eq:scaledNC}, $\Mc{E}(\lambda\varphi\circ f)\le\Mc{E}(\lambda f)$ for any $\lambda>0$, which yields (2).
\end{enumerate}
\end{proof}
\begin{lem}\label{lem:limeval2}
	Let $\{u_n\}_{n\ge 1},\{v_n\}_{n\ge 1}$ be sequences of functions in $V$ and let $u\in V$, satisfying the assumptions of Lemma \ref{lem:limeval} with $\varrho=\Mc{E}$. Then, for any sequence $\{a_n\}_{n\ge 1}$ of nonnegative numbers, $u_n\vee (v_n-a_n),u_n\wedge (v_n+a_n)\in\Mc{D}^i(\Mc{E})$ for sufficiently large $n$, and  $\limsup_{n\to\infty}\Mc{E}(u_n\vee(v_n-a_n))\le\Mc{E}(u)$ and $\limsup_{n\to\infty}\Mc{E}(u_n\wedge(v_n+a_n))\le\Mc{E}(u)$.
\end{lem}
\begin{proof}
	By \ref{cond:rconv}, \ref{cond:rwNC}, \eqref{eq:temp1} and \eqref{eq:temp2} we obtain
		\begin{align}
			\Mc{E}(u\vee (v-a))\le&\frac{1}{2}\left(\Mc{E}(u+v)+\Mc{E}(|u-v+a|-a) \right) \notag\\
		\le &\frac{1}{2}\left(\Mc{E}(u+v)+\Mc{E}(u-v) \right) \label{eq:evee}
		\shortintertext{and similarly,}
			\Mc{E}(u\wedge (v+a))\le&\frac{1}{2}\left(\Mc{E}(u+v)+\Mc{E}(u-v) \right)  \label{eq:ewedge}
			\end{align}
	for any $u,v\in\Mc{M}(\Mc{E})$ and $a\ge0$. Hence, the desired inequalities follow from Lemma \ref{lem:limeval}. In the same way as in the proof of Lemma \ref{lem:limeval}, $u_n\vee(v_n-a_n),u_n\wedge(v_n+a_n)\in\Mc{D}^i(\Mc{E})$ hold for sufficiently large $n$.
\end{proof}
\begin{lem}\label{lem:vanseq}
Let $f\in\Mc{M}(\Mc{E})$. If 
$\{\varphi_n\}_{n\ge 1}$ is a sequence of piecewise affine functions satisfying the following conditions, then $\varphi_n\circ f\to 0$ in $\|\cdot\|_{LN(\Mc{E})}$.
\begin{description}
	\Fcond{cond:fsup} There exists $M>0$ such that $\sup_{n\in\Mb{N}}|\varphi_n'(t)|<M\ \Mc{L}^1$-a.e.~$t$.
	\Fcond{cond:fbseq} There exist $N\in\Mb{N}$ and nested bounded sequences $\{[a_{i,n}, b_{i,n}]\}_{n\ge 1}$ of intervals for $1\le i\le N$ such that $\sup_{1\le i\le N}\lim_{n\to\infty} (b_{i,n}-a_{i,n})=0$ and $\varphi_n'=0$ on $\Mb{R}\setminus\cup_{i=1}^N [a_{i,n},b_{i,n}]$ for any $n$.
	\Fcond{cond:fzero} $\varphi_n(0)=0$ for any $n$.
\end{description}
\end{lem}
\begin{proof}
First note that, for any $h\in\Mc{M}(\Mc{E})$,
\[\lim_{t\downarrow0}\|\eta_{[0,t]}\circ h\|_{LN(\Mc{E})}=0\]
by \ref{cond:rbvanish}, Lemma \ref{lem:continuity} and the identity $\lambda\eta_{[0,t]}\circ h=\eta_{[0,\lambda t]}\circ(\lambda h)$ for $\lambda>0$.
We first assume that $N=1$ and $\varphi_n=\eta_{[a_n,b_n]}$: here we abbreviate the first index $i=1$. We also write $a_\infty=\lim_{n\to\infty}a_n=\lim_{n\to\infty}b_n$. If $a_\infty=0$, the triangle inequality yields
\begin{align*}
\|\eta_{[a_n,b_n]}\circ f\|_{LN(\Mc{E})}
\le{}&\|\eta_{[0,-a_n]}\circ(-f)\|_{LN(\Mc{E})}
+\|\eta_{[0,b_n]}\circ f\|_{LN(\Mc{E})}\\
\to{}&0\quad(n\to\infty).
\end{align*}
If $a_\infty>0$,
\begin{equation*}\eta_{[a_n,b_n]}\circ f
=\eta_{[0,a_\infty-a_n]}\circ f-\eta_{[0,a_\infty-a_n]}\circ g+\eta_{[0,b_n-a_\infty]}\circ h
\end{equation*}
whenever $2a_n\ge a_\infty$, where $g=f\wedge (a_\infty-f)$ and $h=\eta_{[a_\infty,b_1]}\circ f$: note that $g,h\in\Mc{M}(\Mc{E})$ by \ref{cond:rwNC}. Thus, the triangle inequality again gives $\|\eta_{[a_n,b_n]}\circ f\|_{LN(\Mc{E})}\to0$. The case $a_\infty<0$ is reduced to the previous case by $\eta_{[a_n,b_n]}\circ f=-\eta_{[-b_n,-a_n]}\circ(-f)$.

For the general case, our assumptions allow us to choose a piecewise affine function $\psi_n$ with $\psi_n(0)=0$ and $|\psi_n'|\le1\ \Mc{L}^1$-a.e. such that
\[\varphi_n=\psi_n\circ\left(M\sum_{i=1}^N\eta_{[a_{i,n},b_{i,n}]}\right).\]
Hence, by Lemma \ref{lem:LNcontraction} and the case $N=1$,
\begin{align*}
\|\varphi_n\circ f\|_{LN(\Mc{E})}
\le{}&M\sum_{i=1}^N\|\eta_{[a_{i,n},b_{i,n}]}\circ f\|_{LN(\Mc{E})}
\to0\quad(n\to\infty).
\end{align*}
\end{proof}

Hereafter, we write
\[f_n^{g,a}\coloneq (T_n\circ f)\wedge (S_n^a\circ g)\]
for any $f,g\in\Mc{M}(\Mc{E})$ and $a\in\Mb{R}$. Since $S_n^a(0)\ge 0$,
\[f_n^{g,a}=(T_n\circ f)\wedge\left((S_n^a-S_n^a(0))\circ g+S_n^a(0) \right)\in\Mc{M}(\Mc{E})\]
by Lemma \ref{lem:domcont}. 
\begin{prop}\label{prop:pushforward}
For $f,g\in\Mc{M}(\Mc{E})$ and $a\in\Mb{R}$, the limit $\lim_{n\to\infty}\Mc{E}(f_n^{g,a})$ exists whenever $f_n^{g,a}\in\Mc{D}^i(\Mc{E})$ for some $n$. Moreover, for any $f,g\in\Mc{M}(\Mc{E})$, the function $F^g_f$, defined by
\[F^g_f(a)=\begin{cases}
	\lim_{n\to\infty}\Mc{E}(f_n^{g,a})&(\text{if $f_n^{g,a}\in\Mc{D}^i(\Mc{E})$ for some $n$})\\
		\infty&(\text{otherwise})
\end{cases} \] is a nondecreasing right-continuous function with $\lim_{a\to-\infty}F^g_f(a)=0$. In particular, if $T_n\circ f\in\Mc{D}^i(\Mc{E})$ for some $n$ then $F_f^g$ is a distribution function of a finite Radon measure on $\Mb{R}$.
\end{prop}
\begin{proof}
	For any $n,m\in\Mb{N}$ with $f_m^{g,a}\in \Mc{D}^i(\Mc{E})$ and $m<n$, we have
	\[f_n^{g,a}=(T_n\circ f_m^{g,a})\wedge (S_n^a\circ g) \]
	and
	\[\limsup_{n\to\infty}\Mc{E}(f_n^{g,a})\le\Mc{E}(f_m^{g,a})<\infty\]
	by Lemma \ref{lem:limeval2} (because $\Mc{E}(\alpha T_n\circ f)\le\Mc{E}(\alpha f)$ for $\alpha>0$ by \eqref{eq:scaledNC}, and $\|f\|_\infty<\infty$). This ensures existence of $\lim_{n\to\infty}\Mc{E}(f_n^{g,a})$.
	
	Here, fix any $\epsilon>0$, $a\in\Mb{R}$, $m\in\Mb{N}$ such that $\Mc{E}(f_m^{g,a})\le F_f^g(a)+\epsilon<\infty$. Then, for any $b<a+2^{-m}$,
	\[f_n^{g,b}=(T_n\circ f_m^{g,a})\wedge (S_n^b\circ g)\]
	whenever $n\ge m$ and $b+2^{1-n}<a+2^{-m}$. Hence, by Lemma \ref{lem:limeval2}, we obtain
	\[F_f^g(b)=\lim_{n\to\infty}\Mc{E}(f_n^{g,b})\le \Mc{E}(f_m^{g,a})\le F_f^g(a)+\epsilon,\]
	which leads to the right continuity and monotonicity of $F_f^g$. Moreover,
	\[F^g_f(a)=\lim_{n\to\infty}\Mc{E}(T_n\circ f)\le\Mc{E}(f)\]
	for any $a\ge \|g\|_{\infty}$ with $F_f^g(a)<\infty$, and
	\[F^g_f(a)=\lim_{n\to\infty}\Mc{E}(0)=0\]
	for any $a<-\|g\|_{\infty}$.
	Therefore, we conclude that $F_f^g$ is a distribution function of a finite Radon measure on $\Mb{R}$, whenever $T_n\circ f\in\Mc{D}^i(\Mc{E})$ for some $n$.
\end{proof}
In the remainder of this section, we prove some properties of $F_f^g$ for applications later in this paper.

\begin{lem}\label{lem:anotherwave}
 $F_f^g(a)=\lim_{n\to\infty}\Mc{E}((T_{n,\alpha}\circ f)\wedge (\alpha S_n^{a}\circ g))<\infty$ for any $\alpha>0$, if and only if $F_f^g(a)<\infty$.
\end{lem}
\begin{proof}
	By Lemmas \ref{lem:limeval2} and \ref{lem:vanseq}, it suffices to show 
	\[\lim_{n\to\infty}\Mc{E}((T_{n,\alpha}\circ f)\wedge (\alpha S_n^{a}\circ g))\le \Mc{E}((T_{m,\beta}\circ f)\wedge (\beta S_m^{a}\circ g))\]
	for any $f,g\in\Mc{M}(\Mc{E})$, $\alpha,\beta>0$, $a\in\Mb{R}$ and $m\in\Mb{N}$ with $(T_{m,\beta}\circ f)\wedge(\beta S_m^a\circ g)\in\Mc{D}^i(\Mc{E})$.
	Let 
	\begin{align*}
		I_{n}&=
		\left(\Cup_{i=1}^{2^{m-1}N}\left( \frac{2i-1}{2^m}\beta, \frac{2i-1}{2^m}\beta+\frac{\alpha}{2^n}  \right)\right)\Cup
		\left(\Cup_{i=1}^{2^{m-1}N}\left( \frac{2i}{2^m}\beta-\frac{\alpha}{2^n},\frac{2i}{2^m}\beta\right)\right),\\
		J_{n}&=
		\left(\Cup_{i=1}^{2^{m-1}N}\left( \frac{2i-1}{2^m}\beta,\left\lceil \frac{2i-1}{2^m}\beta\right\rceil_{\frac{\alpha}{2^n }} \right)\right)\Cup
		\left(\Cup_{i=1}^{2^{m-1}N}\left(\left\lfloor \frac{2i}{2^m}\beta\right\rfloor_{\frac{\alpha}{2^n}},\frac{2i}{2^m}\beta\right)\right),\\
		u_{n}(t)&=2\int_0^{t}\bigl(\indi_{J_{n}}(s)-\indi_{-J_{n}}(s)\bigr)\,\der s,
	\end{align*}
	where
	\[
	\lfloor x\rfloor_c=\max\{c i\mid i\in\Mb{Z},\ c i\le x \},
	\qquad
	\lceil x \rceil_c=\min\{c i\mid i\in\Mb{Z},\ c i\ge x\},
	\]
	and $N=\left\lceil \|f\|_{\infty}/\beta \right\rceil_1$.  
	
	Then, for any $n$ with $2^{-n}\alpha\le 2^{-m}\beta$, we have $T_{n,\alpha}(t)=T_{n,\alpha}(T_{m,\beta}(t)+u_n(t))$ for $|t|\le N\beta$, $u_n(0)=0$, and $u_n'(t)=0$ for every $t\in\Mb{R}\setminus\overline{I_n\cup(-I_n)}$.
	Then, for sufficiently small $\delta>0$,
	\begin{align*}
		&\limsup_{n\to\infty}\Mc{E}\left((1+\delta) T_{n,\alpha}\circ \left(\bigl((T_{m,\beta}\circ f)\wedge(\beta S_m^a\circ g) \bigr)+(u_n\circ f) \right)\right)\\
		\le&\limsup_{n\to\infty}\Mc{E}\left((1+\delta)\left((T_{m,\beta}\circ f)\wedge(\beta S_m^a\circ g ) \right) +(1+\delta)(u_n\circ f)  \right)\\
		\le& \Mc{E}((1+\delta)\left((T_{m,\beta}\circ f)\wedge(\beta S_m^a\circ g ) \right) )
	\end{align*}
	by \eqref{eq:scaledNC}, Lemmas \ref{lem:limeval} and \ref{lem:vanseq}. Using Lemmas \ref{lem:continuity} and \ref{lem:limeval2}, we obtain
	\begin{align*}
		\lim_{n\to\infty}\Mc{E}\left((T_{n,\alpha}\circ f)\wedge (\alpha S_n^{a}\circ g)\right)\le\Mc{E}((T_{m,\beta}\circ f)\wedge (\beta S_m^a \circ g))
	\end{align*}
	because
	\begin{align*}
		&(T_{n,\alpha}\circ f)\wedge (\alpha S_n^{a}\circ g)\\
		=&\left(T_{n,\alpha}\circ \left(\left((T_{m,\beta}\circ f)\wedge(\beta S_m^a\circ g) \right)+(u_n\circ f) \right)\right)\wedge \left(\alpha S_n^a\circ g \right)
	\end{align*}
	whenever $n>m$ and $2^{-n}\alpha\le\beta2^{-m}$.
\end{proof}
%
%
\begin{lem}\label{lem:localvanishpre}
$F_f^{|f-a|-|a|}(-|a|)=0$ for any $f\in\Mc{M}(\Mc{E})$ and $a\in\Mb{R}$.
\end{lem}
\begin{proof}
	Let 
	\[\varphi_n(t)=T_n(t)\wedge\eta_{[0,2^{-n}]}(2^{1-n}-|t-a| )  .\]
	Then, $f_n^{|f-a|-|a|,-|a|}=\varphi_n\circ f,\ \varphi_n(0)=0$, $|\varphi_n'(t)|\le 1\ \Mc{L}^1$-a.e.~and $\varphi_n'(t)=0$ on $\Mb{R}\setminus[a-2^{1-n},a+2^{1-n}]$. Hence, $F_f^{|f-a|-|a|}(-|a|)=\lim_{n\to\infty}\Mc{E}(\varphi_n\circ f)=0$ by Corollary \ref{cor:Dint} and Lemma \ref{lem:vanseq}.
\end{proof}

\subsection{Strongly local case}
In this subsection, we prove a strong version of \ref{cond:rstrsubad} and another property under a variant of \ref{cond:rstrloc}.
\begin{defi}\label{defi:gsmmisc}
 Let $X,V,\varrho$ be as in Definition \ref{defi:condgsm}. We define the following condition, which we apply with $\Mc{E}=\varrho$.
 \begin{description}
	\rcond[cond:rstrloc]{cond:rstrlocbd} Assume that $f,g\in \Mc{M}(\varrho)\cap C_b(X)$, and that $\supp(f)$ and $\supp(g+a)$ are separated by open sets for some $a\in\Mb{R}$. Then, $\varrho(f+g)=\varrho(f)+\varrho(g)$.
 \end{description}
\end{defi}
If $X$ is Hausdorff and $V\subset C_c(X)$, then \ref{cond:rstrloc} and \ref{cond:rstrlocbd} are equivalent; we use \ref{cond:rstrlocbd} for the same reason as \ref{cond:rbvanish}.

For further consequences, see Appendix \ref{ap:further}.

\begin{lem}\label{lem:MisD}
If $\Mc{E}$ satisfies \ref{cond:rstrlocbd}, then $\Mc{M}(\Mc{E})=\Mc{D}(\Mc{E})=\Mc{D}^i(\Mc{E})$.
\end{lem}
\begin{proof}
Let $f\in\Mc{M}(\Mc{E})$. By Lemma \ref{lem:vanseq}, for any $b\in\Mb{R}$, there exists $\epsilon>0$ such that $\eta_{[b-\epsilon,b+\epsilon]}\circ f\in\Mc{D}(\Mc{E})$. Hence, by a compactness argument, there exist $a_0<\dots<a_n$, $b_1<\dots<b_n$ and $\epsilon_i>0$ for $1\le i\le n$ such that $[-\|f\|_{\infty}, \|f\|_{\infty}]\subset (a_0,a_n)$, $[a_{i-1},a_i]\subset(b_i-\epsilon_i,b_i+\epsilon_i)$ and $\eta_{[b_i-\epsilon_i,b_i+\epsilon_i]}\circ f\in\Mc{D}(\Mc{E})$. Then, by Lemmas \ref{lem:lsc} and \ref{lem:vanseq} we obtain
\begin{multline*}
\Mc{E}(f)\le\liminf_{\epsilon\downarrow 0}\Mc{E}\left(\sum_{i=1}^n \eta_{[a_{i-1}, a_i-\epsilon]}\circ f\right)\\
\overset{\ref{cond:rstrlocbd}}{=}\liminf_{\epsilon\downarrow 0}\sum_{i=1}^n \Mc{E}(\eta_{[a_{i-1}, a_i-\epsilon]}\circ f)\overset{\ref{cond:rwNC}}{\le} \sum_{i=1}^n \Mc{E}(\eta_{[b_{i}-\epsilon_i, b_i+\epsilon_i]}\circ f)<\infty.
\end{multline*}
This yields $\Mc{M}(\Mc{E})=\Mc{D}(\Mc{E})$, and hence it also coincides with $\Mc{D}^i(\Mc{E})$.
\end{proof}	

\begin{lem}\label{lem:locfold}
	Let $\varphi:\Mb{R}\to\Mb{R}$ be a piecewise affine function with $|\varphi'(t)|=1\ \Mc{L}^1$-a.e. and $\varphi(0)=0$. If $\Mc{E}$ satisfies \ref{cond:rstrlocbd}, then $\Mc{E}(\varphi\circ f)=\Mc{E}(f)$ for any $f\in\Mc{M}(\Mc{E})$.
\end{lem}
\begin{proof}
	Fix any $f\in\Mc{M}(\Mc{E})$. Since $\|f\|_{\infty}<\infty$, there exist $n\in\Mb{N}$, $\{a_i\}_{i=0}^n\in\Mb{R}^{n+1}$ and $\{\epsilon_i\}_{i=1}^n\in\{1,-1\}^n$ such that $a_0<a_1<\dots <a_n$, $a_0<-\|f\|_{\infty}$, $a_n>\|f\|_{\infty}$ and
	\[\varphi\circ f=\sum_{i=1}^n \epsilon_i\eta_{[a_{i-1},a_i]}\circ f.\]
	
	Since \ref{cond:rwNC} holds, it is sufficient to prove $\Mc{E}(\varphi\circ f)\ge\Mc{E}(f)$. By Lemma \ref{lem:MisD}, $\varphi\circ f\in \Mc{D}^i(\Mc{E})$. Then, by Lemmas \ref{lem:lsc}, \ref{lem:limeval} and \ref{lem:vanseq}, we obtain
	\begin{multline*}
		\Mc{E}(\varphi\circ f)=\lim_{t\downarrow 0}\Mc{E}\left(\sum_{i=1}^{n} \epsilon_i\eta_{[a_{i-1}+t,a_i]}\circ f \right)\\
		\overset{\ref{cond:rstrlocbd}}{=}\lim_{t\downarrow0}\sum_{i=1}^n\Mc{E}\left( \epsilon_i\eta_{[a_{i-1}+t,a_i]}\circ f \right)
		\overset{\ref{cond:rsym}}{=}\lim_{t\downarrow0}\sum_{i=1}^n\Mc{E}\left( \eta_{[a_{i-1}+t,a_i]}\circ f \right)\\
		\overset{\ref{cond:rstrlocbd}}{=}\lim_{t\downarrow0}\Mc{E}\left( \sum_{i=1}^n\eta_{[a_{i-1}+t,a_i]}\circ f \right)
		\ge\Mc{E}(f).
	\end{multline*}
\end{proof}
\begin{lem}\label{lem:locpushforward}
	Assume that $\Mc{E}$ satisfies \ref{cond:rstrlocbd}. Then, $F_f^g(a)+\lim_{t\downarrow0}F_f^{-g}(-a-t)=\Mc{E}(f)$ for any $f,g\in\Mc{M}(\Mc{E})$ and $a\in\Mb{R}$.
\end{lem}	
	\begin{proof}
		We first prove that 
		\begin{equation}\label{eq:Fevalbasic}
F_f^g(a)+F_f^{|g-c|-|c|}(b-c-|c|)\le F^g_f(b)
		\end{equation}
		for $b>a$ and $c\in(\frac{a+b}{2},b)$. In fact, $\supp(f_n^{g,a})\cap \supp(f_n^{|g-c|-|c|,b-c-|c|})=\emptyset$ and
		\begin{gather*}
		f_n^{g,a}+f_n^{|g-c|-|c|,b-c-|c|}=\left(T_n\circ f_m^{g,b}\right)\wedge (\varphi_n\circ g)
		\shortintertext{for any $n,m$ with $n\ge m$ and $2c-b-a>2^{2-n}$, where}
		\varphi_n(t)=\eta_{[0,2^{-n}]}\left((2^{1-n}+b-c-|t-c|)\vee (a+2^{1-n}-t) \right).
		\end{gather*}
		Hence, since $f_m^{g,b}\in\Mc{D}^i(\Mc{E})$ by Lemma \ref{lem:MisD},
		\[F_f^g(a)+F_f^{|g-c|-|c|}(b-c-|c|)=\lim_{n\to\infty}\Mc{E}(f_n^{g,a}+f_n^{|g-c|-|c|,b-c-|c|})\le\Mc{E}(f_m^{g,b}) \]
		by \ref{cond:rstrlocbd} and Lemmas \ref{lem:limeval2} and \ref{lem:vanseq}. Letting $m\to\infty$, we obtain \eqref{eq:Fevalbasic}.
		
		Let $b>\|g\|_{\infty}$. Then $F_f^g(b)=\lim_{n\to\infty}\Mc{E}(T_n\circ f)=\Mc{E}(f)$ by Lemma \ref{lem:locfold} and $F_f^{|g-c|-|c|}(b-c-|c|)=F_f^{-g}(b-2c)$. Together with \eqref{eq:Fevalbasic}, this proves $F_f^g(a)+\lim_{t\downarrow0}F_f^{-g}(-a-t)\le\Mc{E}(f)$.
		
		We next prove the reverse inequality. We first note that
\begin{align}
	\Mc{E}(\lambda (u \vee v))\overset{\eqref{eq:evee}}{\le}& \frac{1}{2}\left(\Mc{E}(\lambda (u+v))+\Mc{E}(\lambda(u-v)) \right)\notag\\
	\overset{\ref{cond:rconv}}\le& \frac{1}{2}\left(2\lambda\Mc{E}(u)+(1-\lambda)\Mc{E}\left(\frac{-\lambda v}{\lambda-1}\right)+(1-\lambda)\Mc{E}\left(\frac{\lambda v}{\lambda-1}\right)\right)\notag\\
	\overset{\ref{cond:rsym}}{=}&\lambda\Mc{E}(u)+(1-\lambda)\Mc{E}\left(\frac{\lambda}{1-\lambda}v\right)\label{eq:convestm}
\end{align}
for any $u,v\in\Mc{M}(\Mc{E})$ and $\lambda\in(0,1)$.
Fix $f,g\in\Mc{M}(\Mc{E})$, $a\in\Mb{R}$ and $\epsilon\in(0,1)$. By \eqref{eq:Fevalbasic} and the right continuity of $F^g_{\frac{1-\epsilon}{\epsilon}f}$,
\begin{multline*}
0\le F^{|g-a-2t|-|a+2t|}_{\frac{1-\epsilon}{\epsilon}f}(t-|a+2t|)\\
\le F^g_{\frac{1-\epsilon}{\epsilon}f}(a+3t)-F^g_{\frac{1-\epsilon}{\epsilon}f}(a)\to0\quad(t\downarrow0).
\end{multline*}
Choose $t_\epsilon\in(0,\epsilon)$ such that
\[
F^{|g-a-2t_\epsilon|-|a+2t_\epsilon|}_{\frac{1-\epsilon}{\epsilon}f}(t_\epsilon-|a+2t_\epsilon|)<\epsilon/2,
\]
and set
\[
h_{\epsilon,n}\coloneq \frac{1-\epsilon}{\epsilon}f_n^{|g-a-2t_\epsilon|-|a+2t_\epsilon|,\,t_\epsilon-|a+2t_\epsilon|}.
\]
By Lemma \ref{lem:anotherwave}, applied to $\frac{1-\epsilon}{\epsilon}f$ with $\alpha=\frac{1-\epsilon}{\epsilon}$, $\Mc{E}(h_{\epsilon,n})\le\epsilon$ for all sufficiently large $n$. If also $2^{2-n}<2t_\epsilon$, the supports of $f_n^{g,a+t_\epsilon}$ and $f_n^{-g,-a-3t_\epsilon}$ are separated by open sets, and
\[
T_n\circ f=(f_n^{g,a+t_\epsilon}+f_n^{-g,-a-3t_\epsilon})\vee\left(\frac{\epsilon}{1-\epsilon}h_{\epsilon,n}\right).
\]
By Lemma \ref{lem:locfold}, \eqref{eq:convestm} and \ref{cond:rstrlocbd},
\begin{align*}
\Mc{E}((1-\epsilon)f)
&=\Mc{E}((1-\epsilon)T_n\circ f)\\
&\le(1-\epsilon)\Mc{E}(f_n^{g,a+t_\epsilon}+f_n^{-g,-a-3t_\epsilon})+\epsilon\Mc{E}(h_{\epsilon,n})\\
&\le\Mc{E}(f_n^{g,a+t_\epsilon})+\Mc{E}(f_n^{-g,-a-3t_\epsilon})+\epsilon^2.
\end{align*}
Letting $n\to\infty$ first, we obtain
\[
\Mc{E}((1-\epsilon)f)\le F_f^g(a+t_\epsilon)+F_f^{-g}(-a-3t_\epsilon)+\epsilon^2.
\]
Since $t_\epsilon\to0$ as $\epsilon\downarrow0$, \ref{cond:rlc}, right continuity of $F_f^g$ and monotonicity of $F_f^{-g}$ yield
\[
\Mc{E}(f)\le F_f^g(a)+\lim_{t\downarrow0}F_f^{-g}(-a-t).
\]
\end{proof}

\begin{prop}\label{prop:localmeansnldc} Assume that $\Mc{E}$ satisfies \ref{cond:rstrlocbd}.
	\begin{enumerate}
		\item $\Mc{E}(f\wedge g)+\Mc{E}(f\vee g)=\Mc{E}(f)+\Mc{E}(g)$ for any $f,g\in\Mc{M}(\Mc{E})$.
		\item $\Mc{E}((f-a)\vee g\wedge (f+a))+\Mc{E}((g-a)\vee f\wedge (g+a))=\Mc{E}(f)+\Mc{E}(g)$ for any $f,g\in\Mc{M}(\Mc{E})$ and $a\ge0$.
	\end{enumerate}
\end{prop}
\begin{proof}
We first prove $\Mc{E}(f\wedge g)+\Mc{E}(f\vee g)\le\Mc{E}(f)+\Mc{E}(g)$ for $f,g\in\Mc{M}(\Mc{E})$. Note that, if $u,v,w\in\Mc{M}(\Mc{E})$ and $a,b,c\in\Mb{R}$ satisfy $a<b$, $c>0$ and $(u+ck)|_{\{w\le b\}}=v|_{\{w\le b\}}$ for some $k\in\Mb{Z}$, then $T_{n,c}\circ u=T_{n,c}\circ v$ on $\{w\le b\}$ and
\begin{multline}
F_u^w(a)=\lim_{n\to\infty}\Mc{E}\left((T_{n,c}\circ u)\wedge (cS_n^a\circ w )\right)\\
=\lim_{n\to\infty}\Mc{E}\left((T_{n,c}\circ v)\wedge (cS_n^a\circ w )\right)=F_v^w(a) \label{eq:fequiv}
\end{multline}
by Lemma \ref{lem:anotherwave}.
Since $F^{f-g}_{f\vee g}$ and $F^{f-g}_{f\wedge g}$ are monotone, we can choose $a_n>0$ so that $\lim_{n\to\infty}a_n=0$ and $F^{f-g}_{f\vee g}, F^{f-g}_{f\wedge g}$ are continuous at $a_n$ for any $n$. Then,
\[f-((f-a_n)\vee g\wedge f)=(g\vee f\wedge(g+a_n))-g=\eta_{[0,a_n]}(f-g)\]
converges to $0$ in $\|\cdot\|_{LN(\Mc{E})}$ by Lemma \ref{lem:vanseq}. Hence,
\begin{align*}
	&\Mc{E}(f\wedge g)+\Mc{E}(f\vee g)\\
=&\lim_{t\downarrow0}\left(F_{f\wedge g}^{f-g}(a_n-t)+F_{f\wedge g}^{g-f}(-a_n-t)+F_{f\vee g}^{f-g}(a_n-t)+F_{f\vee g}^{g-f}(-a_n-t)\right)\\
=&\lim_{t\downarrow0}\begin{multlined}[t][.9\linewidth]
\Big(F_{(f-a_n)\vee g\wedge f}^{f-g}(a_n-t)+F_{g\vee f\wedge(g+a_n)}^{g-f}(-a_n-t)\\
+F_{g\vee f\wedge(g+a_n)}^{f-g}(a_n-t)+F_{(f-a_n)\vee g\wedge f}^{g-f}(-a_n-t)\Big)
\end{multlined}\\
\le &\Mc{E}((f-a_n)\vee g\wedge f)+\Mc{E}(g\vee f\wedge(g+a_n))
\end{align*}
by Lemma \ref{lem:locpushforward} and \eqref{eq:fequiv}. Therefore, we obtain
\begin{align*}
\Mc{E}(f\vee g)+\Mc{E}(f\wedge g)\le& \lim_{n\to\infty}\left(\Mc{E}((f-a_n)\vee g\wedge f)+\Mc{E}(g\vee f\wedge(g+a_n)) \right)\\
=&\Mc{E}(f)+\Mc{E}(g) 
\end{align*}
by Lemmas \ref{lem:limeval} and \ref{lem:MisD}.

Conversely, if $F_f^{f-g}$ and $F_g^{f-g}$ are continuous at $b_n>0$ for any $n\ge 0$, and $\lim_{n\to\infty}b_n=0$, we obtain
\[\Mc{E}(f)+\Mc{E}(g)\le\Mc{E}(f\wedge (g+b_n))+\Mc{E}((f-b_n)\vee g )\to \Mc{E}(f\wedge g)+\Mc{E}(f\vee g)\ (n\to\infty).\]

For (2), we first note that
\begin{equation}
	F_u^{-|v|}(t)=F_u^v(t)+F_u^{-v}(t)\label{eq:Fsplit}
\end{equation}
for any $u,v\in\Mc{M}(\Mc{E})$ and $t<0$, by definition of $F_u^v$ and \ref{cond:rstrlocbd}. Similarly to (1), using equations  \eqref{eq:fequiv} and \eqref{eq:Fsplit}, and Lemmas \ref{lem:limeval} and \ref{lem:vanseq}, we obtain
\begin{align*}
	&\Mc{E}((f-a)\vee g\wedge (f+a) )+\Mc{E}((g-a)\vee f\wedge (g+a))\\
\le&\begin{multlined}[t][.9\linewidth]
	\Mc{E}\left((g-b_n)\vee\left((f-a)\vee g\wedge (f+a) \right)\wedge (g+b_n) \right)\\
	+\Mc{E}\left((f-b_n)\vee\left((g-a)\vee f\wedge (g+a) \right)\wedge (f+b_n) \right)
\end{multlined}	\\
\to&\Mc{E}(g)+\Mc{E}(f)\ (n\to\infty)
\end{align*}
and
\begin{align*}
&\Mc{E}(f)+\Mc{E}(g)\\
\le&\Mc{E}((f-a-b_n)\vee g\wedge (f+a+b_n) )+\Mc{E}((g-a-b_n)\vee f\wedge (g+a+b_n))\\
\to&\Mc{E}((f-a)\vee g\wedge (f+a) )+\Mc{E}((g-a)\vee f\wedge (g+a))\ (n\to\infty)
\end{align*}
for some $\{b_n\}_{n\ge 1}$.
\end{proof}
\begin{rem}
 Proposition \ref{prop:localmeansnldc} was implicitly given in \cite{Sas2026}, in the case of strongly local, regular $p$-energy forms satisfying $p$-Clarkson inequalities; this is an easy consequence of the existence of $p$-energy measures (see Section \ref{sec:ex} or \cite{Sas2026} for details of these notions). A direct proof without $p$-energy measures was given by Kajino~\cite{Kaj0000} in the case of strongly local $p$-energy forms. The proof above is based on Kajino's.
\end{rem}

\section{Construction and uniqueness of energy measures}\label{sec:construction}
In this section, we impose the following Assumption \ref{asm:secconst} and prove Theorem \ref{thm:main1}. The former half of this section is devoted to constructing the pre-energy measure $\Pi\langle f\rangle$, while the latter part establishes the remaining statements needed for Theorem \ref{thm:main1}, such as the modification of pre-energy measures and ingredients for the uniqueness of energy measures.
\begin{asm}\label{asm:secconst}
	$X$ is a (nonempty) locally compact metrizable space and $\Mc{E}$ is a gsm on $V\subset C_c(X)$ satisfying \ref{cond:rwNC}--\ref{cond:prereg}. We fix a distance $d$ inducing the topology of $X$.
\end{asm}
Under this assumption, Lemma \ref{lem:domcont}~(1) is partially strengthened as follows.
\begin{lem}\label{lem:domsubad}
$f\vee g, f\wedge g\in\Mc{D}^i(\Mc{E})$ whenever $f,g\in\Mc{D}^i(\Mc{E})$.
\end{lem}
\begin{proof}
It easily follows from \ref{cond:rstrsubad}.
\end{proof}

\begin{lem}\label{lem:slider}
	Let $\Mc{G}$ be a dense subset of $(C_c(X),\|\cdot\|_{\infty})$. Then for any compact set $E$ and open set $U$ with $E\subset U$, there exists $g\in\Mc{G}$ with $E\subset \{g\le-2\}$ and $\{g<-1\}\Subset U$.
\end{lem}
\begin{proof}
Let $V$ be open with $E\subset V\Subset U$. By Urysohn's lemma, there exists $f:X\to[0,1]$ such that $f\in C(X)$, $f|_E\equiv 1$ and $\supp(f)\subset \ol{V}$. Since $\Mc{G}$ is dense, there exists $g\in\Mc{G}$ satisfying $\|-3f-g\|_{\infty}\le 1$. Then $g$ satisfies the desired properties.
\end{proof}
\begin{cor}\label{cor:bumpfcn}
	For any compact set $E$ and open set $U$ with $E\subset U$, there exists $h:X\to[0,1]$ such that $h\in\Mc{M}(\Mc{E})$, $\supp(h)\subset U$ and $h|_E\equiv 1$.
\end{cor}
\begin{proof}
Let $g$ be as in Lemma \ref{lem:slider} with $\Mc{G}=\Mc{M}(\Mc{E})$. Then, $h=\eta_{[1,2]}\circ(-g)$ has the desired properties, because Lemma \ref{lem:domcont}~(1) holds.
\end{proof}

\begin{lem}\label{lem:Fsubad}
Let $N\in\Mb{N}$, $f,g,h_i\in\Mc{M}(\Mc{E})$ and $a,b_i\in\Mb{R}$ for any $1\le i\le N$. Then, $F_f^g(a)\le \sum_{i=1}^N F_f^{h_i}(b_i)$ whenever
$\{g\le a \}\subset \cup_{i=1}^N \{h_i\le b_i\}$.
\end{lem}
\begin{proof}
By definition of $F_f^{h_i}$ and Lemmas \ref{lem:vanseq} and \ref{lem:domsubad}, we may assume $f_{m}^{h_i,b_i}\in\Mc{D}^i(\Mc{E})$ for any $i$, for sufficiently large $m>0$, without loss of generality. In the same way as in the proof of \cite[Lemma 5.4]{Sas2026}, we can show that there exists $n_0$ such that $\{g\le a+2^{-n_0}\}\subset \cup_{i=1}^N \{h_i\le b_i+2^{-m}\}$. Then,
\[f_n^{g,a}=\left( T_n\circ \left(\max_{1\le i\le N}f_m^{h_i,b_i} \right)\right)\wedge (S_n^a\circ g) \]
for any $n>n_0$. Therefore,
\begin{align*}
	F_f^{g}(a)\le \Mc{E}\left(\max_{1\le i\le N}f_m^{h_i,b_i}\right)\le \sum_{i=1}^N \Mc{E}(f_m^{h_i,b_i})\to \sum_{i=1}^N F_f^{h_i}(b_i)\ (m\to\infty)
\end{align*}
by Lemmas \ref{lem:limeval2}, \ref{lem:vanseq}, and \ref{lem:domsubad}, and \ref{cond:rstrsubad}.
\end{proof}
\begin{defi}\label{defi:premeas}
For $f\in\Mc{M}(\Mc{E})$, we define
\[\Pi_{\Mc{E}}\langle f\rangle (U)=\sup\{F_f^g(a)\mid g\in\Mc{M}(\Mc{E}), a\in\Mb{R}\text{ with }\{g\le a\}\subset U  \}\]
 for any open set $U\subset X$ and 
\[\Pi_{\Mc{E}}\langle f\rangle (A)=\inf\{\Pi \langle f\rangle (U)\mid A\subset U\text{: open} \}\]
for any $A\subset X$. We abbreviate the subscript $\Mc{E}$ if no confusion may occur.
\end{defi}
\begin{lem}\label{lem:zero-zero}
	$\Pi\langle f\rangle (\{g\le a\})=F_f^g(a)$ for any $f,g\in\Mc{M}(\Mc{E})$ and $a\in\Mb{R}$. In particular, $\Pi\langle f\rangle (\{f= a\})=0$ for any $f\in\Mc{M}(\Mc{E})$ and $a\in\Mb{R}$.
\end{lem}
\begin{proof}
	For any open $U\supset \{g\le a\}$, we have $F_f^g(a)\le\Pi\langle f\rangle(U)$ and hence $\Pi\langle f\rangle (\{g\le a\})\ge F_f^g(a)$. Conversely, for any $a'>a$, $\Pi\langle f\rangle (\{g<a'\})\le F_f^g(a')$ by Lemma \ref{lem:Fsubad} for $N=1$. Thus,
	\[\Pi\langle f\rangle(\{g\le a\})\le \Pi\langle f\rangle (\{g<a'\})\le F_f^g(a')\to F_f^g(a)\quad (a'\to a) \]
	by the right continuity of $F_f^g$. Moreover,
	\[\Pi\langle f\rangle (\{f=a\})=\Pi\langle f\rangle (\{|f-a|-|a|\le -|a|\})=F_f^{|f-a|-|a|}(-|a|)=0\]
	by Lemma \ref{lem:localvanishpre}.
\end{proof}
\begin{thm}\label{thm:Radon}
For each $f\in\Mc{M}(\Mc{E})$, $\Pi\langle f\rangle|_{\Mc{B}(X)}$ is a finite Radon measure.
\end{thm}
\begin{proof}
We first prove that 
\begin{equation}
\Pi\langle f\rangle (U)=\sup\{F_f^g(a)\mid g\in\Mc{M}(\Mc{E}), a<0\text{ with }\{g\le a\}\subset U  \}\label{eq:emminus}
\end{equation}
for any open set $U$. In fact, for any $u,v\in\Mc{M}(\Mc{E})$ and $a,b\in\Mb{R}$ we have
\begin{multline*}
\{u\le a, v\le b\}=\{(u-a)\vee (v-b)\le 0\}\\=\left\{\left(u+(0\wedge(b-a))\right)\vee \left(v+(0\wedge(a-b))\right)\le a\wedge b\right\},
\end{multline*}
and hence, for arbitrarily fixed $g\in\Mc{M}(\Mc{E})$ and $a\in\Mb{R}$ with $\{g\le a\}\subset U$, we obtain $\{h_t\le a\wedge -t\}\subset U$ and
\[ F_f^g(a)\le F_f^{|f|}(t)+F_f^{h_t}(a\wedge -t)\]
for any $t>0$ by Lemma \ref{lem:Fsubad}, where $h_t=(g+(0\wedge(-t-a)))\vee (-|f|+(0\wedge (a+t)))$.
By the right continuity of $F_f^{|f|}$ and Lemma \ref{lem:localvanishpre} with $a=0$, we obtain $F_f^g(a)\le\limsup_{t\downarrow 0}F_f^{h_t}(a\wedge -t)$ and also obtain \eqref{eq:emminus}.

Since \eqref{eq:emminus} and Lemma \ref{lem:Fsubad} hold, the subadditivity of $\Pi\langle f\rangle$ follows in the same way as in the proof of \cite[Theorem 5.6]{Sas2026}. Thus, $\Pi$ is an outer measure.

To show $\Pi\langle f\rangle$-measurability of Borel sets, it is sufficient to prove
\begin{equation}\label{eq:mom}
\Pi\langle f\rangle(A)+\Pi\langle f\rangle(B)\le \Pi\langle f\rangle(A\cup B)\ (A,B\subset X)
\end{equation}
whenever $d(A,B)>0$. (This sufficient condition is well-known: see \cite[Section 2.3.2~(9)]{Fed1969} for example.)
Fix any open sets $U,V$, $g,h\in\Mc{M}(\Mc{E})$ and $a,b<0$ such that $U\cap V=\emptyset,\ \{g\le a\}\subset U$ and $\{h\le b\}\subset V$. Then, by the argument of decreasing compact sets, there exist $a'\in(a,0)$ and $b'\in(b,0)$ such that $\{g\le a'\}\subset U$ and $\{h\le b'\}\subset V$. We also let $u,v'\in\Mc{M}(\Mc{E})$ such that $u,v':X\to[0,1]$, $u|_{\{g\le a' \}}\equiv 1$, $u|_{X\setminus U}\equiv 0$, $v'|_{\{h\le b'\}}\equiv 1$ and $\supp(v')\subset V$. Let
\[w=(g-a+(a\vee b))\wedge (h-b+(a\vee b)), \quad c=a\vee b.\]
Then, $\{w\le c+ t\}= \{g\le a+t\}\cup\{h\le b+t\}$ for any $t\in\Mb{R}$. In particular, $S_n^a\circ g+S_n^b\circ h=S_n^c\circ w$ and $f_n^{g,a}+f_n^{h,b}=f_n^{w,c}$ whenever $2^{1-n}\le (a'-a)\wedge (b'-b)$. If $F_f^{w}(c)=\infty$, \eqref{eq:mom} is obvious. Otherwise, for sufficiently large $n$, we have
\begin{align*}
(\eta_{[0,2^{-n}]}\circ u)\vee f_n^{w,c}&=(\eta_{[0,2^{-n}]}\circ u)+f_n^{h,b},\\
(\eta_{[0,2^{-n}]}\circ u)\wedge f_n^{w,c}&=f_n^{g,a}.
\end{align*}
We also have
\[f_n^{h,b}=\eta_{[2^{-n},2^{1-n}]}\circ\left(\eta_{[0,2^{-n}]}\circ v'+f_n^{h,b}\right)\]
and
\[\left((\eta_{[0,2^{-n}]}\circ u)+f_n^{h,b}\right)\wedge(\eta_{[0,2^{-n}]}\circ v')=f_n^{h,b}.\]
Applying \ref{cond:rstrsubad} twice, we obtain
\begin{align*}
&\Mc{E}(f_n^{g,a})+\Mc{E}(f_n^{h,b})\\
\le{}&\Mc{E}(f_n^{g,a})+\Mc{E}((\eta_{[0,2^{-n}]}\circ u)+f_n^{h,b})
+\Mc{E}(\eta_{[0,2^{-n}]}\circ v')\\
\le{}&\Mc{E}(f_n^{w,c})+\Mc{E}(\eta_{[0,2^{-n}]}\circ u)+\Mc{E}(\eta_{[0,2^{-n}]}\circ v').
\end{align*}
Letting $n\to\infty$ and using Lemmas \ref{lem:vanseq}, \ref{lem:Fsubad} and Proposition \ref{prop:pushforward}, we obtain
\[F_f^g(a)+F_f^h(b)\le F_f^w(c).\]
Taking the supremum over all such $g,h,a$ and $b$, we obtain $\Pi\langle f\rangle(U)+\Pi\langle f\rangle(V)\le \Pi\langle f\rangle(U\cup V)$ by \eqref{eq:emminus}. As in the proof of \cite[Theorem 5.6]{Sas2026}, the case of arbitrary sets is reduced to this case.
Now we fix any $f\in\Mc{M}(\Mc{E})$ and let $a\in\Mb{R}$. By Lemma \ref{lem:zero-zero} and right continuity of $F_f^{|f-a|-|a|}$ at $-|a|$, there exists $\epsilon_{(a)}>0$ such that
\[\Pi\langle f\rangle(\{a-\epsilon_{(a)}\le f\le a+\epsilon_{(a)}\})=\Pi\langle f\rangle(\{|f-a|-|a|\le \epsilon_{(a)}-|a|\})<\infty\]
and hence,
\begin{multline*}
\Pi\langle f\rangle (X)=\Pi\langle f\rangle (\{-\|f\|_{\infty}\le f\le \|f\|_{\infty}\}) \\
\le \sum_{i=1}^N \Pi\langle f\rangle(\{a_i-\epsilon_{(a_i)}\le f\le a_i+\epsilon_{(a_i)}\})<\infty
\end{multline*}
for some $N$ and $\{a_i\}_{i=1}^N$. Moreover, we obtain $\Pi \langle f\rangle (\{f=0\})=0$ and thus $\Pi\langle f\rangle (X)=\Pi\langle f\rangle (\supp(f))$ by Lemma \ref{lem:zero-zero}. Therefore, $\Pi\langle f\rangle$ is Radon.
\end{proof}
\begin{rem}
Lemma \ref{lem:zero-zero} and Theorem \ref{thm:Radon} imply that for any $f\in\Mc{M}(\Mc{E})$, $T_n\circ f\in\Mc{D}^i(\Mc{E})$ for some $n$. \ref{cond:rvanish} and $\Mc{M}(\Mc{E})\subset C_b(X)$ yield this phenomenon, and this may be also the case under Assumption \ref{asm:secpf}.
\end{rem}
\begin{cor}\label{cor:mustfinite}
For any $f,g\in\Mc{M}(\Mc{E})$, $t\mapsto F_f^g(t)$ is a distribution function of a finite Radon measure. Moreover, $\Pi\langle f\rangle(\{g\le a\})=F_f^g(a)=\lim_{n\to\infty}\Mc{E}(f_n^{g,a})$ for any $a\in\Mb{R}$.
\end{cor}
\begin{proof}
It follows from Proposition \ref{prop:pushforward}, Lemma \ref{lem:zero-zero} and Theorem \ref{thm:Radon}.
\end{proof}
Hereafter, we use the same notation $\Pi\langle f\rangle$ for its restriction to $\Mc{B}(X)$ for simplicity.
\begin{lem}\label{lem:setseq}
Let $\{\mu_i\}_{i=1}^n$ be a finite family of $[-\infty,\infty]$-valued monotone functions on $\Mc{B}(X)$, in other words, $\mu_i(A)\le \mu_i(B)$ whenever $A\subset B$.
\begin{enumerate}
	\item Suppose $\mu_i$ is inner regular on open sets, that is,
	\begin{equation}\label{eq:innneopen}
	\mu_i(U)=\sup_{V\Subset U,\ V\text{: open}}\mu_i(V)
	\end{equation}
	for any $i$ and open set $U$. Let $\Mc{G}$ be a dense subset of $(C_c(X),\|\cdot\|_{\infty})$. Then, for any compact set $E$ and open set $U$ with $E\subset U$, there exists $\{g_k\}_{k=1}^\infty\in \Mc{G}^\Mb{N}$ such that
	\[E\subset \{g_k< -1\}\Subset \{g_{k+1}\le -2\}\Subset U\]
	for any $k$ and $\lim_{k\to\infty}\mu_i(\{g_k< -1\})=\lim_{k\to\infty}\mu_i(\{g_k\le -2\})=\mu_i(U)$ for any $i$.
	\item Suppose $\mu_i$ is outer regular, that is,
	\[ \mu_i(A)=\inf_{U\supset A,\ U\text{: open}}\mu_i(U).\]
	Then, for any open set $U$ and $A\in\Mc{B}(X)$ with $A\subset U$, there exists a sequence $\{U_k\}_{k=1}^\infty$ of open sets such that
	\[A\subset U_{k+1}\subset U_k\subset U\]
	for any $k$ and $\lim_{k\to\infty}\mu_i(U_k)=\mu_i(A)$ for any $i$.
\end{enumerate}
\end{lem}
\begin{rem}
Since $X$ is locally compact and Hausdorff, and $\mu_i$ is monotone, \eqref{eq:innneopen} is equivalent to
\[\mu_i(U)=\sup_{K\subset U,\ K\text{: compact}}\mu_i(K).\]
\end{rem}

\begin{proof}
\begin{enumerate}
	\item By the inner regularity, there exists a family $\{V_{i,k}\}_{1\le i\le n}^{1\le k}$ of open sets such that $V_{i,k}\Subset U$ for any $i$ and $k$, and $\lim_{k\to\infty}\mu_{i}(V_{i,k})=\mu_i(U)$. Let $E_1=E\cup (\Cup_{i=1}^n \ol{V_{i,1}})$. Now we inductively define $E_k$ and $g_k$ as follows. For given $E_k$, we choose $g_k\in\Mc{G}$ such that $E_k\subset \{g_k\le -2\}$ and $\{g_k<-1\}\Subset U$, by Lemma \ref{lem:slider}. Then, we also define $E_{k+1}=E_k\cup \ol{\{g_k<-1\}}\cup (\Cup_{i=1}^n \ol{V_{i,k+1}} )$. Here,
	\[E\subset E_k \subset \{g_k< -1\}\Subset E_{k+1}\subset\{g_{k+1}\le -2\}\Subset\{g_{k+1}<-1\}\Subset U\] 
	and
	\begin{multline*}
\mu_i(U)=\lim_{k\to\infty}\mu_i(V_{i,k})\le \lim_{k\to\infty}\mu_i(\{g_k\le -2\})\\ \le \lim_{k\to\infty}\mu_i(\{g_k< -1\})=\mu_i(U)
	\end{multline*}
	for any $i$, by monotonicity.
	\item By the outer regularity, there exists a family $\{U_{i,k}\}_{1\le i\le n}^{1\le k}$ of open sets such that $U_{i,k}\supset A$ for any $i$ and $k$, and $\lim_{k\to\infty}\mu_{i}(U_{i,k})=\mu_i(A)$. Then, by monotonicity, it is evident that $U_k=U\cap(\Cap_{1\le i\le n}^{1\le j\le k}U_{i,j} ) (k\in\Mb{N})$ satisfies the desired inclusion relation and equality.
\end{enumerate}
\end{proof}
\begin{prop}\label{prop:localtoconvex}
If $\Mc{E}$ satisfies \ref{cond:rstrloc}, then the functional $f\mapsto \Pi\langle f\rangle(A)$ is convex on $\Mc{M}(\Mc{E})$, for any $A\in\Mc{B}(X)$.
\end{prop}
\begin{proof}
Since Lemma \ref{lem:setseq} holds, it is sufficient to prove the claim when $A=\{h\le a\}$ for arbitrary $a<0$ and $h\in\Mc{M}(\Mc{E})$. We first note that
\begin{align}
T_n(T_n(t)+T_n(s))&=T_n(t+s)\vee T_n(t-s), \label{eq:localmagic1}\\
|T_n(t)-T_n(s)|&=T_n(t+s)\wedge T_n(t-s) \label{eq:localmagic2}
\end{align}
for any $t,s\in\Mb{R}$ and $n\in\Mb{N}$ (see Figure \ref{fig:localmagic}). 
\begin{figure}
	\input{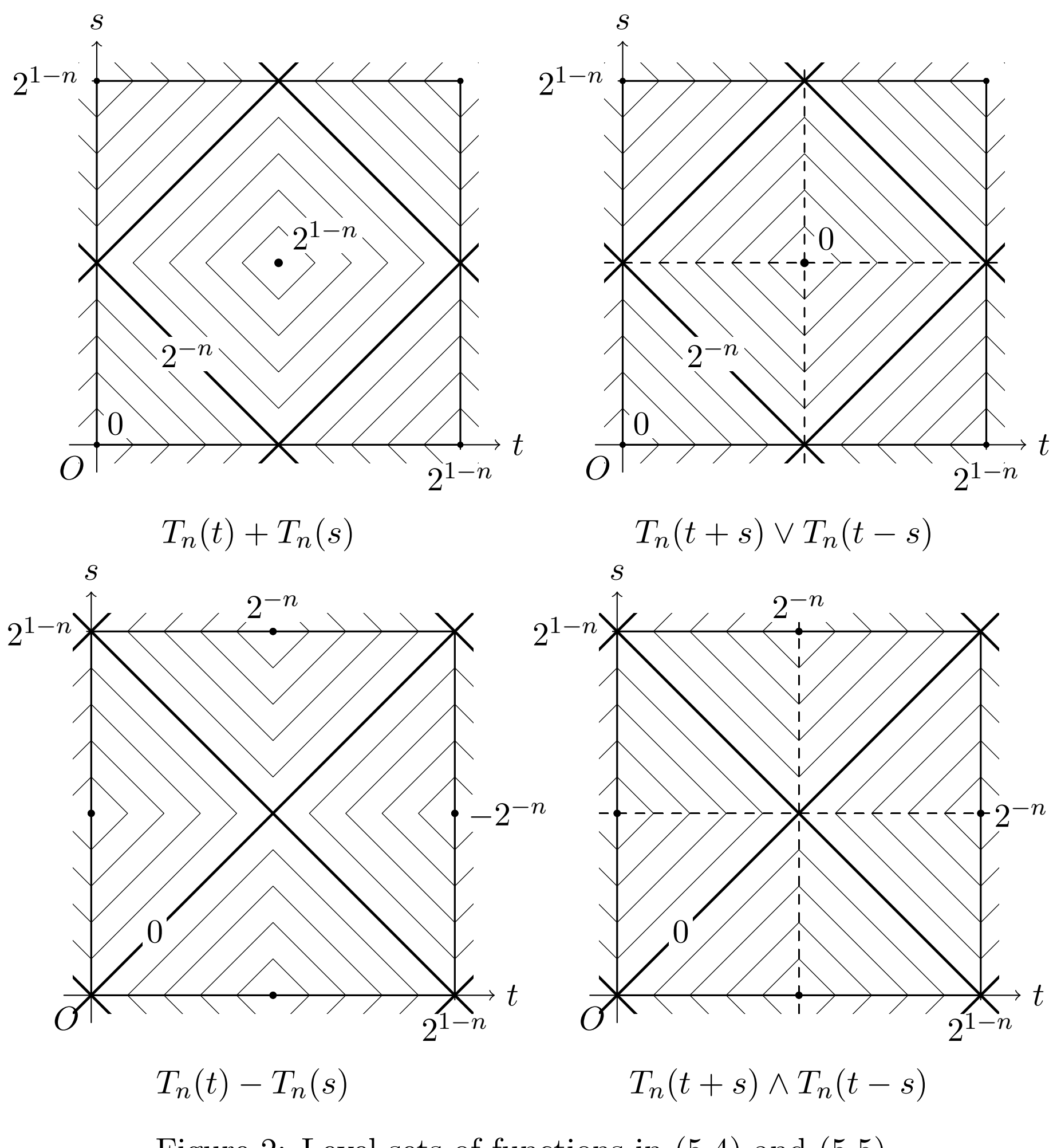}
	\caption{Level sets of functions in \eqref{eq:localmagic1} and \eqref{eq:localmagic2}}\label{fig:localmagic}
\end{figure}
Fix any $f,g\in\Mc{M}(\Mc{E})$ and $a'>a$. By Corollary \ref{cor:bumpfcn}, there exists $u:X\to[0,1]$ such that $u\in\Mc{M}(\Mc{E})$ and  $u|_{E}\equiv 1$ for $E\coloneq \supp(f)\cup\supp(g)\cup\supp(h)$. Let
\begin{align*}
	v_n&= \left(\eta_{[0,2^{-n}]}+\eta_{[1-2^{-n},1]}\right)\circ \left( \left(1+\eta_{[a',a'+2^{-1}]}(a')\right)u-\eta_{[a',a'+2^{-1}]}\circ h  \right)\\
	w_n&=\eta_{[1-2^{1-n},1-2^{-n}]}\circ \left( \eta_{[a',a'+1]}(a'+1)u-\eta_{[a',a'+1]}\circ h  \right)
\end{align*}
(note that $v_n,w_n\in\Mc{M}(\Mc{E})$). Then
\begin{align}
	&T_n\circ((T_n\circ v_n)\vee (T_n\circ f+T_n\circ g)\wedge v_n)\notag\\
	=&\begin{cases}
		(T_n\circ v_n)\vee (T_n\circ f+T_n\circ g)&\text{(on $E\cap \{T_n\circ f+T_n\circ g< 2^{-n}\} $)}\\
		T_n\circ ((T_n\circ f+T_n\circ g)\wedge v_n)&\text{(on $E\cap \{T_n\circ f+T_n\circ g\ge 2^{-n}\} $)}\\
		T_n\circ v_n&\text{(on $X\setminus E $)}
	\end{cases}\notag\\
	=&(T_n\circ v_n)\vee (T_n\circ (T_n\circ f+T_n\circ g))\notag\\
	\overset{\eqref{eq:localmagic1}}{=}&(T_n\circ v_n)\vee(T_n\circ(f+g))\vee(T_n\circ (f-g))\notag\\
	=&\begin{cases}
		(T_n\circ v_n)\vee(T_n\circ(f+g))\vee(T_n\circ (f-g)) &\text{(on $E\cap \{h<a'+2^{-n}\}$)}\\
		2^{-n}&\text{(on $E\cap\{h\ge a'+2^{-n}\} $)}\\
		T_n\circ v_n&\text{(on $X\setminus E$)}
	\end{cases}\notag\\
=&(T_n\circ v_n)\vee (f+g)_n^{h,a'}\vee(f-g)_n^{h,a'}.\label{eq:magicfold1}
\end{align}
Moreover,
\begin{align}
	&\left| (-w_n)\vee(T_n\circ f-T_n\circ g)\wedge (w_n) \right|\notag\\
	=&w_n\wedge|T_n\circ f-T_n\circ g|\notag\\
	\overset{\eqref{eq:localmagic2}}{=}&w_n\wedge(T_n\circ (f+g))\wedge (T_n\circ (f-g))\notag\\
	=&\begin{cases}
		(S_n^{a'}\circ h)\wedge(T_n\circ (f+g))\wedge (T_n\circ (f-g))&\text{(on $E\cap \{h<a'+2^{1-n}\}$)}\\
		0&\text{(otherwise)}\\
	\end{cases}\notag\\
	=&(f+g)_n^{h,a'}\wedge(f-g)_n^{h,a'}.\label{eq:magicfold2}
\end{align}
Therefore, we obtain
\begin{align}
	&\begin{multlined}[t][.9\linewidth]
		\Mc{E}\left((T_n\circ v_n)\vee (T_n\circ f+T_n\circ g)\wedge v_n\right)\\
		+\Mc{E}\left((-w_n)\vee(T_n\circ f-T_n\circ g)\wedge (w_n) \right)
	\end{multlined}\notag\\
	=&\Mc{E}\left((T_n\circ v_n)\vee (f+g)_n^{h,a'}\vee(f-g)_n^{h,a'}\right)+\Mc{E}\left( (f+g)_n^{h,a'}\wedge(f-g)_n^{h,a'} \right)\notag\\
	\overset{\ref{cond:rstrsubad}}{\le}&
		\Mc{E}(T_n\circ v_n)+\Mc{E}\left((f+g)_n^{h,a'}\vee(f-g)_n^{h,a'}\right)
	+\Mc{E}\left((f+g)_n^{h,a'}\wedge(f-g)_n^{h,a'}\right)\notag
	\\
		\overset{\ref{cond:rstrsubad}}{\le}&\Mc{E}(T_n\circ v_n)+\Mc{E}((f+g)_n^{h,a'})+\Mc{E}((f-g)_n^{h,a'})\label{eq:localmagic3}
\end{align}
by Lemma \ref{lem:locfold}.
Here Lemmas \ref{lem:limeval2} and \ref{lem:vanseq} yield
\begin{multline}
	F_f^h(a)\le \Mc{E}\Big(\frac{1}{2}\Big(\left((T_n\circ v_n)\vee (T_n\circ f+T_n\circ g)\wedge v_n\right) \\+\left((-w_n)\vee(T_n\circ f-T_n\circ g)\wedge w_n\right) \Big)\Big) \label{eq:localmagic4}
\end{multline}
because
\begin{multline*}
	f_m^{h,a}=T_m\circ\Big(\frac{1}{2}\Big(\left((T_n\circ v_n)\vee (T_n\circ f+T_n\circ g)\wedge v_n\right) \\+\left((-w_n)\vee(T_n\circ f-T_n\circ g)\wedge w_n\right) \Big)\Big)\wedge (S_m^a\circ h)
\end{multline*}
for sufficiently large $m>n$ with $a+2^{1-m}<a'$, since $S_m^a(0)=0$ for such $m$.
Hence,
\[\Pi\langle f\rangle(\{h\le a\})=F_f^h(a)\le\frac{\Mc{E}(T_n\circ v_n)+\Mc{E}((f+g)_n^{h,a'})+\Mc{E}((f-g)_n^{h,a'})}{2} \]
follows from \eqref{eq:localmagic3}, \eqref{eq:localmagic4} and \ref{cond:rconv}. Letting $n\to\infty$ and $a'\downarrow a$, we obtain
\begin{align*}
\Pi\langle f\rangle(\{h\le a\})&\le \frac{F_{f+g}^h(a)+F_{f-g}^h(a)}{2}\\
&=\frac{\Pi\langle f+g\rangle(\{h\le a\})+\Pi\langle f-g\rangle(\{h\le a\})}{2}.
\end{align*}
By the same argument, we conclude that $t\mapsto \Pi\langle f+tg\rangle(\{h\le a\})$ is midpoint convex on $(0,1)$. Moreover,
\[\Pi\langle f+tg\rangle(\{h\le a\})\le \Pi\langle f+tg\rangle (X)=\Mc{E}(f+tg)\le \Mc{E}(f)+\Mc{E}(f+g)<\infty\]
for $t\in(0,1)$, by Lemmas \ref{lem:MisD} and \ref{lem:locfold}. This yields convexity of $t\mapsto \Pi\langle f+tg\rangle(\{h\le a\})$ on $(0,1)$, by the Bernstein--Doetsch Theorem (see \cite[Theorem 6.4.2]{Kuc2009}\footnote{Note that in this reference, ``convex function'' means a midpoint convex function: see \cite[Subsection 5.3]{Kuc2009}.} for example).
\end{proof}
In general, we do not know the convexity of $f\mapsto \Pi\langle f\rangle(X)$, so we need a ``convex envelope'' of $\Pi\langle f\rangle$.
\begin{defi}\label{defi:themeas}
Let $A\in\Mc{B}(X)$. We define
\begin{align*}
&\Gamma_\Mc{E}^c\langle f\rangle(A)\\
=&\inf\left\{\sum_{i=1}^m \sum_{j=1}^n t_{i,j} \Pi\langle f_{i,j}\rangle (A_i)\Bigg\vert \begin{minipage}{.6\linewidth}
	$m,n\in\Mb{N}, 0\le t_{i,j}\le 1, f_{i,j}\in\Mc{M}(\Mc{E}),\ A_i\in \Mc{B}(X)\ (1\le i\le m, 1\le j\le n)$ such that $A=\bigsqcup_{i=1}^m A_i$,
	$\sum_{j=1}^n t_{i,j}=1$ and $\sum_{j=1}^n t_{i,j} f_{i,j}=f$ for any $i$.
\end{minipage} \right\} 
\end{align*}
for general $f\in\Mc{M}(\Mc{E})$, and $\Gamma^c_{\Mc{E}} \langle f\rangle(A)=0\ (A=\emptyset), \infty\ (\text{otherwise})$ for $f\in V\setminus \Mc{M}(\Mc{E})$. We also define $\Mc{E}^{(c)}(f)=\Gamma^c_{\Mc{E}}\langle f\rangle(X)$ for any $f\in V$.
Similarly to $\Pi$, we abbreviate the subscript $\Mc{E}$, when no confusion can arise.
\end{defi}
Before the proof of Theorem \ref{thm:main1}, we also prepare a lemma for the uniqueness part.

\begin{lem}\label{lem:fundestm}
	Let $\{\mu^c\langle f\rangle\}_{f\in\Mc{M}(\Mc{E})}$ be a family of finite Radon measures with \ref{cond:mene}--\ref{cond:mlocal}.
	For $f,g\in \Mc{M}(\Mc{E})$, $a\in\Mb{R}$ and $A\in\Mc{B}(X)$ we have
	\begin{align}
		\mu^c \langle f_n^{g,a}\rangle(A)&\ge 
			\mu^c \langle f\rangle (A\cap\{g\le a\}),\\
		\mu^c \langle f_n^{g,a}\rangle(A)&\le \begin{multlined}[t][.6\linewidth]
			\mu^c \langle f\rangle (A\cap\{g<a+2^{1-n}\})\\+\mu^c \langle g\rangle (A\cap \{a+2^{-n}<g<a+2^{1-n}\}).
		\end{multlined}\label{eq:fundestm}
	\end{align}
	In particular, $\mu^c\langle f\rangle (\{g\le a\})=\lim_{n\to\infty}\Mc{E}^{(c)}(f_n^{g,a})$.
\end{lem}
\begin{proof} Let $\|f\|_{\infty}<M\in\Mb{N}$. By \ref{cond:mgsm} and \ref{cond:mlocal}, we obtain
	\begin{align*}
		&\mu^c \langle f_n^{g,a}\rangle(A)\\
		\ge &
			\mu^c \langle T_n\circ f\rangle(A\cap \{g\le a\})\\
		=&
			\sum_{k=-2^nM}^{2^nM-1}\mu^c \langle T_n\circ f\rangle (A\cap\{2^{-n}k<f\le2^{-n}(k+1), g\le a\})\\
		=&\mu^c\langle f\rangle (A\cap\{g\le a\})
	\end{align*}
	and
	\begin{align*}
		&\mu^c \langle f_n^{g,a}\rangle(A)\\
		=&\begin{multlined}[t][.9\linewidth]
			\mu^c \langle T_n\circ f\rangle(A\cap \{0<T_n\circ f\le S_n^a\circ g\})\\
			+\mu^c \langle S_n^a\circ g-S_n^a(0)\rangle(A\cap \{0<S_n^a\circ g<T_n\circ f\})\\
			+\mu^c \langle 0\rangle (A\cap\{(T_n\circ f)\wedge (S_n^a\circ g)=0\} )
		\end{multlined}\\
		\le&\begin{multlined}[t][.9\linewidth]
			\sum_{k=-2^nM}^{2^nM-1}\mu^c \langle T_n\circ f\rangle (A\cap\{2^{-n}k<f\le2^{-n}(k+1), g<a+2^{1-n}\})\\+\mu^c\langle -g\rangle (A\cap \{a+2^{-n}<g<a+2^{1-n}\})+0
		\end{multlined}\\
		=&\mu^c\langle f\rangle (A\cap\{g<a+2^{1-n}\})+\mu^c\langle g\rangle (A\cap \{a+2^{-n}<g<a+2^{1-n}\}).
	\end{align*}
	Since $\mu^c\langle g\rangle (\{a<g<a+1\})\le\Mc{E}^{(c)}(g)<\infty$, letting $A=X$ and $n\to\infty$ yields $\mu^c\langle f\rangle (\{g\le a\})=\lim_{n\to\infty}\Mc{E}^{(c)}(f_n^{g,a})$.
\end{proof}

\begin{proof}[Proof of Theorem \ref{thm:main1}~(1) and (2)]
We first prove (2). Since $\Pi\langle 0\rangle\equiv 0$ and $\Pi\langle f\rangle=\Pi\langle -f\rangle$ by $T_n\circ f=T_n\circ(-f)$, it is easy to show that $f\mapsto \Gamma^c\langle f\rangle(A)$ is a gsm, where \ref{cond:rlc} comes from $\Pi\langle f\rangle(X)<\infty$ for any $f$. It also easily follows that $\Gamma^c\langle f\rangle$ is finitely additive.
If $\bigsqcup_{i=1}^\infty A_i=A$ for $A\in\Mc{B}(X)$ then
\[\Gamma^c \langle f\rangle \left(A\setminus \bigsqcup_{i=1}^n A_i \right)\le \Pi \langle f\rangle \left(A\setminus \bigsqcup_{i=1}^n A_i \right)\to 0\quad(n\to\infty),\]
which proves $\Gamma^c\langle f\rangle|_A$ is a Borel measure. Consequently, $\Gamma^c\langle f\rangle\le \Pi\langle f\rangle$ yields that $\Gamma^c\langle f\rangle$ is a finite Radon measure with $\Gamma^c\langle f\rangle(\supp(f))=\Mc{E}^{(c)}(f)$.

Moreover, if $(f-g)|_A\equiv a$ for $A\in\Mc{B}(X)$ and $a\in\Mb{R}$, then $\Gamma^c\langle b(f-g)\rangle(A)\le \Pi \langle b(f-g)\rangle(\{b(f-g)=b a\})= 0$ for any $b\in\Mb{R}$ by Lemma \ref{lem:zero-zero}. Thus, \ref{cond:mlocal} also follows from Lemma \ref{lem:lsc}. Here  Lemmas \ref{lem:setseq} and \ref{lem:fundestm} ensure the uniqueness of $\{\Gamma^c\langle f\rangle\}_{f\in\Mc{M}(\Mc{E})}$.

Since (2) and $\Mc{E}^{(c)}\le \Mc{E}$ hold, it is routine to prove that $\Mc{E}^{(c)}$ satisfies \ref{cond:rconv}--\ref{cond:rstrloc}. For example, for any $f\in\Mc{M}(\Mc{E})$ and piecewise affine function $\varphi(t)=\sum_{k\in\Mb{Z}} a_k\eta_{[b_k,b_{k+1}]}(t)$ with $|a_k|\le 1$, we obtain
\begin{multline*}
\Mc{E}^{(c)}(\varphi\circ f)\overset{\ref{cond:mene},\ref{cond:mlocal}}{=}\sum_{k\in\Mb{Z}} \Gamma^c\langle a_k\eta_{[b_k,b_{k+1}]}\circ f\rangle(\{b_k\le f<b_{k+1}\})\\
\overset{\ref{cond:mgsm}}{\le}\sum_{k\in\Mb{Z}} \Gamma^c\langle \eta_{[b_k,b_{k+1}]}\circ f\rangle(\{b_k\le f<b_{k+1}\})\overset{\ref{cond:mene},\ref{cond:mlocal}}{=}\Mc{E}^{(c)}(f).
\end{multline*}
Next we show the maximality of $\Mc{E}^{(c)}$. Since $\Mc{E}^{(c)}=\infty$ on $V\setminus\Mc{M}(\Mc{E})$, it is sufficient to show the maximality on $\Mc{M}(\Mc{E})$. Let $E$ be a gsm on $\Mc{M}(\Mc{E})$ with $E\le\Mc{E}$ and \ref{cond:rwNC}--\ref{cond:rstrloc}. By Lemma \ref{lem:MisD}, $E$ is finite on $\Mc{M}(\Mc{E})$. Then, we obtain \[E(f)=\lim_{n\to\infty}E(T_n\circ f)=\Pi_{E}\langle  f\rangle(X)=\Gamma^c_{E}\langle f\rangle(X) \]
for $f\in\Mc{D}(\Mc{E}^{(c)})=\Mc{M}(\Mc{E})$ by Lemma \ref{lem:locfold} and Proposition \ref{prop:localtoconvex}. Therefore
\[E(f)=\Gamma^c_{E}\langle f\rangle(X)\le\Gamma^c_{\Mc{E}}\langle f\rangle (X)= \Mc{E}^{(c)}(f) \]
follows from $E\le \Mc{E}$ when $f\in\Mc{M}(\Mc{E})$: otherwise $\Mc{E}^{(c)}(f)=\infty\ge E(f)$.
\end{proof}

\begin{lem}[{cf. \cite[Theorem 3.1]{Sas2026}}]\label{lem:dom}
	Assume that $E$ is a gsm with \ref{cond:rconv}--\ref{cond:rstrsubad} such that $\Mc{M}(\Mc{E})=\Mc{M}(E)$, and let $a\in\Mb{R}$ and $b\ge0$.
	\begin{enumerate}
		\item If $\Mc{E}(a f)\le b E(f)$ for any $f\in\Mc{M}(\Mc{E})$ then $\Gamma_\Mc{E}^c\langle a f\rangle \le b \Gamma_E^c\langle f\rangle$ for any $f\in\Mc{M}(\Mc{E})$.
		\item If $\Mc{E},E$ satisfy \ref{cond:rstrloc} and $|\Mc{E}(a f)-\Mc{E}(a g)|\le b |E(f)-E(g)|$ for any $f,g\in\Mc{M}(\Mc{E})$, then $|\Gamma_\Mc{E}^c\langle a f\rangle-\Gamma_\Mc{E}^c\langle a g\rangle| \le b |\Gamma_E^c\langle f\rangle-\Gamma_E^c\langle g\rangle|$ for any $f,g\in\Mc{M}(\Mc{E})$.
	\end{enumerate}

\end{lem}
\begin{proof}\begin{enumerate}
		\item By \ref{cond:rsym}, it is sufficient to prove the case $a >0$.
		By Lemma \ref{lem:anotherwave} and Corollary \ref{cor:mustfinite}, we obtain
		\begin{multline*}
			\Pi_\Mc{E}\langle a f\rangle(\{g\le c\})=\lim_{n\to\infty}\Mc{E}((T_{n,a}\circ(a f))\wedge (a S_n^c\circ g))\\ \le\lim_{n\to\infty}b E(f_n^{g,c})=b \Pi_E\langle f\rangle (\{g\le c\})
		\end{multline*}
		for any $f,g\in\Mc{M}(\Mc{E})$ and $c\in\Mb{R}$. Hence, the claim follows from the definition of $\Gamma_\Mc{E}^c, \Gamma_E^c$ and Lemma \ref{lem:setseq}.
		\item $\Pi_\Mc{E}=\Gamma_\Mc{E}^c$ and $\Pi_E=\Gamma_E^c$ by Proposition \ref{prop:localtoconvex}. Hence, for any $f,g\in\Mc{M}(\Mc{E})$ we obtain
		\[|\Gamma_\Mc{E}^c\langle a f\rangle(A)-\Gamma_\Mc{E}^c\langle a g\rangle(A)|\le b |\Gamma_E^c\langle  f\rangle(A)-\Gamma_E^c\langle g\rangle(A)| \]
		when $A=\{h\le c\}$ for some $h\in\Mc{M}(\Mc{E})$ and $c\in\Mb{R}$, and thus the claim follows from Lemma \ref{lem:setseq}.
	\end{enumerate} 
\end{proof}
\begin{proof}[Proof of Theorem \ref{thm:main1}~(3)]
	$p$-homogeneity follows from Lemma \ref{lem:dom}~(1). Then, 
	\[\Gamma^c\langle f\rangle(A)^{\frac{1}{p}}=\|f\|_{LN(f\mapsto \Gamma^c\langle f\rangle (A))} \]
	by $p$-homogeneity, and thus the triangle inequality \eqref{eq:tri} holds.
\end{proof}

The next lemma is stated here to avoid imposing unnecessary assumptions in later sections, but it is not used immediately, so the reader can skip it on a first reading.

\begin{lem}\label{lem:boundinherit}
	Let $\alpha>0, \varphi\in C(\Mb{R}^2)$ such that
	\[\Mc{E}(\alpha(f+g))+\Mc{E}(\alpha(f-g))\le \varphi\left(\Mc{E}(f),\Mc{E}(g) \right)\]
	for any $f,g\in\Mc{D}(\Mc{E})$.
	\begin{enumerate}
		\item For any $f,g\in\Mc{M}(\Mc{E})$ and $A\in\Mc{B}(X)$, we obtain
			\begin{equation}
			\Pi\langle \alpha(f+g)\rangle (A)+\Pi\langle \alpha(f-g)\rangle (A)\le \varphi\left(\Pi\langle f\rangle (A),\Pi\langle g\rangle (A) \right). \label{eq:piinherit}
		\end{equation} In particular, if \ref{cond:rstrloc} holds then
		\begin{equation}\label{eq:boundinherit}
			\Gamma^c\langle \alpha(f+g)\rangle (A)+\Gamma^c\langle \alpha(f-g)\rangle (A)\le \varphi\left(\Gamma^c\langle f\rangle (A),\Gamma^c\langle g\rangle (A) \right).
		\end{equation}
		\item If $\varphi$ is superlinear and 1-homogeneous on nonnegative numbers, that is,
		\[\varphi(a+b,c+d)\ge \varphi(a,c)+\varphi(b,d),\quad \varphi(\lambda a,\lambda b)=\lambda\varphi(a,b) \]
		for any $a,b,c,d,\lambda\ge0$, then \eqref{eq:boundinherit} holds for any $f,g\in\Mc{M}(\Mc{E})$ and $A\in\Mc{B}(X)$.
	\end{enumerate}
\end{lem}
\begin{proof}
	Let $f,g\in\Mc{M}(\Mc{E})$. Note that, if $U$ is open and $U\subset \{\varepsilon f=g+C\}$ for some $\varepsilon\in \{\pm 1\}$ and $C\in\Mb{R}$, then $\Pi\langle f\rangle(A)=\Pi\langle g\rangle(A)$ for any $A\subset U$. Indeed, if $h\in\Mc{M}(\Mc{E})$ and $a\in\Mb{R}$ satisfy $\{h\le a\}\subset U$, then $\{h\le a+2^{1-n}\}\subset U$ for sufficiently large $n$. Then, on $U$, $T_{n}\circ f=T_{n}\circ g$ if $C=0$, and $T_{n,|C|}\circ f=T_{n,|C|}\circ g$ otherwise. Hence, Lemma \ref{lem:anotherwave} and definition of $\Pi\langle f\rangle$ yield this equality.
	\begin{enumerate}
			\item Assume that $A=\{h\le a\}$ for some $h\in\Mc{M}(\Mc{E})$ and $a\in\Mb{R}$.  Fix any $f,g\in\Mc{M}(\Mc{E})$, and also fix sufficiently large $m$ such that $f_m^{h,a},g_m^{h,a}\in \Mc{D}^i(\Mc{E})$. For any $l>m$, we write
			\[\delta_{l}(t)=2\sum_{k\in\Mb{Z}}(-1)^{k+1}\left(\eta_{[0,2^{-l}]}(t-2^{-m}k)-\eta_{[2^{-l},2^{1-l}]}(t-2^{-m}k) \right). \]
			We also define
			\[\square^{l}_{j,k}=\left\{
			\frac{j}{2^m}+\frac{1}{2^{l-1}}<f<\frac{j+1}{2^m}+\frac{1}{2^l},\
			\frac{k}{2^m}+\frac{1}{2^{l-1}}<g<\frac{k+1}{2^m}+\frac{1}{2^l}
			\right\}\]
			and
			\[X^{e}_{l}=\Cup_{j,k\in\Mb{Z}: k+j\text{: even}}\square^{l}_{j,k},\quad X^{o}_{l}=\Cup_{j,k\in\Mb{Z}: k+j\text{: odd}}\square^{l}_{j,k}.\]
			Then,
			\[f+g=\varepsilon\left(f_m^{h,a}+g_m^{h,a}+\delta_{l}\circ f+\delta_{l}\circ g \right)+C\]
			on $\{h<a+2^{-m}\}\cap\square^{l}_{j,k}$, for some $\varepsilon\in\{\pm 1\}$ and $C\in\Mb{R}$, when $k+j$ is even. Hence,
			\begin{align*}
			&\Pi\langle \alpha(f+g)\rangle (A\cap X^{e}_{l})\\
			=&\Pi\left\langle \alpha(f_m^{h,a}+g_m^{h,a}+(\delta_{l}\circ f)+(\delta_{l}\circ g))\right\rangle (A\cap X^{e}_{l}).
			\shortintertext{Similarly, we have}
			&\Pi\langle \alpha(f+g)\rangle (A\cap X^{o}_{l})\\
			=&\Pi\left\langle \alpha(f_m^{h,a}-g_m^{h,a}+(\delta_{l}\circ f)-(\delta_{l}\circ g))\right\rangle (A\cap X^{o}_{l}),\\
			&\Pi\langle \alpha(f-g)\rangle (A\cap X^{e}_{l})\\
			=&\Pi\left\langle \alpha(f_m^{h,a}-g_m^{h,a}+(\delta_{l}\circ f)-(\delta_{l}\circ g))\right\rangle (A\cap X^{e}_{l}),\\
			&\Pi\langle \alpha(f-g)\rangle (A\cap X^{o}_{l})\\
			=&\Pi\left\langle \alpha(f_m^{h,a}+g_m^{h,a}+(\delta_{l}\circ f)+(\delta_{l}\circ g))\right\rangle (A\cap X^{o}_{l}).
			\end{align*}
			Since $\|\delta_l\circ f\|_{LN(\Mc{E})}, \|\delta_l\circ g\|_{LN(\Mc{E})}\to 0$ as $l\to\infty$ by Lemma \ref{lem:vanseq}, these yield
			\begin{align*}
				&\Pi\langle \alpha(f+g)\rangle (A\cap (X^e_{l}\cup X^o_{l}))+\Pi\langle \alpha(f-g)\rangle (A\cap (X^e_{l}\cup X^o_{l}))\\
				\le &\begin{multlined}[t][.9\linewidth]
				\Mc{E}\left(\alpha(f_m^{h,a}+g_m^{h,a}+(\delta_{l}\circ f)+(\delta_{l}\circ g))\right)\\
				+\Mc{E}\left(\alpha(f_m^{h,a}-g_m^{h,a}+(\delta_{l}\circ f)-(\delta_{l}\circ g))\right)
				\end{multlined}\\
				\le& \varphi\left(\Mc{E}(f_{m}^{h,a}+(\delta_l\circ f)),\Mc{E}(g_m^{h,a}+(\delta_l\circ g)) \right).
			\end{align*}
			Since $\indi_{X^e_l\cup X^o_l}\to\indi_X$ pointwise, letting $l\to\infty$ and then $m\to\infty$, we obtain \eqref{eq:piinherit} in this case, by the bounded convergence theorem and Lemma \ref{lem:limeval}. The general case follows from Lemma \ref{lem:setseq}. Now, assuming \ref{cond:rstrloc}, \eqref{eq:boundinherit} follows from Proposition \ref{prop:localtoconvex} (and definition of $\Gamma^c$).
	\item Fix any $\epsilon>0$ and let $t_{i,j}, f_{i,j}, s_{k,l}, g_{k,l}$ be as in the definition of $\Gamma^c$ so that
	\begin{gather*}
\sum_{i=1}^m \sum_{j=1}^n t_{i,j} \Pi\langle f_{i,j}\rangle (A_i)<\Gamma^c\langle f\rangle(A)+\epsilon,\\ \sum_{k=1}^M \sum_{l=1}^N s_{k,l} \Pi\langle g_{k,l}\rangle (\tilde{A}_k)<\Gamma^c\langle g\rangle(A)+\epsilon.
	\end{gather*}
		Since $A=\cup_{i,k} (A_i\cap \tilde{A}_k)$, replacing $A_i, \tilde{A}_k$ with $A_i\cap \tilde{A}_k$ if necessary, we may assume that $m=M$ and $A_i=\tilde{A}_i$ without loss of generality. Moreover, $f=\sum_{j,l}t_{i,j}s_{i,l} f_{i,j}$ and $g=\sum_{j,l}t_{i,j}s_{i,l} g_{i,l}$ for any $i$. Using the common coefficients $t_{i,j}s_{i,l}$ and reindexing the pairs $(j,l)$, we can also assume $N=n$ and $t_{i,j}=s_{i,j}$. Then,
	
	\begin{align*}
	&\varphi\left(\sum_{i=1}^m \sum_{j=1}^n t_{i,j} \Pi\langle f_{i,j}\rangle (A_i),\sum_{i=1}^m \sum_{j=1}^n t_{i,j} \Pi\langle g_{i,j}\rangle (A_i) \right)\\
	\ge&\sum_{i=1}^m \sum_{j=1}^n\varphi\left(t_{i,j} \Pi\langle f_{i,j}\rangle(A_i),t_{i,j}\Pi\langle g_{i,j}\rangle (A_{i}) \right)\\
	=&\sum_{i=1}^m \sum_{j=1}^n t_{i,j}\varphi\left( \Pi\langle f_{i,j}\rangle(A_i),\Pi\langle g_{i,j}\rangle(A_i) \right)\\
	\overset{\eqref{eq:piinherit}}{\ge}&\sum_{i=1}^m \sum_{j=1}^n t_{i,j}\Pi\langle \alpha(f_{i,j}+g_{i,j})\rangle (A_i)+\sum_{i=1}^m \sum_{j=1}^n t_{i,j}\Pi\langle \alpha(f_{i,j}-g_{i,j})\rangle (A_i)\\
		&\ge\Gamma^c\langle \alpha(f+g)\rangle(A)+\Gamma^c\langle \alpha(f-g)\rangle (A).
	\end{align*}
	Letting $\epsilon\to0$, we obtain \eqref{eq:boundinherit}.
\end{enumerate}
\end{proof}
For the last part of this section, we remark on the replacement of assumptions of the main theorem.
\begin{cor}\label{cor:mainbutfaulting}
Theorem \ref{thm:main1} holds even if we replace (each appearance of) the property \ref{cond:rstrsubad} with the following two properties, \ref{cond:rsubad} and \ref{cond:rfaulting} (with $\Mc{E}=\varrho$).
\begin{description}
\rcond[cond:rstrsubad]{cond:rsubad} {\rm(}$f\vee g\in V$ and{\rm)} $\varrho(f\vee g)\le\varrho(f)+\varrho(g)$ whenever $f,g\in\Mc{M}(\varrho)$.
	\rcond{cond:rfaulting} {\rm(}$(f-a)\vee g\wedge (f+a),(g-a)\vee f\wedge(g+a)\in V$ and{\rm)} $\varrho((f-a)\vee g\wedge (f+a))+\varrho((g-a)\vee f\wedge(g+a))\le\varrho(f)+\varrho(g)$ whenever $f,g\in\Mc{M}(\varrho)$ and $a\ge 0$.
\end{description}
\end{cor}
\begin{rem}
\ref{cond:rfaulting} appears in a characterization of nonlinear Dirichlet forms (see Definition \ref{defi:NLDF}; cf.~Proposition \ref{prop:localmeansnldc}). If \ref{cond:rwNC} holds, $(f-a)\vee g\wedge (f+a),(g-a)\vee f\wedge(g+a)\in\Mc{M}(\Mc{E})$ whenever $f,g\in\Mc{M}(\Mc{E})$, by Lemma \ref{lem:domcont}.
\end{rem}
\begin{proof} For the conclusion, it is evident that $\Mc{E}^{(c)}$ also satisfies \ref{cond:rsubad} and \ref{cond:rfaulting}: this follows from \ref{cond:mene}--\ref{cond:mlocal}, similarly to the other properties, or directly follows from Proposition \ref{prop:localmeansnldc}.
	
Among the proofs needed for Theorem \ref{thm:main1}, \ref{cond:rstrsubad} is used in the proofs of Lemma \ref{lem:Fsubad}, Theorem \ref{thm:Radon} and Proposition \ref{prop:localtoconvex} and in the proof of Theorem \ref{thm:main1} itself. For Lemma \ref{lem:Fsubad}, the use of \ref{cond:rsubad} instead of \ref{cond:rstrsubad} is sufficient to prove the claim. By virtue of Proposition \ref{prop:localmeansnldc}~(1), the proofs of Theorem \ref{thm:main1} and Proposition \ref{prop:localtoconvex} hold without assuming \ref{cond:rstrsubad} initially. Hence, we only need to prove Theorem \ref{thm:Radon}, in particular \eqref{eq:mom}, without \ref{cond:rstrsubad}.

Let $U,V\subset X$, $a,b,a',b',c<0$ and $g,h,u,v',w\in\Mc{M}(\Mc{E})$ be as in the original proof of Theorem \ref{thm:Radon}. We also let $u''\in\Mc{M}(\Mc{E})$ such that $u'':X\to[0,1]$, $u''|_{\supp(v')}\equiv 1$ and $u''|_{X\setminus V}\equiv 0$. As in the original proof, we may assume $F_f^w(c)<\infty$ without loss of generality. Here we write
$v_n=\eta_{[0,2^{-n}]}\circ u-\eta_{[0,2^{-n}]}\circ u''=\eta_{[-2^{-n},2^{-n}]}\circ(u-u'')$.
For $v_n$, we obtain
\begin{align*}
	(f_n^{w,c}-2^{-n})\vee v_n\wedge (f_n^{w,c}+2^{-n})=&(f_n^{w,c}-2^{-n})\vee v_n\\
	=&\begin{cases}
			(f_n^{g,a}-2^{-n})\vee 2^{-n}&\text{on }\{g\le a' \}\\
			(f_n^{h,b}-2^{-n})\vee(-2^{-n}) &\text{on }\{h\le b'\}\\
			(-2^{-n})\vee v_n& \text{otherwise}
		\end{cases}\\
	=&f_n^{h,b}+v_n
\end{align*}
and similarly,
\[(v_n-2^{-n})\vee f_n^{w,c}\wedge (v_n+2^{-n})=f_n^{w,c}\wedge (v_n+2^{-n})=f_n^{g,a} \]
for sufficiently large $n$. Moreover,
\[(f_n^{h,b}+v_n-2^{-n})\vee (\eta_{[0,2^{-n}]}\circ v')\wedge(f_n^{h,b}+v_n+2^{-n})=f_n^{h,b}\]
also holds for such $n$. Hence, we obtain
\begin{multline*}
	F_f^g(a)+F_f^h(b)\le \lim_{n\to\infty}(\Mc{E}(f_n^{g,a})+\Mc{E}(f_n^{h,b}+v_n)+\Mc{E}(\eta_{[0,2^{-n}]}\circ v'))\\
	\le \lim_{n\to\infty}\left(\Mc{E}(f_n^{w,c})+\Mc{E}(v_n)+\Mc{E}(\eta_{[0,2^{-n}]}\circ v')\right)=F_f^w(c)
\end{multline*}
by \ref{cond:rfaulting} and Lemma \ref{lem:vanseq}. This is sufficient to prove \eqref{eq:mom}, as in the original proof of Theorem \ref{thm:Radon}.
\end{proof}

\section{Lower semicontinuous extensions}\label{sec:lsc}


When we consider Dirichlet-form-like frameworks, it is usual to regard the domain of the energy as a subspace of an $L^q$ space and to assume that the energy is lower semicontinuous. This assumption is needed to apply functional analytical tools, such as nonlinear semigroup theory, and is utilized to prove chain rules, for example. In this section, under Assumption \ref{asm:lsc}, which ensures that $C_c(X)$ is naturally embedded into $L^q(\Mf{m})$ for some $\Mf{m}$, we canonically extend the notions of the strongly local part and energy measures to functions that are not necessarily continuous, and prove some properties of them, for example, the chain rule of energy measures.
\begin{asm}\label{asm:lsc}
$X$ is a (nonempty) locally compact metrizable space, $\Mf{m}$ is a full-support Radon measure on $X$, and $\Mc{E}$ is a gsm on $V=\Mc{M}(\Mc{E})\subset C_c(X)$ satisfying \ref{cond:rwNC}--\ref{cond:prereg}. Moreover, $\{\Gamma^c\langle f\rangle\}_{f\in V}$ is as in Theorem \ref{thm:main1}, and $q\in[1,\infty)$ is fixed. In addition, $f\mapsto \Mc{E}(\iota^{-1}(f))$ is lower semicontinuous on $(\iota(V),\|\cdot\|_{L^q(\Mf{m})})$, where $\iota$ is the natural injection $C_c(X)\ni f\mapsto [f]\in L^q(\Mf{m})$.
\end{asm}
Hereafter, we identify $\iota(f)\in \iota (C_c(X))$ with $f$.  We also use this abbreviation for $\Mc{M}(\Mc{E}), \Mc{D}(\Mc{E})$, etc.: in other words, we identify the image of $\iota|_A$ with $A$ for $A\subset C_c(X)$. We also use the extended notion defined in Remark \ref{rem:abbrinj}.
 \begin{defi}[Lower semicontinuous envelope]
 Let $(V,\|\cdot\|_{V})$ be a seminormed vector space and $F:W\to[-\infty,\infty]$ for $W\subset V$. We define the lower semicontinuous envelope of $F$ by
 \[{}_{sc(\|\cdot\|_{V})}F(x)\coloneq \inf\left\{\liminf_{n\to\infty}F(x_n)\ \Big\vert\ \begin{minipage}{.32\linewidth}
 	$x_n\to x$ in $\|\cdot\|_{V}$ norm and $x_n\in W$ for any $n$
 \end{minipage}\right\} \]
 where we define $\inf\emptyset=\infty$. We omit $\|\cdot\|_{V}$ if no confusion may occur.
 \end{defi}
 It is evident that ${}_{sc}\Mc{E}$ is a gsm on $L^q(\Mf{m})$. The following proposition is an easy consequence of results in \cite{Cla2023} and \cite{SZ2025+}.

 \begin{prop}\label{prop:lscE}
 	\begin{enumerate}
		\item $(\Mc{M}({}_{sc}\Mc{E}), \|\cdot\|_{LN({}_{sc}\Mc{E})}+\|\cdot\|_{L^q(\Mf{m})})$ is a Banach space. Moreover, $\|\cdot\|_{LN({}_{sc}(\Mc{E}_1))}$ and $\|\cdot\|_{LN({}_{sc}\Mc{E})}+\|\cdot\|_{L^q(\Mf{m})}$ are equivalent norms on $\Mc{M}({}_{sc}\Mc{E})=\Mc{M}({}_{sc}\Mc{E}_1)$, where
		$\Mc{E}_1=\Mc{E}+\|\cdot\|^q_{L^q(\Mf{m})}$.
 		In addition, $\|\cdot\|_{LN({}_{sc}\Mc{E})}={}_{sc}(\|\cdot\|_{LN(\Mc{E})} )$ and $\|\cdot\|_{LN({}_{sc}\Mc{E}_1)}={}_{sc}(\|\cdot\|_{LN(\Mc{E}_1)} )$
 		\item If $\Mc{E}$ is reflexive,\footnote{Definition \ref{defi:modref}.} then
 		\begin{align*}
			\Mc{F}(\Mc{E})&\coloneq \left\{ f\in \Mc{M}({}_{sc}\Mc{E})\bigg\vert \begin{minipage}{.6\linewidth}
				There exists a family $\{f_n\}_{n\ge 1}$ of functions in $\Mc{M}(\Mc{E})$ such that the family is $\|\cdot\|_{LN(\Mc{E})}$-Cauchy and $f_n\to f$ in $L^q(\Mf{m})$.
 			\end{minipage} \right\}\\
 			&=\Mc{M}({}_{sc}\Mc{E}).
 		\end{align*}
 	\end{enumerate}
 \end{prop}
\begin{proof} Note that we can consider $\Mc{E}$ as a gsm on $L^q(\Mf{m})$, with $\Mc{E}(f)\coloneq \infty$ on $L^q(\Mf{m})\setminus V$: this does not affect the definitions of $\Mc{M}(\Mc{E}),\Mc{D}(\Mc{E})$ and ${}_{sc}(\Mc{E})$.
\begin{enumerate}
	\item Completeness of the normed space follows from \cite[Theorem A.13]{SZ2025+}: note that the implication (ii)$\Rightarrow$(iii)$\Rightarrow$(iv) in that theorem does not require the reflexivity of $\Mc{E}$ (cf. \cite[Theorem 2.3]{Cla2023}). Equivalence of the norms follows in the same way as \cite[Theorem 2.4]{Cla2023}. \cite[Lemma A.20]{SZ2025+} yields the last part.
	\item We obtain the second equality by \cite[Theorem A.23 (c)]{SZ2025+} with $\varrho=\Mc{E}$.
\end{enumerate}
\end{proof}
\begin{rem}\label{rem:lscincl}
	If $\{f_n\}_{n\ge 1}$ is as in the definition of $\Mc{F}(\Mc{E})$, then $\|f-f_n\|_{LN({}_{sc}\Mc{E})}\to 0$ because $\|\cdot\|_{LN({}_{sc}\Mc{E})}={}_{sc}(\|\cdot\|_{LN(\Mc{E})} )$.
\end{rem}
Now we canonically extend $\Gamma^c\langle\cdot\rangle$ onto $\Mc{F}(\Mc{E})$.
\begin{lem}\label{lem:uniquelyextension}
Let $\varrho$ be a gsm on $W$, and $\tau$ be a gsm on $Y\subset \Mc{M}(\varrho)$ with $\tau\le\varrho$ on $Y$. Then, there exists a unique gsm extension $\tilde\tau$ of $\tau$ to the $\|\cdot\|_{LN(\varrho)}$-closure of $Y$ (in $\Mc{M}(\varrho)$), which is denoted by $\tilde{Y}$, such that $\tilde\tau\le \varrho$ on $\tilde{Y}$.
\end{lem}
\begin{proof}
Since $\|\cdot\|_{LN(\tau)}\le\|\cdot\|_{LN(\varrho)}$ on $Y$, $\tau$ is lower semicontinuous with respect to $\|\cdot\|_{LN(\varrho)}$ by Lemma \ref{lem:lsc}. Thus, $\tilde{\tau}\coloneq ({}_{sc(\|\cdot\|_{LN(\varrho)})}\tau)|_{\tilde Y}$ is a gsm extension of $\tau$ on $\tilde{Y}$ with $\tilde{\tau}\le \varrho$.
Now we assume that a gsm $\hat\tau$ satisfies the same property. Let $y\in\tilde{Y}$ and $y_n\in Y\ (n\ge 1)$ with $\|y-y_n\|_{LN(\varrho)}\to 0$. By Lemmas \ref{lem:lsc} and \ref{lem:continuity},
\[\tilde{\tau}(\lambda y)=\lim_{n\to\infty}\tilde{\tau}(\lambda y_n)=\lim_{n\to\infty}\tau(\lambda y_n)=\lim_{n\to\infty}\hat{\tau}(\lambda y_n)=\hat{\tau}(\lambda y) \quad (\lambda\in(0,1))\]
holds for all $\lambda\in(0,1)$, excluding at most two values, because $\tilde{\tau},\hat{\tau}\le\varrho$. Therefore, \ref{cond:rlc} yields the uniqueness.
\end{proof}
\begin{thm}\label{thm:emextension}
Let 
\[G^c\langle f\rangle(A)\coloneq {}_{sc(\|\cdot\|_{LN({}_{sc(\|\cdot\|_{L^q(\Mf{m})})}\Mc{E})})}\Gamma^c\langle f\rangle (A)\]
for any $A\in\Mc{B}(X)$ and $f\in\Mc{F}(\Mc{E})$. Then, $\{G^c\langle f\rangle\}_{f\in\Mc{F}(\Mc{E})}$ is the unique family of finite Radon measures such that
\begin{description}
	\mcond{cond:mgsm2} The map $f\mapsto G^c\langle f\rangle(A)$ is a gsm on $\Mc{F}(\Mc{E})=\Mc{M}(f\mapsto G^c\langle f\rangle (A))$ for any $A\in\Mc{B}(X)$.
	\mcond{cond:mext} $G^c\langle f\rangle=\Gamma^c\langle f\rangle$ on $\Mc{M}(\Mc{E})$
	\mcond{cond:mlim} $G^c\langle f\rangle(A)=\lim_{n\to\infty}G^c\langle f_n\rangle(A)$ for any $A\in\Mc{B}(X)$ whenever $f_n,f\in\Mc{F}(\Mc{E})$ and $f_n\to f$ in $\|\cdot\|_{LN({}_{sc(\|\cdot\|_{L^q(\Mf{m})})}\Mc{E})}$.
	\mcond{cond:mlsc} If the map $f\mapsto \Gamma^c\langle f\rangle (A)$ is lower semicontinuous with respect to $\|\cdot\|_{L^q(\Mf{m})}$ norm for some $A\in\Mc{B}(X)$, then $G^c\langle f\rangle(A)={}_{sc(\|\cdot\|_{L^q(\Mf{m})})}\Gamma^c\langle f\rangle(A)$ on $\Mc{F}(\Mc{E})$.
\end{description}
\end{thm}
\begin{proof}
	\begin{description}
		\item[\ref{cond:mgsm2}] It is evident that $f\mapsto G^c\langle f\rangle(A)$ is a gsm. Additionally, 
		\[{}_{sc(\|\cdot\|_{L^q(\Mf{m})})}\Mc{E}(f)={}_{sc(\|\cdot\|_{LN({}_{sc(\|\cdot\|_{L^q(\Mf{m})})}\Mc{E})})}\Mc{E}(f) \]
		on $\Mc{F}(\Mc{E})$, by Lemma \ref{lem:uniquelyextension} with $\varrho={}_{sc(\|\cdot\|_{L^q(\Mf{m})})}\Mc{E}$ and $\tau=\Mc{E}$ (see also Remark \ref{rem:lscincl}.) Hence, ${}_{sc(\|\cdot\|_{L^q(\Mf{m})})}\Mc{E}(f)\ge G^c\langle f\rangle (A)$
		and \[\Mc{F}(\Mc{E})\subset \Mc{M}({}_{sc(\|\cdot\|_{L^q(\Mf{m})})}\Mc{E})\cap\Mc{F}(\Mc{E})\subset \Mc{M}(f\mapsto G^c\langle f\rangle (A)).\]
		\item[\ref{cond:mext}] Since $\|\cdot\|_{LN(\Mc{E})}\ge \|\cdot\|_{LN(f\mapsto \Gamma^c\langle f\rangle (A))}$ on $\Mc{M}(\Mc{E})$ for any $A\in\Mc{B}(X)$, $f\mapsto \Gamma^c\langle f\rangle (A)$ is lower semicontinuous with respect to $\|\cdot\|_{LN({}_{sc(\|\cdot\|_{L^q(\Mf{m})})}\Mc{E})}$ by Lemma \ref{lem:lsc}, which leads to $G^c\langle f\rangle=\Gamma^c\langle f\rangle$ on $\Mc{M}(\Mc{E})$.
		\item[\ref{cond:mlim}] Let $f\in\Mc{F}(\Mc{E})$ and choose $f_n\in\Mc{M}(\Mc{E})$ as in the definition of $\Mc{F}(\Mc{E})$. By Lemma \ref{lem:continuity} and Remark \ref{rem:lscincl}, ${}_{sc}\Mc{E}(2(f-f_n))\to0$. Convexity, \ref{cond:mext}, and the bound $G^c\langle\cdot\rangle(X)\le{}_{sc}\Mc{E}$ give
\[
G^c\langle f\rangle(X)
\le \frac{\Gamma^c\langle2f_n\rangle(X)
+{}_{sc}\Mc{E}(2(f-f_n))}{2}<\infty
\]
for sufficiently large $n$, since $\Gamma^c\langle2f_n\rangle$ is finite by Theorem \ref{thm:main1}.
By the monotonicity of lower semicontinuous envelopes,
$0\le G^c\langle h\rangle(A)\le G^c\langle h\rangle(X)<\infty$ for any $A\in\Mc{B}(X)$ and $h\in\Mc{F}(\Mc{E})$. Since $\Mc{F}(\Mc{E})$ is linear, $\Mc{D}^i(h\mapsto G^c\langle h\rangle(A))=\Mc{F}(\Mc{E})$. Together with $\|\cdot\|_{LN(h\mapsto G^c\langle h\rangle(A))}\le\|\cdot\|_{LN({}_{sc}\Mc{E})}$, Lemma \ref{lem:limeval} yields \ref{cond:mlim}.
		\item[\ref{cond:mlsc}] This follows from Lemma \ref{lem:uniquelyextension}.
	\end{description}
	Here, \ref{cond:mext} and \ref{cond:mlim} ensure the uniqueness on $\Mc{F}(\Mc{E})$. We next show that $G^c\langle f\rangle$ is a Borel measure. We first show $G^c\langle f\rangle$ is inner regular on open sets. Let $U$ be open and let
	 \begin{align*}
	 	G'\langle f\rangle(U)
	 	=\sup\{G^c\langle f\rangle(V)\mid V\Subset U,\ V\text{:open}.\}
	 \end{align*}
	 Then, it is easily shown that $f\mapsto G'\langle f\rangle (U)$ is a gsm on $\Mc{F}(\Mc{E})$ with $G'\langle f\rangle (U)\le{}_{sc}\Mc{E}(f)$. Since $\Gamma^c\langle f\rangle$ is Radon for $f\in\Mc{M}(\Mc{E})$, we have $G'\langle f\rangle(U)=G^c\langle f\rangle(U)$ on $\Mc{M}(\Mc{E})$. Thus, $G'\langle f\rangle(U)=G^c\langle f\rangle(U)$ on $\Mc{F}(\Mc{E})$ by Lemma \ref{lem:uniquelyextension}. Hence, $G^c\langle f\rangle$ is inner regular on open sets. Similarly, for $A\in\Mc{B}(X)$,
	 \begin{align*}
	 	G''\langle f\rangle(A)
	 	=\inf\{G^c\langle f\rangle(U)\mid A\subset U,\ U\text{:open}.\}
	 \end{align*}
	 is a gsm on $\Mc{F}(\Mc{E})$ (by Lemma \ref{lem:setseq}), and $G''\langle f\rangle=G^c\langle f\rangle$ is outer regular for any $f\in \Mc{F}(\Mc{E})$, by Lemma \ref{lem:uniquelyextension}.
	 
	 Now fix $f\in\Mc{F}(\Mc{E})$ and choose $f_n\in\Mc{M}(\Mc{E})$ as in the definition of $\Mc{F}(\Mc{E})$. By \ref{cond:mext} and \ref{cond:mlim}, for any open sets $U,V$,
	 \begin{multline*}
G^c\langle f\rangle(U\cup V)=\lim_{n\to\infty}\Gamma^c\langle f_n\rangle(U\cup V)\\
\le\lim_{n\to\infty}\bigl(\Gamma^c\langle f_n\rangle(U)+\Gamma^c\langle f_n\rangle(V)\bigr)
=G^c\langle f\rangle(U)+G^c\langle f\rangle(V).
	 \end{multline*}
	 Equality holds when $U\cap V=\emptyset$. Therefore \cite[Theorems 5.1, 5.6 and Proposition 5.5(a)]{DeG1977} (cf. Theorem \ref{thm:Radon}) yields that $G^c\langle f\rangle$ is a Borel measure.

Set $S=\bigcup_{n\ge1}\supp(f_n)$.
By Remark \ref{rem:lscincl}, \ref{cond:mlim},
\ref{cond:mext} and \ref{cond:mlocal},
\[
G^c\langle f\rangle(X\setminus S)
=\lim_{n\to\infty}
\Gamma^c\langle f_n\rangle(X\setminus S)=0.
\]
The restriction of $G^c\langle f\rangle$ to each compact set
$\bigcup_{k=1}^n\supp(f_k)$ is a finite Borel measure
on a compact metrizable space and hence is Radon.
Since these sets increase to $S$,
$G^c\langle f\rangle$ is a finite Radon measure on $X$.
\end{proof}
Hereafter, we write $\Gamma^c\langle f\rangle=G^c\langle f\rangle$ and $\Mc{E}^{(c)}(f)=G^c\langle f\rangle (X)$ also for $f\in\Mc{F}(\Mc{E})$, to ease notation.


For general $A\in\Mc{B}(X)$, we do not know whether the assumption of \ref{cond:mlsc} holds. However, if $\Mc{E}^{(c)}$ is lower semicontinuous, then that holds for any open $A$.
 
\begin{lem}\label{lem:lscopen}
If $\Mc{E}^{(c)}$ is lower semicontinuous on $(\Mc{M}(\Mc{E}), \|\cdot\|_{L^q})$, then for any open set $U$, the function $f\mapsto \Gamma^c\langle f\rangle (U)$ is also lower semicontinuous on $(\Mc{M}(\Mc{E}), \|\cdot\|_{L^q})$.
\end{lem}
\begin{proof}
Let $f_n,f\in \Mc{M}(\Mc{E})\ (n\in\Mb{N})$ with $f_n\to f$ in $L^q(\Mf{m})$. The claim is immediate if $\liminf_{n\to\infty}\Gamma^c\langle f_n\rangle(U)=\infty$. Otherwise, taking a subsequence if necessary, we may assume $\lim_{n\to\infty}\Gamma^c\langle f_n\rangle (U)<\infty$ without loss of generality. Moreover, taking a subsequence again, if necessary, we also assume that $f_n\to f\ \Mf{m}$-a.e. Then, for any $m\in\Mb{N}, g\in\Mc{M}(\Mc{E})$ and $a\in\Mb{R}$, we have $(f_n)_m^{g,a}\to f_m^{g,a}\ (n\to\infty)\ \Mf{m}$-a.e. In particular, if $a<0$ then $0\le (f_n)_m^{g,a}\le S_m^a\circ g$ and $\|S_m^a\circ g\|_{L^q(\Mf{m})}\le 2^{-m}\Mf{m}(\supp(g))^\frac{1}{q}$ for sufficiently large $m$, and hence $(f_n)_m^{g,a}\to f_m^{g,a}\ (n\to\infty)$ in $L^q(\Mf{m})$, by the dominated convergence theorem.

Let $g\in\Mc{M}(\Mc{E}), a<0$ such that $\{g\le a\}\subset U$. Then, by the lower semicontinuity of $\Mc{E}^{(c)}$ we obtain
\begin{align*}
	&\Mc{E}^{(c)}(f_m^{g,a})\\
	\le &\liminf_{n\to\infty} \Mc{E}^{(c)}((f_n)_m^{g,a})\\
	\overset{\eqref{eq:fundestm}}{\le}&\liminf_{n\to\infty}\left(\Gamma^c\langle f_n\rangle(\{g< a+2^{1-m}\})+\Gamma^c\langle g\rangle(\{a+2^{-m}<g<a+2^{1-m}\})   \right)\\
	\le &\lim_{n\to\infty}\Gamma^c\langle f_n\rangle(U)+\Gamma^c\langle g\rangle(\{a<g<a+2^{1-m}\}).
\end{align*}
Since $\lim_{m\to\infty}\Gamma^c\langle g\rangle(\{a<g<a+2^{1-m}\})=0$, letting $m\to\infty$ gives $\Gamma^c\langle f\rangle (\{g\le a\})\le \lim_{n\to\infty}\Gamma^c\langle f_n\rangle(U)$. Inner regularity of $\Gamma^c\langle f\rangle$ then yields $\Gamma^c\langle f\rangle(U)\le \lim_{n\to\infty}\Gamma^c\langle f_n\rangle(U)$.
\end{proof}

In the rest of this section, we prove the chain rule for energy measures: the basic strategy is the same as that of \cite[Theorem 7.6]{Sas2026}.
\begin{lem}\label{lem:NC}
	If $f\in\Mc{F}(\Mc{E})$, $\varphi\in C(\Mb{R})$ with $\varphi(0)=0$ and $\varphi$ is 1-Lipschitz, then $\varphi\circ f\in \Mc{M}({}_{sc}\Mc{E})$ and ${}_{sc}\Mc{E}(\varphi\circ f)\le {}_{sc}\Mc{E}(f)$. In particular, if $\Mc{E}^{(c)}$ is lower semicontinuous on $(\Mc{M}(\Mc{E}), \|\cdot\|_{L^q})$, then ${}_{sc}\Gamma^c\langle \varphi\circ f \rangle (U) \le {}_{sc}\Gamma^c\langle f \rangle (U)$ for open $U$.
\end{lem}
\begin{proof}
	By \ref{cond:rlc}, the first assertion reduces to the case $f\in\Mc{D}^i({}_{sc}\Mc{E})$ and follows as in \cite[Proposition 4.5]{Sas2026}. By \ref{lem:lscopen}, the second assertion follows by the same approximation argument.
\end{proof}
\begin{defi}
Let $\{x_n\}_{n\ge 1}$ be a sequence of elements in a vector space $W$. We call $\{y_n\}_{n\ge 1}$ a \emph{Mazur subsequence} of $\{x_n\}_{n\ge 1}$ if $y_n$ is a (finite) convex combination of $\{x_k\}_{k\ge n}$ for each $n$. 
\end{defi}
Note that a Mazur subsequence of a subsequence of $\{x_n\}_{n\ge 1}$ is itself a Mazur subsequence of $\{x_n\}_{n\ge 1}$.
\begin{lem}[{cf. \cite[Lemma 2.10]{EM2025+}}]\label{lem:Mazsub}
Assume that $\Mc{E}$ is reflexive and let $\{f_n\}$ be a sequence in $\Mc{M}(\Mc{E})$. If $f_n\to f$ in $L^q(\Mf{m})$ and $\sup_{n\ge 1}\|f_n\|_{LN(\Mc{E})}<\infty$, then $f\in\Mc{F}(\Mc{E})$ and $\Gamma^c\langle f\rangle (A)\le \liminf_{n\to\infty}\Gamma^c\langle f_n\rangle(A)$ for any $A\in\Mc{B}(X)$. Moreover, there exists a Mazur subsequence $\{g_n\}_{n\ge 1}$ of $\{f_n\}_{n\ge 1}$ such that $g_n\to f$ in $\|\cdot\|_{LN({}_{sc}\Mc{E})}+\|\cdot\|_{L^q(\Mf{m})}$-norm.
\end{lem}
\begin{proof}
Fix $A\in\Mc{B}(X)$ and take a subsequence $\{f_{n_k}\}_{k\ge1}$ such that
\[\lim_{k\to\infty}\Gamma^c\langle f_{n_k}\rangle(A)=\liminf_{n\to\infty}\Gamma^c\langle f_n\rangle(A).\]
By the assumed $LN(\Mc{E})$-boundedness, reflexivity, \cite[Theorem 2.4.1]{HKST2015} and Mazur's lemma, there exists a $\|\cdot\|_{LN(\Mc{E})}$-Cauchy Mazur subsequence $\{g_n\}_{n\ge1}$ of $\{f_{n_k}\}_{k\ge1}$. Since $g_n\to f$ in $L^q(\Mf{m})$, the sequence is Cauchy in $\|\cdot\|_{LN({}_{sc}\Mc{E})}+\|\cdot\|_{L^q(\Mf{m})}$-norm. By Proposition \ref{prop:lscE}~(1), it converges in this norm to some $g\in\Mc{M}({}_{sc}\Mc{E})$. Thus $g=f\in\Mc{F}(\Mc{E})$, by the definition of $\Mc{F}(\Mc{E})$. Moreover, $\{g_n\}_{n\ge1}$ is a Mazur subsequence of the original sequence. Hence,
\begin{align*}
\Gamma^c\langle f\rangle(A)&=\lim_{n\to\infty}\Gamma^c\langle g_n\rangle(A)\\
&\le \lim_{k\to\infty}\Gamma^c\langle f_{n_k}\rangle(A)
=\liminf_{n\to\infty}\Gamma^c\langle f_n\rangle(A)
\end{align*}
by \ref{cond:mlim} and \ref{cond:rconv}, since $g_n$ is a convex combination of $\{f_{n_k}\}_{k\ge n}$.
\end{proof}
\begin{cor}\label{cor:Mazsub}
Assume that $\Mc{E}$ is reflexive. Let $\varphi:\Mb{R}\to\Mb{R}$ be a Lipschitz function with $\varphi(0)=0$, and $f\in\Mc{M}(\Mc{E})$. Then, there exists a sequence $\{\varphi_n\}_{n\ge 1}$ of piecewise affine functions with the following properties.
\begin{itemize}
	\item $\varphi_n\to\varphi$ pointwise on $\Mb{R}$ and $\varphi_n'\to\varphi'$ $\Mc{L}^1$-a.e.~on $\Mb{R}$. In particular, if $\varphi$ is piecewise $C^1$, then  $\{\varphi_n\}_{n\ge 1}$ can be chosen so that $\varphi_n'(t)\not\to\varphi'(t)$ only on an at most countable set.
	\item $\varphi_n(0)=0$ and $\Lip(\varphi_n)\le \Lip(\varphi)$ for any $n$.
	\item $\varphi_n\circ f\to\varphi\circ f$ in $\|\cdot\|_{LN({}_{sc}\Mc{E})}+\|\cdot\|_{L^q(\Mf{m})}$-norm.
\end{itemize}
\end{cor}
\begin{proof}
Let $L=\Lip(\varphi)$. By \cite[(4.5) in Theorem 4.1.4]{GP2020}, there exists a sequence $\{\psi_n\}_{n\ge1}$ of piecewise affine functions with the first two properties, with $\varphi_n$ replaced by $\psi_n$.\footnote{The assertion for piecewise $C^1$ functions follows from the construction.} Then, $\psi_n\circ f\to\varphi\circ f$ pointwise, $|\psi_n\circ f|\le L|f|$ and $\|Lf\|_{L^q(\Mf{m})}<\infty$. Thus, $\psi_n\circ f\to\varphi\circ f$ in $L^q(\Mf{m})$ by the dominated convergence theorem. Lemma \ref{lem:LNcontraction} yields
\[\sup_{n\ge1}\|\psi_n\circ f\|_{LN(\Mc{E})}\le L\|f\|_{LN(\Mc{E})}<\infty.\]
Applying Lemma \ref{lem:Mazsub} to $\{\psi_n\circ f\}_{n\ge1}$ and using the same convex coefficients for $\{\psi_n\}_{n\ge1}$, we obtain a Mazur subsequence $\{\varphi_n\}_{n\ge1}$ such that $\varphi_n\circ f\to\varphi\circ f$ in $\|\cdot\|_{LN({}_{sc}\Mc{E})}+\|\cdot\|_{L^q(\Mf{m})}$-norm. The Mazur subsequence inherits the first two properties.
\end{proof}
\begin{lem}\label{lem:eidp}
Assume that $\Mc{E}$ is reflexive. Then $f_*\Gamma^c\langle f\rangle \ll \Mc{L}^1|_{\Mc{B}(\Mb{R})}$ for any $f\in\Mc{M}(\Mc{E})$.
\end{lem}
\begin{proof}
Fix $f\in\Mc{M}(\Mc{E})$. Since $f_*\Gamma^c\langle f\rangle$ is a finite Borel measure on $\Mb{R}$ and hence inner regular, it is sufficient to show $f_*\Gamma^c\langle f\rangle(K)=0$ for any Lebesgue-null compact set $K$. Let $M>0$ be such that $K\subset[-M,M]$. Note that the $n^{-1}$-neighborhood of $K$, denoted by $K_n$, is a union of at most $1+\lceil nM\rceil_1$ intervals. Let $\varphi_n(t)\coloneq \int_0^t \indi_{K_n}(s)\der s$. Then, $|\varphi_n\circ f|\le \Mc{L}^1(K_n)\to 0$ uniformly and $|\varphi_n\circ f|\le |f|$, hence $\varphi_n\circ f\to 0$ in $L^q(\Mf{m})$ by the dominated convergence theorem. Moreover, Lemma \ref{lem:LNcontraction} applied to $t\mapsto t-\varphi_n(t)$ gives
\[\sup_n\|f-\varphi_n\circ f\|_{LN(\Mc{E})}\le\|f\|_{LN(\Mc{E})}<\infty.\]
Since $f_*\Gamma^c\langle f\rangle(\Mb{R})=\Mc{E}^{(c)}(f)<\infty$, we obtain
\begin{multline*}
f_*\Gamma^c\langle f\rangle (K)=\lim_{n\to\infty}f_*\Gamma^c\langle f\rangle (K_n)\\=\Mc{E}^{(c)}(f)-\lim_{n\to\infty}\Mc{E}^{(c)}(f-\varphi_n\circ f)\le \Mc{E}^{(c)}(f)-\Mc{E}^{(c)}(f)=0
\end{multline*}
by \ref{cond:mlocal} and Lemma \ref{lem:Mazsub}.

\end{proof}

\begin{rem}
This lemma is an easy version of the \emph{energy image density property}, introduced in \cite{BH1986}. In the $p$-homogeneous case, Eriksson-Bique and Murugan proved its multi-variable version for general functions (and also gave a positive answer to the long-standing conjecture of Bouleau and Hirsch): see \cite{EM2025+}.
\end{rem}

\begin{thm}[Chain rule for 1 variable]\label{thm:chain1}
Assume that $\Mc{E}$ is reflexive, and that $E$ is a gsm with \ref{cond:rconv}--\ref{cond:rstrsubad} such that $\Mc{M}(\Mc{E})=\Mc{M}(E)$. Let $\Phi:\Mb{R}\to[0,\infty)$ be continuous.
\begin{enumerate}
	\item If $\Mc{E}(tf)\le \Phi(t)E(f)$ for any $f\in\Mc{M}(\Mc{E})$ and $t\in\Mb{R}$, then
	\begin{equation}
\Gamma_\Mc{E}^c\langle\varphi\circ f\rangle
\le \bigl(\Phi\circ(\varphi')\circ f\bigr)\der\Gamma_E^c\langle f\rangle\notag
	\end{equation}
	for any $f\in\Mc{M}(\Mc{E})$ and any piecewise $C^1$ function $\varphi:\Mb{R}\to\Mb{R}$ with $\varphi(0)=0$. Moreover, if $E$ is lower semicontinuous with respect to $\|\cdot\|_{L^q(\Mf{m})}$ norm and reflexive, then the same inequality holds for any locally Lipschitz $\varphi$ with $\varphi(0)=0$.
	\item If $\Mc{E}(tf)\ge \Phi(t)E(f)$ for any $f\in\Mc{M}(\Mc{E})$ and $t\in\Mb{R}$, then
	\begin{equation}
	\bigl(\Phi\circ(\varphi')\circ f\bigr)\der\Gamma_E^c\langle f\rangle
	\le \Gamma_\Mc{E}^c\langle\varphi\circ f\rangle\notag 
	\end{equation}
	for any $f\in\Mc{M}(\Mc{E})$ and any locally Lipschitz function $\varphi:\Mb{R}\to\Mb{R}$ with $\varphi(0)=0$.
\end{enumerate}
\end{thm}
\begin{rem}\label{rem:measwelldef}
Let $M>\|f\|_\infty$. Set $\tilde\varphi(t)\coloneq \varphi((-M)\vee t\wedge M)$. Then $\varphi\circ f=\tilde\varphi\circ f$ and $\varphi'=\tilde\varphi'$ on $[-\|f\|_{\infty},\|f\|_{\infty}]$, so we may assume that $\varphi$ is Lipschitz.

In (1), by \ref{cond:mlocal}, $\Gamma_E^c\langle f\rangle (f^{-1}(\{\varphi'\text{ does not exist}\}))=0$ if $\varphi$ is a piecewise $C^1$ function. This also holds for locally Lipschitz $\varphi$ under the additional assumptions on $E$, by Lemma \ref{lem:eidp}. Hence, $\bigl(\Phi\circ(\varphi')\circ f\bigr)\der\Gamma_E^c\langle f\rangle$ is well-defined. For convenience, we set $\varphi'=0$ where the derivative does not exist.

In (2), if $\Phi(t_0)\ne 0$ for some $t_0\in\Mb{R}$, then Lemmas \ref{lem:dom} and \ref{lem:eidp} yield \[f_*\Gamma_E^c\langle f\rangle \le \Phi(t_0)^{-1} f_*\Gamma_\Mc{E}^c\langle t_0 f\rangle \ll  \Mc{L}^1,\]
so $\bigl(\Phi\circ(\varphi')\circ f\bigr)\der\Gamma_E^c\langle f\rangle$ is well-defined, similarly to (1). Otherwise, $\Phi\equiv0$ and $\bigl(\Phi\circ(\varphi')\circ f\bigr)\der\Gamma_E^c\langle f\rangle=0$.
\end{rem}
Typically, $\Mc{E}=E$ and $\Phi(t)=C(|t|^p\vee |t|^{p'})$ with some $p,p'>1$ and $C>0$ fulfill the assumption of this theorem.

\begin{proof}[Proof of Theorem \ref{thm:chain1}]
(1)\ We first assume that $\varphi$ is piecewise affine, and write
\[\varphi\circ f=\sum_{i=1}^k \alpha_i\eta_{[a_{i-1},a_i]}\circ f.\]
Then by \ref{cond:mlocal} and Lemma \ref{lem:dom},
\begin{multline*}
\Gamma_\Mc{E}^c\langle \varphi\circ f\rangle=\sum_{i=1}^k \Gamma_\Mc{E}^c\langle \alpha_i \eta_{[a_{i-1},a_i]}\circ f\rangle\\
\le\sum_{i=1}^k \Phi(\alpha_i)\Gamma_E^c\langle \eta_{[a_{i-1},a_i]}\circ f\rangle
=\sum_{i=1}^k \Phi(\alpha_i)\Gamma_E^c\langle f\rangle|_{f^{-1}([a_{i-1},a_i))},
\end{multline*}
and thus the claim follows. 

For general $\varphi$, let $\varphi_n$ be as in Corollary \ref{cor:Mazsub}. By \ref{cond:mlim},
\[\lim_{n\to\infty}\Gamma_\Mc{E}^c\langle\varphi_n\circ f\rangle(A)
=\Gamma_\Mc{E}^c\langle\varphi\circ f\rangle(A)\qquad(A\in\Mc{B}(X)).\]
Since $\sup_{|s|\le \Lip(\varphi)}\Phi(s)<\infty$, $\Gamma_E^c\langle f\rangle(X)<\infty$ and $(\varphi_n')\circ f\to(\varphi')\circ f$ $\Gamma_E^c\langle f\rangle$-a.e.~(see Remark \ref{rem:measwelldef}), the dominated convergence theorem yields
\begin{multline*}
	\Gamma_\Mc{E}^c\langle\varphi\circ f\rangle(A)
=\lim_{n\to\infty}\Gamma_\Mc{E}^c\langle\varphi_n\circ f\rangle(A)\\
	\le\lim_{n\to\infty}\int_A\Phi\circ(\varphi_n')\circ f\der\Gamma_E^c\langle f\rangle
	=\int_A\Phi\circ(\varphi')\circ f\der\Gamma_E^c\langle f\rangle
\end{multline*}
for every $A\in\Mc{B}(X)$, which proves (1).

(2)\ The proof of (2) is similar, with the comparison inequalities reversed.
\end{proof}
\begin{proof}[{Proof of Theorem \ref{thm:main2}}]
(In this proof alone, we do not use the notation of Assumption \ref{asm:lsc}.) Let $f\in\Mc{D}(\Mc{E})$ and choose $f_n\in\Mc{C}$ converging to $f$ in $\|\cdot\|_{LN(\Mc{E})}+\|\cdot\|_{L^q(\Mf{m})}$. For $0<\lambda<1$, convexity gives
\[
\Mc{E}(\lambda f_n)\le\lambda\Mc{E}(f)
+(1-\lambda)\Mc{E}\left(\frac{\lambda}{1-\lambda}(f_n-f)\right).
\]
By Lemma \ref{lem:continuity}, the last term tends to zero. Taking the $L^q$-lower semicontinuous envelope and then letting $\lambda\uparrow1$, we obtain
\[
{}_{sc(\|\cdot\|_{L^q(\Mf{m})})}(\Mc{E}|_{\Mc{C}})=\Mc{E},
\]
where the reverse inequality follows from the lower semicontinuity of $\Mc{E}$. The assumed regularity also gives $\Mc{F}(\Mc{E}|_{\Mc{C}})=\Mc{M}(\Mc{E})$.
Applying Theorems \ref{thm:main1} and \ref{thm:emextension} to $\Mc{E}|_{\Mc{C}}$, we obtain the core properties in (1) and the finite Radon measures in (2). Set $\Mc{E}^{(c)}=\infty$ on $L^q(\Mf{m})\setminus\Mc{M}(\Mc{E})$; then $\Mc{M}(\Mc{E}^{(c)})=\Mc{M}(\Mc{E})$. Maximality in (1) follows from Theorem \ref{thm:main1} on $\Mc{C}$, Lemma \ref{lem:uniquelyextension} and the monotonicity of lower semicontinuous envelopes. Uniqueness in (2) follows from uniqueness on $\Mc{C}$ and \ref{cond:mconti}. Part (3) follows from the corresponding core properties and uniqueness, and (4) follows from Theorem \ref{thm:chain1}, with the composed functions in $\Mc{F}(\Mc{E}|_{\Mc{C}})=\Mc{M}(\Mc{E})$.

Assume now that $\Mc{E}$ is reflexive. Its restriction to $\Mc{C}$ is reflexive by the assumed regularity of $\Mc{C}$ in $LN(\Mc{E})$. Let $f,g\in\Mc{M}(\Mc{E})$, and choose $f_n,g_n\in\Mc{C}$ converging respectively to $f,g$ in $\|\cdot\|_{LN(\Mc{E})}+\|\cdot\|_{L^q(\Mf{m})}$-norm. By \ref{cond:mlim},
\[
\Mc{E}^{(c)}(f_n)\to\Mc{E}^{(c)}(f),\qquad
\Mc{E}^{(c)}(g_n)\to\Mc{E}^{(c)}(g).
\]
By Lemma \ref{lem:LNcontraction}, the triangle inequality and $f_n\vee g_n,f_n\wedge g_n=(f_n+g_n\pm|f_n-g_n|)/2$, we have
\[
\|f_n\vee g_n\|_{LN(\Mc{E})}\vee\|f_n\wedge g_n\|_{LN(\Mc{E})}
\le\|f_n\|_{LN(\Mc{E})}+\|g_n\|_{LN(\Mc{E})}.
\]
Thus both sequences are bounded in $LN(\Mc{E})$ and converge respectively to $f\vee g,f\wedge g$ in $L^q(\Mf{m})$. Applying Lemma \ref{lem:Mazsub} for $\Mc{E}|_{\Mc{C}}$ with $A=X$ to these sequences gives $f\vee g,f\wedge g\in\Mc{M}(\Mc{E})$. Using \ref{cond:rstrsubad} on $\Mc{C}$, we obtain
\[
\begin{aligned}
\Mc{E}^{(c)}(f\vee g)+\Mc{E}^{(c)}(f\wedge g)
&\le\liminf_{n\to\infty}\{\Mc{E}^{(c)}(f_n\vee g_n)+\Mc{E}^{(c)}(f_n\wedge g_n)\}\\
&\le\lim_{n\to\infty}\{\Mc{E}^{(c)}(f_n)+\Mc{E}^{(c)}(g_n)\}\\
&=\Mc{E}^{(c)}(f)+\Mc{E}^{(c)}(g).
\end{aligned}
\]
For a piecewise affine $\varphi$ as in \ref{cond:rwNC}, Lemma \ref{lem:LNcontraction} gives
\[
\|\varphi\circ f_n\|_{LN(\Mc{E})}\le\|f_n\|_{LN(\Mc{E})}.
\]
Since $\varphi\circ f_n\to\varphi\circ f$ in $L^q(\Mf{m})$, Lemma \ref{lem:Mazsub} for $\Mc{E}|_{\Mc{C}}$ with $A=X$ gives $\varphi\circ f\in\Mc{M}(\Mc{E})$ and, together with \ref{cond:rwNC} on $\Mc{C}$, yields
\[
\Mc{E}^{(c)}(\varphi\circ f)
\le\liminf_{n\to\infty}\Mc{E}^{(c)}(\varphi\circ f_n)
\le\lim_{n\to\infty}\Mc{E}^{(c)}(f_n)
=\Mc{E}^{(c)}(f).
\]
\end{proof}

\section{Clarkson-type and Ball--Carlen--Lieb-type inequalities and detailed chain rules}\label{sec:CI}
To discuss reflexivity and stronger properties of $\Mc{E}$, in this section we introduce Clarkson-type lower bounds and Ball--Carlen--Lieb-type upper bounds for gsms (see Remark \ref{rem:clk}~(4)--(5)). In particular, we prove more detailed chain rules than those in Theorem \ref{thm:chain1} with these inequalities. Throughout this section, we continue to assume Assumption \ref{asm:lsc}.
\begin{defi}\label{defi:rCLK}
Let $\varrho$ be a gsm on a vector space $W$.
\begin{enumerate}
	\item We say $\varrho$ satisfies a \emph{Clarkson-type lower bound} with indices $\alpha,\beta\in(0,\infty)$ if
	\begin{equation}\label{eq:lCLK}
		2(\varrho(u)^\alpha+\varrho(v)^\alpha)^{1/\alpha}\le \varrho(u+v)+\varrho(u-v)\ \text{and}\   \varrho(tu)\le t^\beta \varrho(u) \tag{$\lCla(\alpha,\beta)$}
	\end{equation}
	for any $t\in[1,\infty)$ and $u,v\in\Mc{M}(\varrho)$.
	\item We say $\varrho$ satisfies a \emph{Ball--Carlen--Lieb-type upper bound} with indices $\alpha,\beta\in(0,\infty)$ with $\alpha\beta>1$ if
	\begin{equation}\label{eq:uCLK}
			2(\varrho(u)^\alpha+c\varrho(v)^\alpha)^{1/\alpha}\ge \varrho(u+v)+\varrho(u-v)\ \text{and}\  \varrho(tu)\le Ct^\beta \varrho(u) \tag{$\uBCL(\alpha,\beta)$}
	\end{equation}
	for some $c,C\ge1$ and any $u,v\in\Mc{M}(\varrho)$ and $t\in[0,1]$.
\end{enumerate}
\end{defi}
\begin{rem}\label{rem:clk}
\begin{enumerate}
	\item Since $\alpha\mapsto(a^\alpha+cb^\alpha)^{1/\alpha}$ is nonincreasing on $(0,\infty)$ for fixed $a,b\ge0$ and $c\ge1$, ($\lCla(\alpha_1,\beta)$) implies ($\lCla(\alpha_2,\beta)$) for $\alpha_2>\alpha_1$, and ($\uBCL(\alpha_1,\beta)$) implies ($\uBCL(\alpha_2,\beta)$) for $\alpha_2\in (\beta^{-1},\alpha_1)$ with the same $c,C$. We usually consider ($\lCla(\alpha,\beta)$) for $\alpha\in[1,\infty)$ and ($\uBCL(\alpha,\beta)$) for $\alpha\in(0,1]$.
	\item If \eqref{eq:lCLK} holds and $\varrho(u)>0$ for some $u\in\Mc{M}(\varrho)$, then $1+\frac{1}{\alpha}\le \beta$, which forces $\alpha\beta>1$.
	\item Both \eqref{eq:lCLK} and \eqref{eq:uCLK} imply the following $\Delta_2$ condition:
	\begin{equation}\label{eq:delta2}
		\varrho(2u)\le C\varrho(u) \quad\text{for some $C>1$ and any $u\in\Mc{M}(\varrho)$}\tag{$\Delta_2$}
	\end{equation}
	In particular, \eqref{eq:uCLK} with $v=u$ gives $\varrho(2u)\le2(1+c)^{1/\alpha}\varrho(u)$. The $\Delta_2$ condition is one of the standard assumptions on modular functions and gsms. We put some standard and easy consequences of \eqref{eq:delta2} in Lemma \ref{lem:delta2} without proof, for the reader's convenience.
	\item For a $p$-homogeneous $\varrho$ with $p\in(1,\infty)$, the \emph{$p$-Clarkson inequalities} mean $\lCla\left(\frac{1}{p-1}\vee 1,p\right)$ and $\uBCL\left(\frac{1}{p-1}\wedge 1,p\right)$ with $c=1$. Clarkson~\cite{Cla1936} proved that $\|\cdot\|^p_{L^p(\Mf{m})}$ satisfies the $p$-Clarkson inequalities for any $L^p$ space. In this example, the ambient exponent is $q=p$. These inequalities are used in recent studies of $p$-energy forms (see \cite{KS2024+,Sas2026}, for example).
	\item The upper bound \eqref{eq:uCLK} is motivated by the inequality of Ball and Pisier for $L^p$ spaces, reported as an unpublished result in \cite[Theorem 1]{BCL1994}, where Ball, Carlen and Lieb proved its Schatten-norm counterpart. These inequalities yield sharp bounds for the moduli of convexity and smoothness. In view of Lemma \ref{lem:boundinherit}, we adopt this generalized form only for the upper bound.
\end{enumerate}
\end{rem}
\begin{lem}\label{lem:delta2}
Let $\varrho$ be a gsm satisfying \eqref{eq:delta2}. Then the following holds.
\begin{enumerate}
	\item $\Mc{M}(\varrho)=\Mc{D}(\varrho)=\Mc{D}^i(\varrho)$.
	\item $t\mapsto \varrho(tu)$ is continuous on $\Mb{R}$ for any $u\in\Mc{M}(\varrho)$.
	\item $\varrho(u)=0$ (if and) only if $\|u\|_{LN(\varrho)}=0$.
	\item $\varrho(u)=1$ if and only if $\|u\|_{LN(\varrho)}=1$.
	\item $C\ge2$ for the constant $C$ in \eqref{eq:delta2} if $\varrho(u)>0$ for some $u\in\Mc{M}(\varrho)$.
	\item There exist $C'\ge1$ and $\gamma\ge \frac{\log C}{\log 2}$ depending only on $C$ such that $\varrho(tu)\le \left(C't^\frac{\log C}{\log 2}\wedge t^\gamma\right) \varrho(u)$ for any $t\ge 1$.
\end{enumerate}
\end{lem}
\begin{lem}\label{lem:midestm}
If $\varrho$ satisfies \eqref{eq:uCLK} then 
\begin{equation}\label{eq:midestm}
	0\le \varrho(u+v)+\varrho(u-v)-2\varrho(u)\le C_{(\varrho(u))}\|v\|^{\alpha \beta}_{LN(\varrho)}
\end{equation}
for any $u,v\in\Mc{M}(\varrho)$ with $\varrho(u)\ne 0$ and $\varrho(v)\le 1$, where $C_{(\varrho(u))}$ is the constant depending only on constants in \eqref{eq:uCLK} and $\varrho(u)$. Moreover, if $\alpha\le 1$ then the same estimate also holds for $u$ with $\varrho(u)=0$, and $C_{(\varrho(u))}$ can be chosen so that $t\mapsto C_{(t)}$ is nondecreasing on $[0,\infty)$.
\end{lem}
\begin{proof}
The first inequality follows from convexity. By \eqref{eq:uCLK}, we also obtain
\begin{align*}
	\varrho(u+v)+\varrho(u-v)-2\varrho(u)&\le 2\left((\varrho(u)^\alpha+c\varrho(v)^\alpha )^\frac{1}{\alpha}-\varrho(u)\right)\\
	&\le 2\left((\varrho(u)^\alpha+cC^\alpha\|v\|_{LN(\varrho)}^{\alpha \beta} )^\frac{1}{\alpha}-\varrho(u)\right)\\
	&\le C_{(\varrho(u))}\|v\|^{\alpha \beta}_{LN(\varrho)},
\end{align*}
by the mean value theorem, where we may take
\[
C_{(s)}=\begin{cases}
\dfrac{2cC^\alpha}{\alpha}(s^\alpha+cC^\alpha)^{1/\alpha-1},&0<\alpha\le1,\ s\ge0,\\
\dfrac{2cC^\alpha}{\alpha}s^{1-\alpha},&\alpha>1,\ s>0.
\end{cases}
\]
The first expression is nondecreasing in $s\ge0$.
\end{proof}
\begin{cor}\label{cor:diff}
Assume that $\varrho$ satisfies \eqref{eq:uCLK}.
\begin{enumerate}
	\item \[D_\varrho(u)[v]\coloneq \frac{\der}{\der t}\varrho(u+tv)\big|_{t=0} \]
	is well-defined for $u,v\in\Mc{M}(\varrho)$ whenever $\varrho(u)\ne 0$ or $\alpha\le 1$, and $D_\varrho(u):v\mapsto D_\varrho(u)[v]$ is a bounded linear operator for each $u$.
	\item For such $u$ and $v$ with $\varrho(v)\le 1$ we have
	\begin{equation}\label{eq:Bergestm}
		0\le \varrho(u+v)-\varrho(u)-D_\varrho(u)[v] \le C_{(\varrho(u))}\|v\|^{\alpha \beta}_{LN(\varrho)}
	\end{equation}
	\item We also have, for such $u$,
	\begin{equation}\label{eq:diagestm}
	\varrho(u)\le D_\varrho(u)[u]\le \left(2(1+c)^{1/\alpha}-1\right)\varrho(u)
	\end{equation}
\end{enumerate} 
\end{cor}
\begin{proof}
	\begin{enumerate}
		\item Since $\varrho$ is convex, $t\mapsto \varrho(u+tv)$ has both right and left derivatives, and thus \eqref{eq:midestm} yields the differentiability. Linearity also easily follows from the convexity of $\varrho$ (see \cite[B, Proposition 4.2.1]{HL2001}, for example). 
		
		Now we prove the boundedness of $D_{\varrho}(u)$. If $\varrho(u)=\varrho(v)=1$, then by convexity,
		$\varrho(v)\ge \varrho(u)+D_{\varrho}(u)[v-u] \text{ and thus }D_\varrho(u)[u]\ge D_\varrho(u)[v]$; this ensures $\| D_\varrho(u)\|_{op}\le D_\varrho(u)[u]$. For general $u$ and $v$ with $\varrho(u)\ne0$, we obtain
		\begin{align}
			|D_\varrho(u)[v]|=&|\varrho(u)D_{\frac{\varrho}{\varrho(u)}}(u)[v]|\notag\\
			\le&\varrho(u)D_{\frac{\varrho}{\varrho(u)}}(u)[u]\|v\|_{LN(\varrho/\varrho(u))}\notag\\
			\le&D_\varrho(u)[u]\left(1\vee \left(\frac{\varrho(u)}{C}\right)^{-\frac{1}{\beta}} \right)\|v\|_{LN(\varrho)}.\label{eq:boundedeva}
		\end{align}
		If $\varrho(u)=0$ and $\alpha\le 1$, then $\alpha\beta>1$ forces $\beta>1$. Moreover, $\|u\|_{LN(\varrho)}=0$ by \eqref{eq:delta2} and hence $\varrho(u+tv)=\varrho(tv)$. Therefore
		$\varrho(u+tv)=\varrho(tv)=\varrho(-tv)=\varrho(u-tv)$ and differentiability forces $D_\varrho(u)[v]=0$.
		
		\item $\varrho(u-v)\ge \varrho(u)-D_\varrho(u)[v]$ by convexity. Together with \eqref{eq:midestm} and convexity, this yields \eqref{eq:Bergestm}.
		\item By convexity and \eqref{eq:uCLK}, we obtain
		\begin{multline*}
		\varrho(u)=\frac{\varrho(u)-\varrho(0)}{1}\le D_\varrho(u)[u]\\
		\le \frac{\varrho(2u)-\varrho(u)}{1}\le \left(2(1+c)^{1/\alpha}-1\right)\varrho(u).
		\end{multline*}		
	\end{enumerate}
\end{proof}
\begin{prop}\label{prop:uconvsmooth}
\begin{enumerate}
	\item If a gsm $\varrho$ satisfies \eqref{eq:lCLK} for some $\alpha, \beta\in(0,\infty)$, then the completion of $(\Mc{M}(\varrho)/\{\varrho=0\}, \|\cdot\|_{LN(\varrho)})$ is uniformly convex.
	\item If a gsm $\varrho$ satisfies \eqref{eq:uCLK} for some $\alpha,\beta\in(0,\infty)$, then the completion of $(\Mc{M}(\varrho)/\{\varrho=0\}, \|\cdot\|_{LN(\varrho)})$ is uniformly smooth.
\end{enumerate}
\end{prop}
\begin{proof}
\begin{enumerate}
	\item Note that, for $u\in\Mc{M}(\varrho)$ with $0<\|u\|_{LN(\varrho)}\le 1$ we have
	 \[1=\varrho\left(\frac{u}{\|u\|_{LN(\varrho)}} \right)\le \|u\|^{-\beta}_{LN(\varrho)}\varrho(u).\]
	 
	 Let $u,v\in\Mc{M}(\varrho)$ with $\|u\|_{LN(\varrho)}=\|v\|_{LN(\varrho)}=1$; recall that $\varrho(u)=\varrho(v)=1$. Then, by \eqref{eq:lCLK} we obtain
	\begin{align*}
	\varrho\left(\frac{u+v}{2}\right)\le&\left(\left(\frac{\varrho(u)+\varrho(v)}{2} \right)^\alpha-\varrho\left(\frac{u-v}{2}\right)^\alpha\right)^\frac{1}{\alpha}\\
	\le&\left(1-\left\|\frac{u-v}{2}\right\|^{\alpha \beta}_{LN(\varrho)}  \right)^\frac{1}{\alpha}
	\end{align*}
	because $\|\frac{u-v}{2}\|_{LN(\varrho)}\le 1$. Here,
	\[\left\|\frac{u+v}{2}\right\|_{LN(\varrho)}\le\left(1-2^{-\alpha \beta}\left\|u-v\right\|^{\alpha \beta}_{LN(\varrho)}  \right)^\frac{1}{\alpha \beta}\]
	holds, which ensures the desired uniform convexity.
	\item Let $u,v\in\Mc{M}(\varrho)$ with $\|u\|_{LN(\varrho)}=1$ and $\|v\|_{LN(\varrho)}\le 1/3$. Then,
	\begin{align*}
\left\|\frac{u+v}{\|u+v\|_{LN(\varrho)}}-u\right\|_{LN(\varrho)}&\le\frac{\left|1-\|u+v\|_{LN(\varrho)}\right|\|u\|_{LN(\varrho)}+\|v\|_{LN(\varrho)}}{\|u+v\|_{LN(\varrho)}}\\
&\le \frac{3}{2}\cdot 2\|v\|_{LN(\varrho)}\le 1.
	\end{align*}
	Hence, by \eqref{eq:Bergestm} we obtain
\[\varrho\left(\frac{u+v}{\|u+v\|_{LN(\varrho)}} \right)-\varrho(u)-D_\varrho(u)\left[\frac{u+v}{\|u+v\|_{LN(\varrho)}}-u \right]\le C_{(1)}3^{\alpha \beta}\|v\|^{\alpha \beta}_{LN(\varrho)},\]
which leads to 
\begin{equation}\label{eq:uv1}
(\|u+v\|_{LN(\varrho)}-1)D_\varrho(u)[u]-D_\varrho(u)[v]\le \frac{4}{3}C_{(1)} 3^{\alpha \beta}\|v\|^{\alpha \beta}_{LN(\varrho)}.
\end{equation}
By the same inequality for $-v$, we also have
\begin{equation}\label{eq:uv2}
(\|u-v\|_{LN(\varrho)}-1)D_\varrho(u)[u]+D_\varrho(u)[v]\le \frac{4}{3}C_{(1)} 3^{\alpha \beta}\|v\|^{\alpha \beta}_{LN(\varrho)}. 
\end{equation}
Since \eqref{eq:diagestm} holds, adding \eqref{eq:uv1} to \eqref{eq:uv2} we obtain
\[\|u+v\|_{LN(\varrho)}+\|u-v\|_{LN(\varrho)}-2\le 8C_{(1)}3^{\alpha \beta-1}\|v\|^{\alpha \beta}_{LN(\varrho)}. \]
This yields the desired uniform smoothness.
\end{enumerate}
\end{proof}

\begin{cor}\label{cor:clareflex}
If a gsm $\varrho$ satisfies \eqref{eq:lCLK} or \eqref{eq:uCLK} for some $\alpha,\beta\in(0,\infty)$, then $\varrho$ is reflexive.
\end{cor}
\begin{proof}
It follows from Proposition \ref{prop:uconvsmooth} with \cite[Proposition 1.e.3]{LT1979b}.
\end{proof}
\begin{cor}\label{cor:refsep}
	\begin{enumerate}
		\item If $\Mc{E}$ satisfies \eqref{eq:lCLK} for some $\alpha,\beta\in(0,\infty)$ and $q\ne1$\footnote{$q$ is the fixed number in Assumption \eqref{asm:lsc}.}, then $(\Mc{F}(\Mc{E}),\|\cdot\|_{LN({}_{sc}\Mc{E})}+\|\cdot\|_{L^q(\Mf{m})})$ is reflexive.
		\item If $\Mc{E}$ satisfies \eqref{eq:uCLK} for some $\alpha,\beta\in(0,\infty)$ with
		$(\frac{1}{q-1}\wedge 1\wedge \alpha)(\beta\wedge q)>1$ and $q\ne1$, then $(\Mc{F}(\Mc{E}),\|\cdot\|_{LN({}_{sc}\Mc{E})}+\|\cdot\|_{L^q(\Mf{m})})$ is reflexive.
	\end{enumerate}
In either case, this space is separable if $L^q(\Mf{m})$ is separable, in particular if $X$ is separable.
\end{cor}
\begin{proof}
	\begin{enumerate}
		\item Since $\|\cdot\|^q_{L^q(\Mf{m})}$ satisfies $(\lCla(\frac{1}{q-1}\vee 1,q))$ (see \cite{Cla1936} for example), $\Mc{E}_1$ satisfies $(\lCla(\frac{1}{q-1}\vee 1\vee\alpha,q\vee\beta))$ by the Minkowski inequality, and so does ${}_{sc}\Mc{E}_1$ by definition. Hence, $(\Mc{F}(\Mc{E}),\|\cdot\|_{LN({}_{sc}\Mc{E})}+\|\cdot\|_{L^q(\Mf{m})})$ is reflexive by Proposition \ref{prop:lscE} and Corollary \ref{cor:clareflex}.
		\item Since $\|\cdot\|^q_{L^q(\Mf{m})}$ satisfies $(\uBCL(\frac{1}{q-1}\wedge 1,q))$ with $c=C=1$ (see \cite{Cla1936} for example), $\Mc{E}_1$ satisfies $(\uBCL(\frac{1}{q-1}\wedge 1\wedge\alpha,q\wedge\beta))$ with the same $c,C$ as $\Mc{E}$ by the reverse Minkowski inequality. The rest of the proof is the same as (1).
	\end{enumerate}
If $L^q(\Mf{m})$ is separable, the continuous injection $\Mc{F}(\Mc{E})\hookrightarrow L^q(\Mf{m})$ and reflexivity imply separability, as in \cite[Corollary 3.16]{KS2024+}.
\end{proof}
We next prove detailed chain rules under \eqref{eq:uCLK}. We first derive some consequences of \eqref{eq:lCLK} and \eqref{eq:uCLK} for $\Mc{E}$, following the strategy of \cite[Section 7]{Sas2026}.
\begin{lem}\label{lem:Holdercontinuity}
	Let $\varrho$ be a gsm satisfying \eqref{eq:uCLK} for some $\alpha,\beta\in(0,\infty)$. Then, for $u,v,w\in\Mc{M}(\varrho)$ with $\|u-v\|_{LN(\varrho)}\le 1/2$ we have
	\begin{equation}\label{eq:holder}
		|D_\varrho(u)[w]-D_\varrho(v)[w]|\le (C_{(\varrho(u))}\vee C_{(\varrho(v))})2^{\alpha \beta}\|u-v\|^{\alpha \beta-1}_{LN(\varrho)}\|w\|_{LN(\varrho)}
	\end{equation}
	whenever $\varrho(u)\varrho(v)\ne 0$ or $\alpha\le 1$, where the constants  $C_{(\varrho(u))},C_{(\varrho(v))}$ are as in \eqref{eq:midestm}.
\end{lem}
We will give the proof of Lemma \ref{lem:Holdercontinuity} in a forthcoming paper.

\begin{lem}\label{lem:clinherit}
\begin{enumerate}
	\item If $\Mc{E}$ satisfies $(\lCla(1,\beta))$ for some $\beta$, then $\Mc{F}(\Mc{E})\ni f\mapsto\Gamma^c\langle f\rangle (A)$ satisfies $(\lCla(1,\beta))$ for any $A\in\Mc{B}(X)$, and $\Pi\langle f\rangle=\Gamma^c\langle f\rangle$ for any $f\in\Mc{M}(\Mc{E})$.
	\item  If $\Mc{E}$ satisfies \eqref{eq:uCLK} for $\alpha\in(0,1]$ and some $\beta$, then $\Mc{F}(\Mc{E})\ni f\mapsto\Gamma^c\langle f\rangle (A)$ satisfies \eqref{eq:uCLK} for any $A\in\Mc{B}(X)$, with the same constants $c,C$.
\end{enumerate}
\end{lem}
\begin{proof}
\begin{enumerate}
	\item Since the first condition of ($\lCla(1,\beta)$) for $\Mc{E}$ is equivalent to
	\[\Mc{E}\left(\frac{f+g}{2}\right)+\Mc{E}\left(\frac{f-g}{2}\right)\le \frac{\Mc{E}(f)+\Mc{E}(g)}{2} \]
	for any $f,g\in\Mc{M}(\Mc{E})$, we obtain
	\begin{equation}
		\Pi\left\langle\frac{f+g}{2}\right\rangle(A)+\Pi\left\langle\frac{f-g}{2}\right\rangle(A)\le \frac{\Pi\langle f\rangle(A) +\Pi\langle g\rangle(A)}{2}\label{eq:picla}
	\end{equation}
	for any $A\in\Mc{B}(X)$ by \eqref{eq:piinherit}. Moreover, this yields that $t\mapsto \Pi\langle f+tg\rangle(A)$ is midpoint convex. Since
	\[\Pi\langle f+tg\rangle(A)\le \Mc{E}(f+tg)\le\Mc{E}(f)+\Mc{E}(f+g)<\infty\]
	on $[0,1]$ by \ref{cond:rconv} and $\Mc{M}(\Mc{E})=\Mc{D}(\Mc{E})$, $t\mapsto \Pi\langle f+tg\rangle(A)$ is convex on $[0,1]$ by the Bernstein--Doetsch Theorem, as in Proposition \ref{prop:localtoconvex}. This yields $\Pi\langle f\rangle=\Gamma^c\langle f\rangle$. Hence, the first inequality of ($\lCla(1,\beta)$) for $(f\mapsto\Gamma^c\langle f\rangle (A))|_{\Mc{M}(\Mc{E})}$ follows from \eqref{eq:picla}, and the second one follows from Lemma \ref{lem:dom}. By the definition of $\Gamma^c$ for $f\in\Mc{F}(\Mc{E})$, the inequalities in ($\lCla(1,\beta)$) extend to $\Mc{F}(\Mc{E})$.
	\item For $\alpha\in(0,1]$, the function $\varphi(a,b)\coloneq 2(|a|^\alpha+c|b|^\alpha)^{1/\alpha}$ satisfies the assumptions of Lemma \ref{lem:boundinherit}~(2) by the reverse Minkowski inequality. Hence, \eqref{eq:uCLK} for $(f\mapsto\Gamma^c\langle f\rangle (A))|_{\Mc{M}(\Mc{E})}$ follows from Lemmas \ref{lem:dom} and \ref{lem:boundinherit}, with the same $c,C$, and extends to $\Mc{F}(\Mc{E})$ as in (1).
\end{enumerate}
\end{proof}
For a finite signed Borel measure $\eta$, we write $|\eta|$ for its total variation measure and $\|\eta\|_{\mathrm{TV}}\coloneq |\eta|(X)$.
\begin{lem}\label{lem:deriveHolder}
Assume that $\Mc{E}$ satisfies \eqref{eq:uCLK} for $\alpha\in(0,1]$ and $\beta\in(0,\infty)$. For any $f,g\in\Mc{F}(\Mc{E})$,
\[
D\Gamma^c\langle f;g\rangle(A)
\coloneq D_{\Gamma^c\langle\cdot\rangle(A)}(f)[g]
=\left.\frac{\der}{\der t}\Gamma^c\langle f+tg\rangle(A)\right|_{t=0}
\qquad(A\in\Mc{B}(X))
\]
defines a finite signed Borel measure satisfying $D\Gamma^c\langle f;g\rangle\ll\Gamma^c\langle f\rangle$ and $D\Gamma^c\langle f;g\rangle\ll\Gamma^c\langle g\rangle$.

For each fixed $f\in\Mc{F}(\Mc{E})$, the maps $g\mapsto D\Gamma^c\langle f;g\rangle$ and $g\mapsto D\Gamma^c\langle g;f\rangle$ are continuous from $(\Mc{F}(\Mc{E})/\{\|\cdot\|_{LN({}_{sc}\Mc{E})}=0\},\|\cdot\|_{LN({}_{sc}\Mc{E})})$ to the space of finite signed Borel measures equipped with $\|\cdot\|_{\mathrm{TV}}$. The first map is bounded and linear, and the second is locally $(\alpha\beta-1)$-H\"older continuous.
\end{lem}
\begin{proof}
Corollary \ref{cor:diff} and Lemma \ref{lem:clinherit}~(2) show that $D\Gamma^c\langle f;g\rangle$ is well-defined and finitely additive, and is linear in $g$. If $\Gamma^c\langle f\rangle(A)>0$, \eqref{eq:diagestm}, \eqref{eq:boundedeva} and
$\|g\|_{LN(\Gamma^c\langle\cdot\rangle(A))}\le\|g\|_{LN({}_{sc}\Mc{E})}$ give
\begin{multline}\label{eq:boundestmmeas}
|D\Gamma^c\langle f;g\rangle(A)|
\le\left(2(1+c)^{1/\alpha}-1\right)\\
\cdot\left(\Gamma^c\langle f\rangle(A)\vee
C^{1/\beta}\Gamma^c\langle f\rangle(A)^{1-1/\beta}\right)
\|g\|_{LN({}_{sc}\Mc{E})}.
\end{multline}
If $\Gamma^c\langle f\rangle(A)=0$, Corollary \ref{cor:diff} gives $D\Gamma^c\langle f;g\rangle(A)=0$, so the same bound holds. Since $\alpha\beta>1$ and $\alpha\le1$ imply $\beta>1$, this bound gives continuity at the empty set and boundedness on $\Mc{B}(X)$; here $\Gamma^c\langle f\rangle(X)<\infty$ by \eqref{eq:delta2}. Thus $D\Gamma^c\langle f;g\rangle$ is a finite signed Borel measure absolutely continuous with respect to $\Gamma^c\langle f\rangle$. If $\Gamma^c\langle g\rangle(A)=0$, Lemma \ref{lem:delta2} gives $\|g\|_{LN(\Gamma^c\langle\cdot\rangle(A))}=0$, and the boundedness of $D_{\Gamma^c\langle\cdot\rangle(A)}(f)$ gives $D\Gamma^c\langle f;g\rangle(A)=0$. This proves $D\Gamma^c\langle f;g\rangle\ll\Gamma^c\langle g\rangle$.

Applying \eqref{eq:boundestmmeas} to the two sets of a Hahn decomposition yields
\begin{align}\label{eq:boundestmTV}
&\|D\Gamma^c\langle f;g\rangle\|_{\mathrm{TV}}\notag\\
\le&2\left(2(1+c)^{1/\alpha}-1\right)
\left({}_{sc}\Mc{E}(f)\vee
C^{1/\beta}{}_{sc}\Mc{E}(f)^{1-1/\beta}\right)
\|g\|_{LN({}_{sc}\Mc{E})}.
\end{align}
Similarly, if $\|f-h\|_{LN({}_{sc}\Mc{E})}\le1/2$, apply Lemma \ref{lem:Holdercontinuity} to the two sets of a Hahn decomposition of $D\Gamma^c\langle f;g\rangle-D\Gamma^c\langle h;g\rangle$. Since the constants $C_{(s)}$ in \eqref{eq:holder} are independent of the Borel set and nondecreasing in $s$, we obtain
\begin{multline*}
\|D\Gamma^c\langle f;g\rangle-D\Gamma^c\langle h;g\rangle\|_{\mathrm{TV}}\\
\le2^{\alpha\beta+1}
\left(C_{({}_{sc}\Mc{E}(f))}\vee C_{({}_{sc}\Mc{E}(h))}\right)
\|f-h\|_{LN({}_{sc}\Mc{E})}^{\alpha\beta-1}
\|g\|_{LN({}_{sc}\Mc{E})}.
\end{multline*}
These estimates give the asserted continuity and show that both maps descend to the quotient. The local H\"older assertion follows because ${}_{sc}\Mc{E}$ is bounded on bounded sets for its Luxemburg seminorm, by \eqref{eq:delta2}.
\end{proof}
\begin{thm}[Chain rule for 2 variables]\label{thm:chain2}
Let $\varphi:\Mb{R}\to\Mb{R}$ be locally Lipschitz with $\varphi(0)=0$, and suppose that $\Mc{E}$ satisfies \eqref{eq:uCLK} for $\alpha\in(0,1]$ and $\beta\in(0,\infty)$. Then the following statements hold as identities or inequalities of Borel measures.
\begin{enumerate}
\item For any $f\in\Mc{F}(\Mc{E})$,
\[
\Gamma^c\langle f\rangle
\le D\Gamma^c\langle f;f\rangle
\le\left(2(1+c)^{1/\alpha}-1\right)\Gamma^c\langle f\rangle.
\]
\item For any $f,g\in\Mc{M}(\Mc{E})$,
\[
D\Gamma^c\langle f;\varphi\circ g\rangle
=(\varphi'\circ g)\,\der D\Gamma^c\langle f;g\rangle.
\]
\item Assume that $\Mc{E}$ satisfies \ref{cond:rstrloc}, and let $\Phi:\Mb{R}\to[0,\infty)$ be continuous.
\begin{enumerate}
\item If $|\Mc{E}(tf)-\Mc{E}(tg)|\le |t|\Phi(t)|\Mc{E}(f)-\Mc{E}(g)|$ for all $t\in\Mb{R}$ and $f,g\in\Mc{M}(\Mc{E})$, then, for any $f,g\in\Mc{M}(\Mc{E})$,
\begin{equation}\label{eq:2chian1val}
|D\Gamma^c\langle\varphi\circ f;g\rangle|
\le(\Phi\circ\varphi'\circ f)\,\der|D\Gamma^c\langle f;g\rangle|.
\end{equation}
\item If $|\Mc{E}(tf)-\Mc{E}(tg)|\ge |t|\Phi(t)|\Mc{E}(f)-\Mc{E}(g)|$ for all $t\in\Mb{R}$ and $f,g\in\Mc{M}(\Mc{E})$, then the reverse inequality in \eqref{eq:2chian1val} holds.
\end{enumerate}
\item If $\Mc{E}$ is $p$-homogeneous for $p\in(1,\infty)$, then
\[
D\Gamma^c\langle f;f\rangle=p\Gamma^c\langle f\rangle
\qquad(f\in\Mc{F}(\Mc{E})),
\]
and, for any $f,g\in\Mc{M}(\Mc{E})$,
\[
D\Gamma^c\langle\varphi\circ f;g\rangle
=(\varphi'\circ f)|\varphi'\circ f|^{p-2}\,\der D\Gamma^c\langle f;g\rangle,
\]
where $0|0|^{p-2}\coloneq 0$.
\end{enumerate}
\end{thm}
To recover the density formulation in \cite{Sas2026}, we additionally need a common dominating measure for the energy measures; see Appendix \ref{sec:medm}.
\begin{rem}
\begin{enumerate}
\item We set $\varphi'=0$ where the derivative does not exist. By Lemmas \ref{lem:eidp}, \ref{lem:deriveHolder} (specifically, the two absolute continuity assertions in the latter) and Corollary \ref{cor:clareflex},
\begin{gather*}
|D\Gamma^c\langle f;g\rangle|\bigl(f^{-1}(N)\cup g^{-1}(N)\bigr)=0,
\\ N\coloneq \{t:\varphi'(t)\text{ does not exist}\}.
\end{gather*}
Thus the weighted measures above are independent of this convention. Their weights are bounded on the ranges of $f,g\in\Mc{M}(\Mc{E})\subset C_c(X)$.
\item As in Theorem \ref{thm:chain1}, one can also formulate a version for two energies $\Mc{E}$ and $E$; we omit it here.
\end{enumerate}
\end{rem}
\begin{proof}[Proof of Theorem \ref{thm:chain2}]
For the composition statements, take $M>\|f\|_\infty\vee\|g\|_\infty$ and replace $\varphi(t)$ by $\varphi((-M)\vee t\wedge M)$. This preserves the compositions and the derivatives on the ranges of $f$ and $g$, and makes $\varphi$ Lipschitz. Corollaries \ref{cor:Mazsub} and \ref{cor:clareflex} then ensure that these compositions belong to $\Mc{F}(\Mc{E})$.
\begin{enumerate}
\item This follows from \eqref{eq:diagestm} and Lemma \ref{lem:clinherit}~(2).
\item Suppose first that $\varphi$ is piecewise affine, with $\varphi\circ g=\sum_{i=1}^k\alpha_i\eta_{[a_{i-1},a_i]}\circ g$ and $\|g\|_\infty<(-a_0)\wedge a_k$. Because $D\Gamma^c\langle f;g\rangle|_{g^{-1}(a)}=0$ for any $a\in\Mb{R}$ by \ref{cond:mlocal} and Lemma \ref{lem:deriveHolder}, differentiating \ref{cond:mlocal} and using linearity in the second variable gives
\[
D\Gamma^c\langle f;\varphi\circ g\rangle
=\sum_i\alpha_i D\Gamma^c\langle f;g\rangle\big|_{g^{-1}((a_{i-1},a_i))}
=(\varphi'\circ g)\,\der D\Gamma^c\langle f;g\rangle.
\]

For general $\varphi$, choose $\varphi_n$ as in Corollary \ref{cor:Mazsub}, applied to $g$. Lemma \ref{lem:deriveHolder} gives
\[
D\Gamma^c\langle f;\varphi_n\circ g\rangle
\longrightarrow D\Gamma^c\langle f;\varphi\circ g\rangle
\quad\text{in total variation}.
\]
By Lemma \ref{lem:eidp} and $D\Gamma^c\langle f;g\rangle\ll\Gamma^c\langle g\rangle$, we have $\varphi_n'\circ g\to\varphi'\circ g$ $|D\Gamma^c\langle f;g\rangle|$-a.e. Since $|\varphi_n'|\le\Lip(\varphi)$, $(\varphi_n'\circ g)\,\der D\Gamma^c\langle f;g\rangle\to(\varphi'\circ g)\,\der D\Gamma^c\langle f;g\rangle$ in total variation, by the dominated convergence theorem. This proves (2).
\item We first obtain the comparison for $\varphi(t)=at$.
Under the hypothesis of (a), the proof of Lemma \ref{lem:dom}~(2) gives
\[
|\Gamma^c\langle tf+stg\rangle(A)-\Gamma^c\langle tf\rangle(A)|
\le |t|\Phi(t)|\Gamma^c\langle f+sg\rangle(A)-\Gamma^c\langle f\rangle(A)|.
\]
Dividing by $|s|$ and letting $s\to0$ yield
\[
|D\Gamma^c\langle tf;g\rangle(A)|
\le\Phi(t)|D\Gamma^c\langle f;g\rangle(A)|
\qquad(t\ne0)
\]
for any $A\in\Mc{B}(X)$, while $D\Gamma^c\langle0f;g\rangle=0$. Taking total variations concludes (a) for this case.

For (b), if $\Mc{E}(u)>0$ for some $u\in\Mc{M}(\Mc{E})$, its hypothesis with $g=0$ and \eqref{eq:uCLK} imply
\[
|t|\Phi(t)\Mc{E}(u)\le\Mc{E}(tu)
\le C|t|^\beta\Mc{E}(u)\qquad(0<|t|\le1).
\]
Hence $\Phi(0)=0$ by continuity and $\beta>1$. If $\Mc{E}$ vanishes on $\Mc{M}(\Mc{E})$, all the measures in the comparison vanish.
Thus (b) follows in the same way when $\Phi(t)\ne0$, and is immediate otherwise.

For piecewise affine $\varphi$, apply the scalar comparisons on $f^{-1}((a_{i-1},a_i))$, using \ref{cond:mlocal} and Lemma \ref{lem:deriveHolder} as in (2). For general $\varphi$, take $\varphi_n$ from Corollary \ref{cor:Mazsub}, applied to $f$. Lemma \ref{lem:deriveHolder} implies convergence of $D\Gamma^c\langle\varphi_n\circ f;g\rangle$, and hence of their total variation measures, in total variation. On the right, $\Phi(\varphi_n'\circ f)\to\Phi(\varphi'\circ f)$ $|D\Gamma^c\langle f;g\rangle|$-a.e., and
\[
\sup_n|\Phi(\varphi_n'\circ f)|\le
\sup_{|t|\le\Lip(\varphi)}\Phi(t)<\infty.
\]
Dominated convergence therefore proves both (a) and (b).
\item If $\Mc{E}$ does not vanish on $\Mc{M}(\Mc{E})$, its $p$-homogeneity and \eqref{eq:uCLK} imply $p\ge\beta>1$. By Theorem \ref{thm:main1}~(3) and \ref{cond:mlim}, $\Gamma^c$ is $p$-homogeneous on $\Mc{F}(\Mc{E})$. Consequently,
\[
D\Gamma^c\langle f;f\rangle
=\left.\frac{\der}{\der t}(1+t)^p\Gamma^c\langle f\rangle\right|_{t=0}
=p\Gamma^c\langle f\rangle,
\]
and
\[
tD\Gamma^c\langle tf;g\rangle
=D\Gamma^c\langle tf;tg\rangle
=|t|^p D\Gamma^c\langle f;g\rangle
\qquad(t\in\Mb{R}).
\]
Thus $D\Gamma^c\langle tf;g\rangle=t|t|^{p-2}D\Gamma^c\langle f;g\rangle$ for $t\in\Mb{R}$. Locality and the same approximation argument as in (3), using continuity of $t\mapsto t|t|^{p-2}$, prove the composition identity.\qedhere
\end{enumerate}
\end{proof}
\begin{cor}[{cf. \cite[Corollary 7.9 and Theorem 7.10]{Sas2026}}]\label{cor:chain3}
Under the same assumptions as in Theorem \ref{thm:chain2}, the following statements hold.
\begin{enumerate}
	\item Let $n\in\Mb{N}$, $f,g_1,\dots,g_n\in\Mc{M}(\Mc{E})$, and let $\varphi:\Mb{R}^n\to\Mb{R}$ be a polynomial with $\varphi(\zero)=0$. Set $\Mf{g}=(g_1,\dots,g_n)$. Then $\varphi\circ\Mf{g}\in\Mc{F}(\Mc{E})$ and
	\begin{equation*}
		D\Gamma^c\langle f;\varphi\circ\Mf{g}\rangle
		=\sum_{i=1}^n\left(\partial_i\varphi\circ\Mf{g}\right)\der D\Gamma^c\langle f;g_i\rangle.
	\end{equation*}
	\item Assume in addition that, for every $n\in\Mb{N}$ and $R>0$, there exists $C_{n,R}<\infty$ such that
	\begin{equation}\label{eq:multicomposition}
		{}_{sc}\Mc{E}\bigl(\psi(g_1,\dots,g_n)\bigr)\le C_{n,R}
	\end{equation}
	whenever $g_1,\dots,g_n\in\Mc{M}(\Mc{E})$, $\max_i\Mc{E}(g_i)\le R$ and $\psi:\Mb{R}^n\to\Mb{R}$ is a polynomial satisfying $\psi(\zero)=0$ and $\Lip(\psi|_Q)\le 1$ for some closed cube $Q\subset \Mb{R}^n$ with $\{\zero\}\cup (g_1,\dots,g_n)(X)\subset Q$, with respect to the Euclidean norm. Then, the conclusion of (1) holds for every $\varphi\in C^1(\Mb{R}^n)$ with $\varphi(\zero)=0$.
	In particular, this assumption holds whenever $\Mc{E}$ satisfies \ref{cond:rstrloc}.
	\item If $\Mc{E}$ is $p$-homogeneous for $p\in(1,\infty)$, then for any $f,g\in\Mc{M}(\Mc{E})$,
	\[
	pg\,\der\Gamma^c\langle f\rangle
	=D\Gamma^c\langle f;fg\rangle
	-\left(\frac{p-1}{p}\right)^{p-1}D\Gamma^c\bigl\langle |f|^{\frac{p}{p-1}};g\bigr\rangle.
	\]
	In particular,
	\begin{equation*}
		\int_X g\,\der\Gamma^c{\langle f\rangle}
		=\frac{1}{p}D_{\Mc{E}^{(c)}}(f)[fg]-\frac{(p-1)^{p-1}}{p^p}D_{\Mc{E}^{(c)}}(|f|^{\frac{p}{p-1}})[g].
	\end{equation*}
\end{enumerate}
\end{cor}
\begin{proof}
By Lemma \ref{lem:NC} and $(f+g)^2-f^2-g^2=2fg$, we have $fg\in\Mc{F}(\Mc{E})$ for $f,g\in\Mc{M}(\Mc{E})$.
\begin{enumerate}
\item The claim follows from the identity for $fg$ above, Theorem \ref{thm:chain2}~(2) and linearity in the second variable; cf.\ \cite[Proof of Corollary 7.9 and Theorem 7.10]{Sas2026}.
\item Choose a closed cube $Q$ containing $\zero$ and $\Mf{g}(X)$, and polynomials $P_k$ with $P_k(\zero)=0$ such that $P_k\to\varphi$ and $\partial_iP_k\to\partial_i\varphi$ $(1\le i\le n)$ uniformly on $Q$, as in \cite[Proof of Theorem 7.10]{Sas2026}. Fix $R>0$ with $\max_i\Mc{E}(g_i)\le R$. Applying \eqref{eq:multicomposition} to the normalized polynomial differences gives
\[
\|(P_k-P_\ell)\circ\Mf{g}\|_{LN({}_{sc}\Mc{E})}
\le (1\vee C_{n,R})\Lip\bigl((P_k-P_\ell)|_Q\bigr)
\longrightarrow0
\;(k,\ell\to\infty).
\]
Now $P_k\circ \Mf{g}\to\varphi\circ \Mf{g}\ (k\to\infty)$ by the dominated convergence theorem. By completeness in Proposition \ref{prop:lscE}, $\varphi\circ\Mf{g}\in\Mc{F}(\Mc{E})$ and the convergence also holds in $\|\cdot\|_{LN({}_{sc}\Mc{E})}$. Part (1), Lemma \ref{lem:deriveHolder} and uniform convergence of the partial derivatives now yield the identity.

It remains to verify \eqref{eq:multicomposition} under \ref{cond:rstrloc}. In this case ${}_{sc}\Mc{E}=\Mc{E}^{(c)}$ on $\Mc{F}(\Mc{E})$ by Theorems \ref{thm:main1} and \ref{thm:emextension}. Lemma \ref{lem:deriveHolder} extends (1) to first arguments in $\Mc{F}(\Mc{E})$. For a polynomial $P$ with $P(\zero)=0$ and $\Lip(P|_Q)\le1$, set
\[
t=\left[4n(1\vee R)\left(2(1+c)^{1/\alpha}-1\right)\right]^{-1}.
\]
Since $|\partial_iP|\le1$ on $Q$ and $\|g_i\|_{LN(\Mc{E})}\le1\vee\Mc{E}(g_i)\le1\vee R$, part (1), Theorem \ref{thm:chain2}~(1) and \eqref{eq:boundestmTV} give
\begin{align*}
{}_{sc}\Mc{E}(tP\circ\Mf{g})
&\le D\Gamma^c\langle tP\circ\Mf{g};tP\circ\Mf{g}\rangle(X)\\
&\le t\sum_{i=1}^n
       \|D\Gamma^c\langle tP\circ\Mf{g};g_i\rangle\|_{\mathrm{TV}}\\
&\le \frac12\left({}_{sc}\Mc{E}(tP\circ\Mf{g})\vee
       C^{1/\beta}{}_{sc}\Mc{E}(tP\circ\Mf{g})^{1-1/\beta}\right).
\end{align*}
Hence ${}_{sc}\Mc{E}(tP\circ\Mf{g})\le C$, and Lemma \ref{lem:delta2} yields \eqref{eq:multicomposition} since $t$ depends only on $n,R$ and the fixed constants in \eqref{eq:uCLK}.
\item This follows in the same way as \cite[Corollary 7.9]{Sas2026}.\qedhere
\end{enumerate}
\end{proof}

\section{Examples}\label{sec:ex}
\subsection{$p$-energy forms}
We first introduce the notion of $p$-energy forms.
\begin{defi}[{$p$-energy form, \cite[Definition 3.1]{KS2024+}}]
Let $(X,\Mf{m})$ be a measure space, and denote by $L^0(X,\Mf{m})$ the space of $\Mb{R}$-valued measurable functions on $X$ modulo $\Mf{m}$-a.e.~equality. For $p\in(1,\infty)$, a $p$-homogeneous gsm $\Mc{E}$ on $\Mc{F}\subset L^0(X,\Mf{m})$ with $\Mc{F}=\Mc{M}(\Mc{E})$ is called a \emph{$p$-energy form}.
\end{defi}
\begin{rem}
Our definition and that in \cite{KS2024+} are equivalent. In fact, $\|f\|_{LN(\Mc{E})}=\Mc{E}(f)^{1/p}$ yields the triangle inequality for $\Mc{E}^{1/p}$, and conversely, the convexity of $t\mapsto |t|^p$ yields that a $p$-energy form in the sense of \cite{KS2024+} is a gsm. $p$-homogeneity follows from $\Mc{E}(a f)^{1/p}=|a|\Mc{E}(f)^{1/p}$.
\end{rem}
\begin{ex}\label{ex:pform}
	\begin{enumerate}
		\item Most typical examples are energies of the type
		\begin{multline*}
		\Mc{E}(f)=\int_{\Omega}\left(\sum_{1\le i,j\le d}a_{i,j}(x)\partial_i f\partial_j f\right)^\frac{p}{2}\der x\\
		+\int_{\Omega\times \Omega\setminus \mathrm{diag}}|f(x)-f(y)|^p j(x,y)\der x\der y+\int_\Omega |f(x)|^p \kappa(x)\der x,
		\end{multline*}
		where $\Omega$ is a domain in $\Mb{R}^d$, $p\in(1,\infty)$, $f\in\Mc{F}\coloneq W_0^{1,p}(\Omega)$, and $\mathrm{diag}=\{(x,x)\mid x\in\Omega\}$. We take $\Mf{m}=\der x$ and $q=p$. Assume that $A(x)=\{a_{i,j}(x)\}_{1\le i,j\le d}$ is a symmetric matrix-valued Borel function satisfying $\Lambda^{-1}I\le A\le\Lambda I$ a.e.~for some $\Lambda\ge1$, that $0\le\kappa\in L^\infty(\Omega)$, and that $j$ is a nonnegative symmetric Borel function satisfying
		\[
		\operatorname*{ess\,sup}_{x\in\Omega}\int_\Omega j(x,y)\der y<\infty.
		\]
		Then $\Mc{E}$ is a $p$-energy form on $\Mc{F}$, lower semicontinuous with respect to the $L^p(\Omega)$ norm, satisfying \ref{cond:rwNC}--\ref{cond:prereg}\footnote{in the sense of Section \ref{sec:lsc}}. Moreover, $\Mc{F}\cap C_c(\Omega)$ is dense in $\Mc{F}$ with respect to the $\|\cdot\|_{LN(\Mc{E}_1)}$ norm.
		\item Let $(\Mc{E},\Mc{F})$ be a (symmetric, bilinear) Dirichlet form on $L^2(X,\Mf{m})$, where $X$ is a locally compact separable metrizable space and $\Mf{m}$ is a full-support Radon measure (see \cite{FOT2011} for the definition). Then, $\Mc{E}(f)\coloneq \Mc{E}(f,f)$ is a lower semicontinuous $2$-energy form on $\Mc{F}$ with \ref{cond:rwNC}--\ref{cond:rstrsubad}. $(\Mc{E},\Mc{F})$ is called regular if \ref{cond:prereg} holds and $\Mc{F}\cap C_c(X)$ is dense with respect to the $\|\cdot\|_{LN(\Mc{E}_1)}$ norm.
		\item As mentioned in the introduction, strongly local $p$-energy forms with properties similar to those of Dirichlet forms are constructed on metric spaces, mainly on fractals (e.g.~\cite{AES2025+,CGQ2022,CGYZ2026,HPS2004,KS2026,Kig2023,KO2025+,MS2025,Shi2024}). The study of their properties is also an active area of research (e.g. \cite{Eri2026+,EM2025+,KS2026+,KS2024+,Sas2026,Yan2025+,Yan2025+2,Yan2025+3,Yan2025+4,Yan2026+2}).
		\item An important class of $p$-energy forms without \ref{cond:rstrloc} consists of \emph{traces of $p$-resistance forms}: see Subsection \ref{subsec:res}. Other (possibly) nonlocal examples appeared in \cite{Bau2024} (cf. \cite{CGYZ2024,MS2025,Yan2023+}), \cite{GYZ2023} and \cite[Theorem 6.22]{MS2025}. Properties of specific classes of $p$-energy forms without \ref{cond:rstrloc} are considered in \cite{CY2025+,GYZ2026,KS2024+,Yan2026+}.
	\end{enumerate}
\end{ex}
To discuss applicability of our theory, we recall the notion of \emph{generalized $p$-contraction property}, introduced in \cite{KS2024+}.
\begin{defi}\label{defi:GCP}
We say a $p$-energy form $\Mc{E}$ satisfies the generalized $p$-contraction property, $\Mr{(GC)}_p$ for short, if $(T_i (f))_{i=1}^m \in (\Mc{M}(\Mc{E}))^m$ and 
\[\left\|\left(\Mc{E}(T_i(f))^\frac{1}{p} \right)_{i=1}^m\right\|_{\ell^\alpha(\Mb{R}^m)}\le \left\|\left(\Mc{E}(f_j)^\frac{1}{p}\right)_{j=1}^n\right\|_{\ell^\beta(\Mb{R}^n)}  \]
for any $f=(f_j)_{j=1}^n\in(\Mc{M}(\Mc{E}))^n$, whenever $m,n\in\Mb{N}, \alpha\in[p,\infty], \beta\in(0,p]$ and $T=(T_i)_{i=1}^m:\Mb{R}^n\to\Mb{R}^m$ with $T(0)=0$ and $\|T(x)-T(y)\|_{\ell^\alpha(\Mb{R}^m)}\le \|x-y\|_{\ell^\beta (\Mb{R}^n)}$.
\end{defi}
Note that $\Mr{(GC)}_p$ implies the extra assumption in Corollary \ref{cor:chain3}~(2).

\begin{defi}
Let $X,V,\varrho$ be as in Definition \ref{defi:condgsm}. We define the following condition, which we apply with $\Mc{E}=\varrho$.
\begin{description}
\rcond{cond:rweakd2} $\varrho(2x_n)\to 0\ (n\to\infty)$ whenever  $\varrho(x_n)\to 0\ (n\to\infty)$.
\end{description}
\end{defi}
Condition \ref{cond:rweakd2} is called the \emph{weak $\Delta_2$-condition} and, as the name suggests, is weaker than the $\Delta_2$ condition \eqref{eq:delta2}.

\begin{lem}\label{lem:condaccum}
	 Let $X$ be a locally compact metrizable space and $\Mf{m}$ be a full-support Radon measure on $X$. We also let $\Mc{E}$ be a lower semicontinuous gsm on $V\subset L^q(\Mf{m})\ (q\in[1,\infty))$ satisfying $\Mc{M}(\Mc{E})=\Mc{D}(\Mc{E})$ or \ref{cond:rweakd2}. If $\Mc{E}$ satisfies \ref{cond:rwNC} and \ref{cond:rstrloc}, or \ref{cond:rwNC} and \ref{cond:rfaulting}, then it also enjoys \ref{cond:rvanish}.
\end{lem}
\begin{proof} We first assume $\Mc{M}(\Mc{E})=\Mc{D}(\Mc{E})$. If $f\in\Mc{M}(\Mc{E})\cap C_c(X)$, then
	\begin{multline*}
	\Mc{E}(f-\eta_{[-2t,2t]}\circ f)+\Mc{E}(\eta_{[0,t]}\circ f)\\=\Mc{E}(f-\eta_{[-2t,2t]}\circ f+\eta_{[0,t]}\circ f )
	\le \Mc{E}(f)<\infty
	\end{multline*}
	for $t>0$, whenever \ref{cond:rwNC} and \ref{cond:rstrloc} hold. Since $\|\eta_{[-2t,2t]}\circ f\|_{L^q(\Mf{m})}\to 0$ as $t\downarrow0$ by the dominated convergence theorem, we obtain $\lim_{t\downarrow0}\Mc{E}(\eta_{[0,t]}\circ f)=0$ by the lower semicontinuity of $\Mc{E}$. Similarly, if \ref{cond:rwNC} and \ref{cond:rfaulting} hold, then we obtain
	\begin{multline*}
	 \Mc{E}(\eta_{[0,t]}\circ f)+\Mc{E}(f-\eta_{[-t,t]}\circ f)\le \Mc{E}(\eta_{[-t,t]}\circ f)+\Mc{E}(f-\eta_{[-t,t]}\circ f)\\
	 =\Mc{E}(-t\vee f\wedge t)+\Mc{E}(f-t\vee 0\wedge f+t)
	 \le \Mc{E}(f)+\Mc{E}(0)=\Mc{E}(f)<\infty
	\end{multline*}
	for $t>0$. Similarly to the former case, this yields $\lim_{t\downarrow0}\Mc{E}(\eta_{[0,t]}\circ f)=0$.
	
	Now we turn to the case where \ref{cond:rweakd2} holds. Let $f\in\Mc{M}(\Mc{E})\cap C_c(X)$ and $\lambda\in(0,1)$ with $\lambda f\in\Mc{D}(\Mc{E})$. Then, we have $\lim_{t\downarrow0}\Mc{E}(\eta_{[0,t]}\circ (\lambda f))=\lim_{t\downarrow0}\Mc{E}(\lambda\eta_{[0,\lambda^{-1} t]}\circ f)=0$, as above. This leads to $\lim_{t\downarrow0}\Mc{E}(\eta_{[0,t]}\circ f)=0$ by \ref{cond:rweakd2}.
\end{proof}
\begin{cor}\label{cor:pefapply}
 Let $X$ be a locally compact metrizable space, $\Mf{m}$ be a full-support Radon measure on $X$, and $q\in[1,\infty)$. Let $\Mc{E}$ be a $p$-energy form on a dense linear subspace $\Mc{F}=\Mc{M}(\Mc{E})$ of $L^q(\Mf{m})$, lower semicontinuous on $\Mc{F}$ with respect to the $L^q(\Mf{m})$ norm, and let ${}_{sc}\Mc{E}$ denote its lower semicontinuous envelope on $L^q(\Mf{m})$. Assume that an $\Mb{R}$-linear subspace $\Mc{C}\subset\Mc{M}({}_{sc}\Mc{E})\cap C_c(X)$ is dense in $(\Mc{M}({}_{sc}\Mc{E}),\|\cdot\|_{LN({}_{sc}\Mc{E})}+\|\cdot\|_{L^q(\Mf{m})})$ and that $({}_{sc}\Mc{E})|_{\Mc{C}}$ satisfies \ref{cond:prereg}. If $({}_{sc}\Mc{E})|_{\Mc{C}}$ also satisfies one of the following conditions, then ${}_{sc}\Mc{E}$ with core $\Mc{C}$ satisfies the assumptions of Theorem \ref{thm:main2}.
 \begin{enumerate}
 	\item \ref{cond:rwNC} and \ref{cond:rstrloc}.
 	\item $\Mr{(GC)}_p$.
 \end{enumerate}
 If $\Mc{E}$ is reflexive, it suffices instead to take an $\Mb{R}$-linear subspace $\Mc{C}\subset\Mc{F}\cap C_c(X)$, dense in $(\Mc{F},\|\cdot\|_{LN(\Mc{E})}+\|\cdot\|_{L^q(\Mf{m})})$, such that $\Mc{E}|_{\Mc{C}}$ satisfies \ref{cond:prereg} and either (1) or (2).
\end{cor}
\begin{proof}
	The envelope ${}_{sc}\Mc{E}$ is $p$-homogeneous. First apply the following arguments to $({}_{sc}\Mc{E})|_{\Mc{C}}$.
	\begin{enumerate}
		\item We need to prove \ref{cond:rvanish} and \ref{cond:rstrsubad}; these follow from Proposition \ref{prop:localmeansnldc} and Lemma \ref{lem:condaccum}.\footnote{Apply the statements in reverse order. Note that every $p$-energy form satisfies \ref{cond:rweakd2}.}
		\item \ref{cond:rwNC} and \ref{cond:rstrsubad} follow from \cite[Proposition 2.3~(b),(c)]{KS2024+}. Moreover, \cite[Proposition 2.3~(c)]{KS2024+} also leads to \ref{cond:rfaulting} with $\varphi(t)=-a\vee t\wedge a$. Hence, Lemma \ref{lem:condaccum} yields \ref{cond:rvanish}.
	\end{enumerate}
	Now assume that $\Mc{E}$ is reflexive and that $\Mc{C}\subset\Mc{F}\cap C_c(X)$ satisfies the stated conditions. The same arguments show that $\Mc{E}|_{\Mc{C}}$ satisfies \ref{cond:rwNC}--\ref{cond:prereg}. Relative lower semicontinuity gives ${}_{sc}\Mc{E}=\Mc{E}$ on $\Mc{F}$, so these conditions also hold for $({}_{sc}\Mc{E})|_{\Mc{C}}$. The assumed regularity gives ${}_{sc}(\Mc{E}|_{\Mc{C}})={}_{sc}\Mc{E}$ and makes $\Mc{E}|_{\Mc{C}}$ reflexive. Proposition \ref{prop:lscE}~(2) and Remark \ref{rem:lscincl}, applied to $\Mc{E}|_{\Mc{C}}$, therefore yield the required regularity of $\Mc{C}$ in $(\Mc{M}({}_{sc}\Mc{E}),\|\cdot\|_{LN({}_{sc}\Mc{E})}+\|\cdot\|_{L^q(\Mf{m})})$.
\end{proof}
\begin{rem}\label{lem:pemapply}
\begin{enumerate}
	\item $p$-energy forms mentioned in Example \ref{ex:pform}~(1)-(3), extended by $\infty$ outside their domains and with regularity assumed in (2), satisfy the assumptions of Theorem \ref{thm:main2}, and their $p$-energy measures, though they were already known, are recovered. In particular, \cite[Theorem 2.4]{Sas2026} is essentially recovered by Theorems \ref{thm:main2} and \ref{thm:chain2}~(4). For the lower $p$-Clarkson inequality when $1<p<2$, we use the argument for (CI1) in \cite[Proof of Theorem 6.2]{Sas2026}. Additionally, the assumptions for a Local $p$-Dirichlet space in \cite{Eri2026+,EM2025+} can be stated solely in terms of conditions on $p$-energy forms, if $\Mc{E}$ is reflexive (cf.~\cite[Lemma 2.10]{EM2025+}).
	\item Thanks to \cite[Remark 8.20]{KS2024+}, the $p$-energy form $\Mc{E}_p^\Gamma$ in \cite[Theorem 6.22]{MS2025} satisfies the assumption of Theorem \ref{thm:main2}.
	In the case of the standard planar Sierpi\'nski carpet $K$,
	equipped with its normalized Euclidean metric and uniform
	self-similar probability measure $\Mf{m}$, the strongly local
	part is nonzero for every $p\in(1,\infty)$. Indeed, let
	$\Mc{E}_p$ and $\Mc{F}_p\subset L^p(K,\Mf{m})$ be as in
	\cite[Theorem 1.1]{MS2025}, and extend both energies by $\infty$
	outside $\Mc{F}_p$.
	By \cite[Theorems 1.1~(iii) and 6.22~(i)]{MS2025}, there exists
	$a>0$ such that $a\Mc{E}_p\le\Mc{E}_p^\Gamma$.
	Since $\Mc{E}_p$ is strongly local, Theorem \ref{thm:main2}~(1)
	yields
	\[
	a\Mc{E}_p(f)\le(\Mc{E}_p^\Gamma)^{(c)}(f)
	\qquad(f\in\Mc{F}_p).
	\]
	Thus $(\Mc{E}_p^\Gamma)^{(c)}(f)>0$ whenever
	$f\in\Mc{F}_p$ is not constant $\Mf{m}$-a.e.,
	by \cite[Theorem 1.1~(vi)]{MS2025}.

	To the author's knowledge, this is the first construction of a strongly local part and its energy measures for a possibly nonlocal $p$-energy with $p\ne2$ whose strongly local part is not apparent from its definition. For integral-type $p$-energies, such as $[u]^p_{B_{p,p}^\sigma}$ in \cite{GYZ2023,GYZ2026} or $\Mc{E}^{(\Upsilon,p)}$ in \cite{Yan2026+}, the corresponding strongly local part must be trivial, because $\lim_{n\to\infty}\Mc{E}(T_n\circ f)=0$ follows from the dominated convergence theorem. Similarly, for \emph{mixed local and nonlocal $p$-energy forms} defined in \cite{CY2025+}, their strongly local part coincides with $\Mc{E}^{(L)}$ in \cite[Definition 2.3]{CY2025+}.
	\item We discuss the application of our theory to the $p$-resistance forms and their traces in the following subsection. In other cases in Example \ref{ex:pform}~(4) which are not mentioned above, \ref{cond:rstrsubad} may fail.
\end{enumerate}
\end{rem}
\subsection{$p$-resistance forms}\label{subsec:res}
\begin{defi}[{$p$-resistance form, \cite[Definition 6.1]{KS2024+}}]
Let $p\in(1,\infty)$, $\Mc{F}\subset \Mb{R}^X$ and $\Mc{E}$ be a $p$-homogeneous gsm on $\Mc{F}=\Mc{M}(\Mc{E})$. The pair $(\Mc{E},\Mc{F})$ is called a \emph{$p$-resistance form}, if the following properties hold.
\begin{description}
	\RFcond{cond:RFconst} $\indi_X\in\Mc{F}$ and $\{\Mc{E}=0\}=\Mb{R}\indi_X$.
	\RFcond{cond:RFbanach} $(\Mc{F}/\Mb{R}\indi_X, \Mc{E}^{1/p})$ is a Banach space. \footnote{$\|f\|_{LN(\Mc{E})}=\Mc{E}(f)^{1/p}$ by homogeneity.}
	\RFcond{cond:RFsep} Let $x,y\in X$. If $x\ne y$, then $f(x)\ne f(y)$ for some $f\in\Mc{F}$.
	\RFcond{cond:RFres} For any $(x,y)\in X^2$,
	\[R_\Mc{E}(x,y)\coloneq \sup\left\{\frac{|f(x)-f(y)|^p}{\Mc{E}(f)}\ \bigg\vert\ f\in\Mc{F}, \Mc{E}(f)\ne 0 \right\}<\infty. \]
	\RFcond{cond:RFgcp} $\Mc{E}$ enjoys $\Mr{(GC)}_p$ as a $p$-energy form on $(X,\#)$, where $\#$ is the counting measure.
\end{description}

\end{defi}
\begin{rem}
For the original definition, $p$-homogeneity follows from \ref{cond:RFgcp}.
\end{rem}
The notion of ($2$-)resistance form was introduced and developed by Kigami~\cite{Kig1995,Kig2001,Kig2012}. It is an abstraction of the effective resistance of electrical networks, and is now widely used as a basic tool in analysis of ``low-dimensional'' fractals (e.g.~planar fractals, $n$-dimensional Sierpi\'nski gaskets) and random graphs/trees. In discrete settings, the notion of $p$-resistance for $p\ne2$ has already been studied in research on networks and/or graph theory (e.g.~\cite{AL2011,ACFP2019,Her2010}).
\begin{thm}[{\cite[Theorem 1.6]{KS2024+}, cf. \cite{ACFP2019,Her2010}}]
Let $(\Mc{E},\Mc{F})$ be a $p$-resistance form on $X$. Then, $(X,R_\Mc{E}^{1/(p-1)})$ is a metric space.
\end{thm}
Note that $\Mc{F}\subset C(X,R_\Mc{E}^{1/(p-1)})$ because
\begin{equation}\label{eq:resestm}
|f(x)-f(y)|\le\Mc{E}(f)^\frac{1}{p}R_{\Mc{E}}(x,y)^\frac{1}{p}
\end{equation}
follows from \ref{cond:RFres}.
\begin{prop}\label{prop:appliRF}
Let $(\Mc{E},\Mc{F})$ be a $p$-resistance form on $X$, and $\Mf{m}$ be a full-support Radon measure on $(X,R_\Mc{E}^{1/(p-1)})$. If $(X,R_\Mc{E}^{1/(p-1)})$ is locally compact and $\Mc{E}$ satisfies \ref{cond:prereg} (as a gsm on $(X,R_\Mc{E}^{1/(p-1)})$), then $\Mc{E}|_{C_c(X)}$ satisfies the assumptions of Theorem \ref{thm:main1}. Moreover, $\Mc{E}|_{C_c(X)}$ is reflexive and lower semicontinuous with respect to $\|\cdot\|_{L^q(\Mf{m})}$ for $q\in[1,\infty)$ (in the sense of identification in Remark \ref{rem:abbrinj}).
\end{prop}
\begin{proof}
As in Corollary \ref{cor:pefapply}~(2), \ref{cond:RFgcp} yields \ref{cond:rwNC}, \ref{cond:rstrsubad} and \ref{cond:rfaulting}. It also yields reflexivity, by \cite[Proposition 2.3~(e),(f)]{KS2024+} and Corollary \ref{cor:clareflex}. It remains to prove lower semicontinuity and \ref{cond:rvanish}.

Let $f_n,f\in\Mc{F}\cap C_c(X)$ with $f_n\to f$ in $L^q(\Mf{m})$. Passing to a subsequence, we may assume that $\Mc{E}(f_n)\to\liminf_{k\to\infty}\Mc{E}(f_k)<\infty$ and $f_n\to f$ $\Mf{m}$-a.e. Choose $x_0$ such that $f_n(x_0)\to f(x_0)$. By reflexivity of $(\Mc{F}/\Mb{R}\indi_X,\Mc{E}^{1/p})$, Mazur's lemma and \ref{cond:RFbanach}, there exist a Mazur subsequence $\{g_n\}_{n\ge1}$ of $\{f_n\}_{n\ge1}$ and $g\in\Mc{F}$ such that $g(x_0)=f(x_0)$ and $\Mc{E}(g_n-g)\to0$. Since $g_n(x_0)\to f(x_0)$, \eqref{eq:resestm} gives
\[
|g_n(y)-g(y)|
\le |g_n(x_0)-f(x_0)|+\Mc{E}(g_n-g)^{1/p}R_{\Mc{E}}(x_0,y)^{1/p}
\longrightarrow0
\]
for every $y\in X$. This forces $f=g$, and therefore convexity yields
\[
\Mc{E}(f)=\lim_{n\to\infty}\Mc{E}(g_n)
\le\liminf_{n\to\infty}\Mc{E}(f_n).
\]
Finally, Lemma \ref{lem:condaccum}, applied to $\Mc{E}|_{C_c(X)}$, yields \ref{cond:rvanish}.
\end{proof}

\begin{defi}[Trace of a $p$-resistance form]\label{defi:trace}
Let $(\Mc{E},\Mc{F})$ be a $p$-resistance form on $X$, and $B\subset X$ be nonempty. We define
\begin{align*}
	\Mc{F}|_B&=\{f\in \Mb{R}^B\mid g|_B=f\text{ for some }g\in\Mc{F}\},\\
	\Mc{E}|_B(f)&=\inf\{\Mc{E}(g)\mid g\in\Mc{F}\text{ such that }g|_B=f \}\quad (f\in\Mc{F}|_B).
\end{align*}
The pair $(\Mc{E}|_B,\Mc{F}|_B)$ is called the trace of $(\Mc{E},\Mc{F})$ on $B$.
\end{defi}
\begin{rem}
Here we do not try to define the trace of general $p$-energy forms and gsms in our framework, because this notion probably requires a more delicate treatment of capacities, which is beyond the scope of this paper.
\end{rem}
\begin{thm}[{Part of \cite[Theorem 6.13]{KS2024+}}]\label{thm:trace}
$(\Mc{E}|_B,\Mc{F}|_B)$ is a $p$-resistance form on $B$. Moreover, $R_{\Mc{E}|_B}=R_{\Mc{E}}|_{B\times B}$.
\end{thm}
\begin{lem}\label{lem:traceinherit}
Let $(\Mc{E},\Mc{F})$ be a $p$-resistance form on $X$, and $B\subset X$ be a nonempty closed subset. If $(\Mc{E},\Mc{F})$ satisfies \ref{cond:prereg} and $(X,R_\Mc{E}^{1/(p-1)})$ is locally compact, then $(\Mc{E}|_B,\Mc{F}|_B)$ also satisfies \ref{cond:prereg} and $(B,R_{\Mc{E}|_B}^{1/(p-1)})$ is locally compact.
\end{lem}
\begin{proof}
By Theorem \ref{thm:trace}, the inheritance of local compactness is trivial. For \ref{cond:prereg}, let $u\in C_c(B)$ be arbitrary. By Tietze's extension theorem, there exists $v\in C(X)$ with $v|_B=u$. Since $\supp(u)$ is also compact as a subset of $X$ and $X$ is locally compact, there exists $\tau\in C_c(X)$ such that $\tau|_{\supp(u)}=1$. Then, $\tau v\in C_c(X)$ and $\tau v|_B=u$. For any $\epsilon>0$, there exists $f\in\Mc{F}\cap C_c(X)$ with $\|\tau v-f\|_{\infty}<\epsilon$, by \ref{cond:prereg} for $(\Mc{E},\Mc{F})$. Then, $f|_B\in(\Mc{F}|_B)\cap C_c(B)$, since $\supp(f|_B)\subset B\cap\supp(f)$, and $\|u-f|_B\|_{\infty}=\|\tau v|_B-f|_B\|_{\infty}<\epsilon$, which concludes the proof.
\end{proof}
By Proposition \ref{prop:appliRF}, Theorem \ref{thm:trace} and Lemma \ref{lem:traceinherit} (under the assumption of Proposition \ref{prop:appliRF}) we can apply Theorems \ref{thm:main1}, \ref{thm:emextension}, \ref{thm:chain1} and \ref{thm:chain2} to $p$-resistance forms and their traces on nonempty closed sets (cf.~\cite[Chapter 9]{Kig2012}). These applications do not require $X$ or the energy seminorm quotient to be separable. In particular, thanks to \cite[Theorems 8.30 and 8.43]{KS2024+}, our theory is also applicable to the trace of $p$-resistance forms in \cite{CGQ2022,Kig2023}.
\subsection{Nonlinear Dirichlet forms}
As another nonlinear generalization of Dirichlet forms, here we introduce the notion of nonlinear Dirichlet forms, defined by Cipriani and Grillo~\cite[Definition 3.1]{CG2003}. We retain their density requirement and use the equivalent characterization in \cite[Theorem 1.3]{BH2024}.
\begin{defi}[{Nonlinear Dirichlet form\footnote{In \cite{CG2003}, this notion is simply called ``Dirichlet form'', and the classical Dirichlet form is called ``quadratic Dirichlet form''. This terminology has the advantage that our ``nonlinear Dirichlet form'' can be bilinear, but it is confusing. In \cite{BH2024}, this notion is called ``non-bilinear Dirichlet form''. We adopt the terminology of \cite{Bri2026}.} (cf. \cite[Theorem 1.3]{BH2024})}]\label{defi:NLDF}
Let $(X,\Omega,\Mf{m})$ be a $\sigma$-finite measure space. A lower semicontinuous functional $\Mc{E}:L^2(\Mf{m})\to[0,\infty]$ is called a \emph{nonlinear Dirichlet form} if $\{\Mc{E}<\infty\}$ is dense in $L^2(\Mf{m})$ and
\begin{align*}
	\Mc{E}(f\vee g)+\Mc{E}(f\wedge g)&\le \Mc{E}(f)+\Mc{E}(g),\\
	\Mc{E}(g-a\vee f\wedge g+a)+\Mc{E}(f-a\vee g\wedge f+a)&\le \Mc{E}(f)+\Mc{E}(g)
\end{align*}
for any $a>0$ and $f,g\in L^2(\Mf{m})$.
\end{defi}
Note that a nonlinear Dirichlet form $\Mc{E}$ is a gsm if $\Mc{E}(0)=0$ and $\Mc{E}$ is symmetric, namely, $\Mc{E}(f)=\Mc{E}(-f)$ for every $f\in L^2(\Mf{m})$.
\begin{prop}\label{prop:AppliNLDF}
Let $X$ be a locally compact metric space, $\Mf{m}$ be a $\sigma$-finite, full-support Radon measure and $\Mc{E}$ be a nonlinear Dirichlet form on $L^2(\Mf{m})$. Assume that $\Mc{E}$ is symmetric with $\Mc{E}(0)=0$, that \ref{cond:prereg} holds, and that $\Mc{M}(\Mc{E})\cap C_c(X)$ is dense in $(\Mc{M}(\Mc{E}), \|\cdot\|_{LN(\Mc{E})}+\|\cdot\|_{L^2(\Mf{m})})$.
Then, $\Mc{E}$ satisfies the assumptions of Theorem \ref{thm:main2} with $q=2$ whenever \ref{cond:rweakd2} holds or $\Mc{M}(\Mc{E})=\Mc{D}(\Mc{E})$.
\end{prop}
\begin{proof}
\ref{cond:rwNC} follows from \cite[Theorem 1.2]{BH2024}. \ref{cond:rstrsubad} and \ref{cond:rfaulting} hold by definition. Thus, Lemma \ref{lem:condaccum} yields \ref{cond:rvanish}.
\end{proof}

\begin{ex}\label{ex:NLDF}
\begin{enumerate}
\item Let $\Omega$ be a bounded domain. The energy in
Example \ref{ex:pform}~(1), with the coefficients specified there,
defines a symmetric nonlinear Dirichlet form on $L^2(\Omega)$ with
$\Mc{D}(\Mc{E})=W_0^{1,p}(\Omega)\cap L^2(\Omega)$ and
$\Mc{E}=\infty$ elsewhere. Proposition \ref{prop:AppliNLDF} applies.
Thus, these standard examples overlap with $p$-energy forms, with
a different choice of ambient space and effective domain.

\item Let $X=(0,1)$ with Lebesgue measure, and define
\[
\Mc{E}(u)=\int_0^1\bigl(|u'|^2+|u'|^3\bigr)\der x,
\qquad \Mc{D}(\Mc{E})=W_0^{1,3}(0,1),
\]
with $\Mc{E}=\infty$ elsewhere in $L^2(0,1)$.
This is an inhomogeneous, symmetric, strongly local nonlinear
Dirichlet form with Orlicz growth $t\mapsto t^2+t^3$
(cf.~\cite[Theorem 4.11]{CG2003}).
We can check that Proposition \ref{prop:AppliNLDF} applies and that
$\uBCL(2/3,2)$ holds with $c=2$ and $C=1$.
Its energy measure is
$\Gamma^c\langle u\rangle=(|u'|^2+|u'|^3)\der x$.

\item Let $\Omega\subset\Mb{R}^d$ be a bounded $C^1$ domain,
and let $p:\overline\Omega\to(1,\infty)$ satisfy
$1<p_-\coloneq \inf_{\overline\Omega}p\le p_+\coloneq \sup_{\overline\Omega}p<\infty$,
with $1/p$ locally log-H\"older continuous on $\overline\Omega$
in the sense of \cite[Definition 4.1.1]{DHHR2011}.
Define on $L^2(\Omega)$
\[
\Mc{E}(u)=\int_\Omega\frac{|\nabla u(x)|^{p(x)}}{p(x)}\der x,
\qquad
\Mc{D}(\Mc{E})=W_0^{1,p(\cdot)}(\Omega)\cap L^2(\Omega),
\]
and set $\Mc{E}=\infty$ elsewhere.
This is a symmetric, strongly local nonlinear Dirichlet form
\cite[Proposition 6.1]{SZ2025+}, \cite[Example 7.1]{CC2025}.
The smooth core furnished by \cite[Corollary 11.2.4]{DHHR2011},
together with bounded truncation, gives the core condition in
Proposition \ref{prop:AppliNLDF}, which applies since
$\Mc{E}(2u)\le2^{p_+}\Mc{E}(u)$.
Pointwise estimates and their preservation under integration yield
\eqref{eq:uCLK} with
\[
\alpha=\min\{1,2/p_+\},\qquad
\beta=p_-,\qquad c=\max\{1,p_+-1\},\qquad C=1,
\]
whenever $\alpha p_->1$.

\item Analogous energies satisfying \eqref{eq:uCLK} can be considered
on compact smooth sub-Riemannian manifolds
\cite[Theorem 4.5 and Example 4.6]{CG2003} and on compact smooth
reversible, strongly convex Finsler manifolds
\cite[Section 2.3]{Bri2026}.
Further examples of nonlinear Dirichlet forms include Cheeger energies
\cite[Theorem 2.4]{Bri2026}, discrete $p$-energies
\cite[Proposition 3.9]{Mug2013}, and nonlocal energies \cite{CW2026+}.
For each realization, the domain and core assumptions in
Proposition \ref{prop:AppliNLDF} must be checked.
\end{enumerate}
\end{ex}

\subsection{Nonlinear resistance forms}\label{subsec:nonlinearres}
Finally, we consider the nonlinear resistance forms introduced by Puchert and Schmidt.
\begin{defi}[{Nonlinear resistance form, \cite[Definition 4.11]{PS2025+}}]\label{defi:NRF}
Let $X$ be a nonempty set. A convex functional $\Mc{E}:\Mb{R}^X\to[0,\infty]$ with $\Mc{E}(0)=0$ is called a \emph{nonlinear resistance form} if it is lower semicontinuous with respect to pointwise convergence and satisfies
\[
\sup\{f(x)-f(y)\mid f\in\Mb{R}^X,\ \Mc{E}(f)\le1\}<\infty
\]
for any $x,y\in X$, and
\begin{equation}\label{eq:NRFcontraction}
\Mc{E}(u+C\circ v)+\Mc{E}(v-C\circ u)
\le\Mc{E}(u+v)+\Mc{E}(u-v)
\end{equation}
for any $u,v\in\Mb{R}^X$ and normal contraction $C:\Mb{R}\to\Mb{R}$, i.e., $C(0)=0$ and $\Lip(C)\le1$.
\end{defi}
\begin{prop}\label{prop:appliNRF}
Let $X$ be a locally compact metrizable space and $\Mc{E}$ be a symmetric nonlinear resistance form on $X$ satisfying \ref{cond:prereg}. If $\Mc{M}(\Mc{E})=\Mc{D}(\Mc{E})$ or \ref{cond:rweakd2} holds, then $\Mc{E}|_{C_c(X)}$ satisfies the assumptions of Theorem \ref{thm:main1}.
\end{prop}
\begin{proof}
Convexity, symmetry and lower semicontinuity imply that $\Mc{E}$ is a gsm. Taking $(u,v)=(0,f)$ in \eqref{eq:NRFcontraction} yields \ref{cond:rwNC}, and taking $(u,v)=(\frac{f+g}{2},\frac{f-g}{2})$ with $C(s)=|s|$ in \eqref{eq:NRFcontraction} yields \ref{cond:rstrsubad}. The corresponding domain inclusions follow by scaling.

Similarly, for $f\in\Mc{D}(\Mc{E})$, applying \eqref{eq:NRFcontraction} with $(u,v)=(\frac{0\vee f}{2},\frac{0\vee f}{2})$ and $C(s)=|s-2/t|-2/t$ gives
\[
\Mc{E}(0\vee (f-t))+\Mc{E}(0\vee f\wedge t)\le\Mc{E}(0\vee f)\le\Mc{E}(f)<\infty
\]
for $t>0$. Since $0\vee(f-t)\to 0\vee f$ pointwise as $t\downarrow0$, lower semicontinuity yields ${\Mc{E}(0\vee f\wedge t)}\to0$. This proves \ref{cond:rvanish} when $\Mc{M}(\Mc{E})=\Mc{D}(\Mc{E})$. If \ref{cond:rweakd2} holds, apply this argument to $2^{-k}f\in\Mc{D}(\Mc{E})$ for $f\in\Mc{M}(\Mc{E})\cap C_c(X)$ and sufficiently large $k$, and then use \ref{cond:rweakd2}.
\end{proof}
By \cite[Proposition 5.11]{PS2025+} and Theorem \ref{thm:trace}, the $p$-resistance forms and their traces in Subsection \ref{subsec:res}, extended by $\infty$ outside their domains, are nonlinear resistance forms. 

\section{Remaining problems and future directions}\label{sec:arg}
\subsection{Properties of jump and killing parts}
We are still far from obtaining the desired nonlinear counterparts of the Beurling--Deny decomposition. Regrettably, we do not even know whether $\Mc{E}-\Mc{E}^{(c)}$ is convex or not.

A related problem is whether $\Pi=\Gamma^c$ holds. If it is true, we can avoid convex modification in the definition of the strongly local part, and we obtain $\Gamma^c\langle f\rangle(\{g\le a\})=\lim_{n\to\infty}\Mc{E}(f_n^{g,a})$ in typical cases. This is more reasonable and manageable in analysis. In this paper, we have only shown $\Pi=\Gamma^c$ in the strongly local case (Proposition \ref{prop:localtoconvex}) and the case where $(\lCla(1,\beta))$ holds (Lemma \ref{lem:clinherit}~(1)).

Solving the second problem may help to resolve the first one. Conversely, to investigate properties of $\Mc{E}-\Mc{E}^{(c)}$ will give clues to resolve the second problem. For example, if we obtain $\Mc{E}(f)-\Mc{E}^{(c)}(f)\le C\|f\|_{\infty}^2$ uniformly for some $C>0$, this yields $\Pi=\Gamma^c$, and, under Assumption \ref{asm:lsc}, the lower semicontinuity of $\Mc{E}^{(c)}$ on $(\Mc{M}(\Mc{E}),\|\cdot\|_{L^q(\Mf{m})})$. (When $p=2$, this uniform bound holds when $\Mf{k}, J$ in \eqref{eq:BDdecomp} are finite measures, but is not true in general, even when $p=2$.)

\subsection{Further relaxation of assumptions of main results}
We achieved considerable relaxation of conditions to obtain a canonical energy measure of the strongly local part, compared to \cite{Sas2026}. Can we pursue this direction more generally? To consider this problem, we introduce a toy model
\[\Mc{E}(f)\coloneq \left(\int_0^1 |\nabla f|^2\right)^\frac{p}{2} \quad (f\in W^{1,2}([0,1]))\]
for $p\in[1,\infty)$. Then $\Mc{E}$ is a gsm on $L^2([0,1])$ with \ref{cond:rwNC}, \ref{cond:rvanish} and \ref{cond:prereg}, but fails \ref{cond:rstrsubad} when $p\ne 2$. In particular, we cannot decompose $\Mc{E}$ into
\begin{align*}
\Mc{E}(f)={}&\Mc{E}^{(c)}(f)+\int_{[0,1]^2\setminus\Mr{diag}}|f(x)-f(y)|^p J(\der x,\der y)+\int_{[0,1]}|f(x)|^p\Mf{k}(\der x)
\end{align*}
for some positive Borel measures $J$ on $[0,1]^2\setminus\Mr{diag}$ and $\Mf{k}$ on $[0,1]$ and a gsm $\Mc{E}^{(c)}$ with \ref{cond:rstrloc} and $\Mc{E}^{(c)}(T_n\circ f)=\Mc{E}^{(c)}(f)$ when $p\ne2$. Indeed, if this holds, then the jump and killing terms for $T_n\circ f$ tend to zero by the dominated convergence theorem. Hence, $\Mc{E}(f)=\lim_{n\to\infty}\Mc{E}(T_n\circ f)=\lim_{n\to\infty}\Mc{E}^{(c)}(T_n\circ f)=\Mc{E}^{(c)}(f)$, which contradicts the fact that $\Mc{E}$ fails \ref{cond:rstrloc} when $p\ne2$.

Simply put, the ``local part'' of $\Mc{E}$ might be 
\[\tilde\Gamma^c\langle f\rangle(A)\coloneq \left(\int_A |\nabla f|^2 \right)^\frac{p}{2}, \]
though $\tilde\Gamma^c\langle f\rangle$ is not a measure when $p\ne2$. To justify this intuition may enable further generalization of our theory.

\subsection{Orlicz--Sobolev spaces on fractals}
Using the methods of analysis on fractals and nonlinear Dirichlet forms, it may be possible to consider Orlicz--Sobolev spaces on fractals, to which our theory can be applied. Let $(X,d)$ be a locally compact metric space and $\Mf{m}$ be a Radon measure, as in Section \ref{sec:ex}. For example, in \cite{KS2026}, the Korevaar--Schoen-type $p$-energy $\Mc{E}^{\Mr{KS}}_p$ for $p\in(1,\infty)$ is defined as
\[\Mc{E}^{\Mr{KS}}_p(f)= \lim_{n\to\infty}\int_X\int_{B_d(x,r_n)}\frac{|f(x)-f(y)|^p}{r_n^{\beta_p}\Mf{m}(B_d(x,r_n))}\der \Mf{m}(y)\der\Mf{m}(x)\]
for some $\beta_p$ and $\{r_n\}_{n\ge 1}$ with $r_n\to0\ (n\to\infty)$, under a ``weak monotonicity estimate'', including fractal settings. Since $\{r_n\}_{n\ge 1}$ can be chosen as a subsequence of any given positive sequence approaching zero, we may naively consider the following candidate, provided that a common sequence can be chosen along which the mixed limit exists on the intended domain:
\begin{align*}
&\Mc{E}^{\Mr{KS}}_{p_1,p_2}(f)\\
=&\lim_{n\to\infty}\int_X\int_{B_d(x,r_n)}\max_{p=p_1,p_2}\frac{|f(x)-f(y)|^p}{r_n^{\beta_p}\Mf{m}(B_d(x,r_n))}\der\Mf{m}(y)\der\Mf{m}(x).
\end{align*}
However, this energy possibly becomes trivial if $\{\Mc{E}^{\Mr{KS}}_{p_1}<\infty\}\cap \{\Mc{E}^{\Mr{KS}}_{p_2}<\infty\}=\{\text{a.e.~constant functions}\}$, and it actually occurs on some fractals (see \cite{AES2025+}). So, to show nontriviality may be difficult on fractals.

Similarly, we may also consider a variable-exponent Sobolev space on fractals, which is another candidate that we can study using canonical energy measures.

\subsection{On inheritance of $(\Mr{GC})_p$}
In \cite{KS2026+}, the author and Kajino proved that, under the assumptions of \cite[Theorem 2.4]{Sas2026} except for the $p$-Clarkson inequalities, $f\mapsto\Gamma^c\langle f\rangle(A)$ satisfies the full $(\Mr{GC})_p$ for every $A\in\Mc{B}(X)$, provided that $\Mc{E}$ satisfies the condition in Definition \ref{defi:GCP} for $\alpha=\beta=p$ and linear $T$. However, one of the key techniques relies heavily on the property \ref{cond:rstrloc}, and a similar proof cannot be applied to the nonlocal case. So, whether $\Mc{E}^{(c)}$ and $f\mapsto \Gamma^c\langle f\rangle(A)$ satisfy $(\Mr{GC})_p$ under the assumption that $\Mc{E}$ satisfies $(\Mr{GC})_p$ remains open.

\appendix
\section{Further relations among conditions on gsms}\label{ap:further}
In this appendix, we check some equivalences or implications between conditions on gsms and discuss related results.

\begin{lem}[{cf. \cite[Lemmas 3.2 and 4.3]{Sas2026}}]\label{lem:ewUC}
Let $V\subset C_b(X)$ be a vector space for some topological space $X$, and let $\varrho$ be a gsm on $V$ with \ref{cond:rbvanish} and \ref{cond:rstrlocbd}. Then $\varrho$ satisfies \ref{cond:rwNC}.
\end{lem}
\begin{proof}
Fix any $f\in\Mc{M}(\varrho)$, and let $\varphi$ be as in \ref{cond:rwNC}. For $a<b$ with $a,b\ne0$, we first note that
\[\eta_{[a,b]}(t)=\begin{cases}
	\eta_{[0,b]}(t)-\eta_{[0,a]}(t)& (0<a<b),\\
	\eta_{[0,b]}(t)-\eta_{[0,-a]}(-t) &(a<0<b),\\
	\eta_{[0,-b]}(-t)-\eta_{[0,-a]}(-t)&(a<b<0).
\end{cases}\]
Thus, by \ref{cond:rbvanish}, applied to $f$ and $-f$, $\eta_{[a,b]}\circ f\in \Mc{M}(\varrho)$ for any $-\infty<a<b<\infty$. Hence $\varphi\circ f\in\Mc{M}(\varrho)$.

We prove $\varrho(\varphi\circ f)\le\varrho(f)$ similarly to the proof of Lemma \ref{lem:locfold}. Write
\[\varphi\circ f=\sum_{i=1}^n \alpha_i\eta_{[a_{i-1},a_i]}\circ f,\]
with $a_0<a_1<\dots <a_n$, $a_0<-\|f\|_{\infty}$, $\|f\|_{\infty}<a_n$ and $|\alpha_i|\le 1$ for any $i$. We may assume $a_m=0$ for some $m\in\{1,\dots,n-1\}$ without loss of generality.
Note that, for sufficiently small $t>0$,
\begin{multline*}
\left\|\sum_{i=1}^{m} \alpha_i\eta_{[a_i-t,a_i]}\circ f+\sum_{i=m}^{n-1}\alpha_{i+1}\eta_{[a_i,a_i+t]}\circ f\right\|_{LN(\varrho)}\\
\le \sum_{i=1}^{m}|\alpha_i|\left\|\eta_{[a_i-t,a_i]}\circ f\right\|_{LN(\varrho)}\\
+\sum_{i=m}^{n-1} |\alpha_{i+1}|\left\|\eta_{[-a_i-t,-a_i]}\circ (-f)\right\|_{LN(\varrho)}\to 0\quad(t\downarrow0)
\end{multline*}
by \ref{cond:rbvanish} and Lemma \ref{lem:continuity}. In the same way,
\[\lim_{t\downarrow 0}\left\|\sum_{i=1}^m \eta_{[a_i-t,a_i]}\circ f+\sum_{i=m}^{n-1} \eta_{[a_i,a_i+t]}\circ f\right\|_{LN(\varrho)}=0.\]
By \ref{cond:rlc}, it is sufficient to consider $f\in \Mc{D}^i(\varrho)$. Then,
\begin{align*}
	\varrho(\varphi\circ f)\le&\liminf_{t\downarrow 0}\varrho\left(\sum_{i=1}^{m} \alpha_i\eta_{[a_{i-1},a_i-t]}\circ f+\sum_{i=m}^{n-1} \alpha_{i+1}\eta_{[a_i+t,a_{i+1}]}\circ f \right)\\
	\overset{\ref{cond:rstrlocbd}}{=}&\liminf_{t\downarrow0}\left(\sum_{i=1}^m\varrho\left( \alpha_i\eta_{[a_{i-1},a_i-t]}\circ f \right)+\sum_{i=m}^{n-1}\varrho\left( \alpha_{i+1}\eta_{[a_i+t,a_{i+1}]}\circ f \right)\right)\\
	\overset{\ref{cond:rconv},\ref{cond:rsym}}{\le}&\liminf_{t\downarrow0}\left(\sum_{i=1}^m\varrho\left(\eta_{[a_{i-1},a_i-t]}\circ f \right)+\sum_{i=m}^{n-1}\varrho\left(\eta_{[a_i+t,a_{i+1}]}\circ f \right)\right)\\
	\overset{\ref{cond:rstrlocbd}}{=}&\lim_{t\downarrow 0}\varrho\left(\sum_{i=1}^{m}\eta_{[a_{i-1},a_i-t]}\circ f+\sum_{i=m}^{n-1}\eta_{[a_i+t,a_{i+1}]}\circ f \right)=\varrho(f)
\end{align*}
by Lemmas \ref{lem:lsc} and \ref{lem:limeval}.
\end{proof}

\begin{cor}
Let $X$ be a locally compact metrizable space, $\Mf{m}$ be a $\sigma$-finite, full-support Radon measure on $X$, and $\Mc{E}$ be a gsm on $L^2(\Mf{m})$. If all the following conditions hold, then $\Mc{E}$ is a symmetric nonlinear Dirichlet form.
\begin{itemize}
	\item $\Mc{E}$ is lower semicontinuous with respect to the $L^2(\Mf{m})$ norm.
	\item $\Mc{D}(\Mc{E})$ is dense in $L^2(\Mf{m})$.
	\item $\Mc{D}(\Mc{E})\cap C_b(X)$ is dense in $\Mc{D}(\Mc{E})$ with respect to $\|\cdot\|_{LN(\Mc{E}_1)}$.
	\item $\Mc{E}|_{L^2(\Mf{m})\cap C_b(X)}$ satisfies \ref{cond:rstrlocbd} and \ref{cond:rbvanish}.
	\end{itemize}
\end{cor}
\begin{rem}
	For the definitions of $\Mc{E}_1$ and extensions of conditions for gsms on $L^2(\Mf{m})$, see Section \ref{sec:lsc}.
\end{rem}
\begin{proof}
By Definition \ref{defi:NLDF}, it suffices to show \ref{cond:rstrsubad} and \ref{cond:rfaulting}. First restrict $\Mc{E}$ to $L^2(\Mf{m})\cap C_b(X)$. Lemma \ref{lem:ewUC} gives \ref{cond:rwNC}, and Proposition \ref{prop:localmeansnldc} gives both required inequalities.

By \ref{cond:rlc}, it suffices to consider $f,g\in\Mc{D}^i(\Mc{E})$. Take $f_n,g_n\in\Mc{D}(\Mc{E})\cap C_b(X)$ such that $f_n\to f$ and $g_n\to g$ in $\|\cdot\|_{LN(\Mc{E}_1)}$. These convergences also hold in $L^2(\Mf{m})$, so lower semicontinuity and the inequalities for bounded functions yield
\begin{multline*}
\Mc{E}(f\vee g)+\Mc{E}(f\wedge g)\le\liminf_{n\to\infty}\left(\Mc{E}(f_n\vee g_n)+\Mc{E}(f_n\wedge g_n) \right)\\
\le\lim_{n\to\infty}\left(\Mc{E}(f_n)+\Mc{E}(g_n) \right)=\Mc{E}(f)+\Mc{E}(g).
\end{multline*}
Similarly, for every $a>0$,
\begin{multline*}
\Mc{E}((f-a)\vee g\wedge(f+a))+\Mc{E}((g-a)\vee f\wedge(g+a))\\
\le\liminf_{n\to\infty}\left(\Mc{E}((f_n-a)\vee g_n\wedge(f_n+a))+\Mc{E}((g_n-a)\vee f_n\wedge(g_n+a))\right)\\
\le\lim_{n\to\infty}\left(\Mc{E}(f_n)+\Mc{E}(g_n)\right)=\Mc{E}(f)+\Mc{E}(g).
\end{multline*}
The last equalities follow from Lemma \ref{lem:limeval} and $\Mc{E}\le\Mc{E}_1$.
\end{proof}


\section{Minimal energy-dominant measures}\label{sec:medm}
Throughout this appendix, we assume Assumption \ref{asm:lsc}. We give a sufficient condition for the existence of a common dominating measure for the energy measures.

\begin{defi}[{Minimal energy-dominant measure; cf.\ \cite[Definition 2.1]{Hin2010}}]\label{defi:medm}
Let $\nu$ be a $\sigma$-finite Borel measure on $X$. We call $\nu$ a \emph{minimal energy-dominant measure} of $(\Mc{E}^{(c)},\Mc{F}(\Mc{E}))$ if the following conditions hold:
\begin{description}
\DMcond{cond:medmdomination} $\Gamma^c\langle f\rangle\ll\nu$ for every $f\in\Mc{F}(\Mc{E})$.
\DMcond{cond:medmminimality} If a $\sigma$-finite Borel measure $\mu$ on $X$ satisfies $\Gamma^c\langle f\rangle\ll\mu$ for every $f\in\Mc{F}(\Mc{E})$, then $\nu\ll\mu$.
\end{description}
If only \ref{cond:medmdomination} is required, we call $\nu$ an \emph{energy-dominant measure}.
\end{defi}
Any two minimal energy-dominant measures are mutually absolutely continuous.

\begin{lem}[{cf.\ \cite[Lemmas 2.2 and 2.3]{Hin2010}}]\label{lem:medmexist}
If
\[
\left(\Mc{F}(\Mc{E})/\{\|\cdot\|_{LN({}_{sc}\Mc{E})}=0\},
\|\cdot\|_{LN({}_{sc}\Mc{E})}\right)
\]
is separable, then there exists a finite minimal energy-dominant measure of $(\Mc{E}^{(c)},\Mc{F}(\Mc{E}))$.
\end{lem}
\begin{proof}
By Remark \ref{rem:lscincl}, the image of $\Mc{M}(\Mc{E})$ is dense in the given quotient. We may therefore choose $f_i\in\Mc{M}(\Mc{E})$, $i\ge1$, whose images form a countable dense subset. Theorem \ref{thm:main1}~(2) gives
\[
\Gamma^c\langle f_i\rangle(X)=\Mc{E}^{(c)}(f_i)<\infty.
\]
Define
\begin{equation}\label{eq:medmseries}
\nu\coloneq \sum_{i=1}^{\infty}
\frac{2^{-i}}{1+\Mc{E}^{(c)}(f_i)}\Gamma^c\langle f_i\rangle.
\end{equation}
Then $\nu$ is a Borel measure with $\nu(X)\le1$ and $\Gamma^c\langle f_i\rangle\ll\nu$ for every $i$.

For $f\in\Mc{F}(\Mc{E})$, choose a sequence $\{f_{i_k}\}_{k\ge1}$ from this family such that $\|f-f_{i_k}\|_{LN({}_{sc}\Mc{E})}\to0$. If $A\in\Mc{B}(X)$ and $\nu(A)=0$, then \ref{cond:mlim} yields
\[
\Gamma^c\langle f\rangle(A)
=\lim_{k\to\infty}\Gamma^c\langle f_{i_k}\rangle(A)=0.
\]
Thus \ref{cond:medmdomination} holds. If $\mu$ is any $\sigma$-finite energy-dominant measure, then \eqref{eq:medmseries} yields $\nu\ll\mu$. This proves \ref{cond:medmminimality}.
\end{proof}
\begin{rem}
\begin{enumerate}
\item Lemma \ref{lem:medmexist} does not require \eqref{eq:delta2}. The finite measures used in \eqref{eq:medmseries} are chosen from the original core $\Mc{M}(\Mc{E})$.
\item If $L^q(\Mf{m})$ is separable and either hypothesis of Corollary \ref{cor:refsep} is satisfied, the hypothesis of Lemma \ref{lem:medmexist} follows: the quotient map from
\[
\left(\Mc{F}(\Mc{E}),\|\cdot\|_{LN({}_{sc}\Mc{E})}+\|\cdot\|_{L^q(\Mf{m})}\right)
\]
onto the energy-seminorm quotient is continuous.
\end{enumerate}
\end{rem}

If an energy-dominant measure exists, the Radon--Nikodym density versions of Lemma \ref{lem:deriveHolder} and Theorem \ref{thm:chain2} follow immediately under their respective assumptions.

\subsection*{Acknowledgments}
The author would like to thank Sylvester Eriksson-Bique and Mathav Murugan for their stimulating suggestions concerning the non-$p$-homogeneous extension of the results in \cite{Sas2026} and the relaxation of the $p$-Clarkson inequality assumption.
The author is also grateful to Ryosuke Shimizu for helpful comments on an earlier version of the manuscript.
The author would also like to thank the organizers and participants of the Workshop on Potential Theory and Geometric Function Theory in Metric Spaces at ICMAT for stimulating discussions, with particular thanks to the participant who pointed out the significance of condition \ref{cond:rvanish}.
This work was supported by JSPS KAKENHI Grant Numbers JP24KJ0022 and JP26K17012.
\subsection*{Declaration of use of AI tools}
The author wrote most of the first draft of this paper and developed its logical structure and proofs, excluding the initial proof of Proposition \ref{prop:uconvsmooth}~(2), which was generated by OpenAI's GPT-5.6 Sol.
The first draft of Example \ref{ex:NLDF} was generated by OpenAI's GPT-6 Astra.

In refining the manuscript, the author used GPT-6 Astra for a range of tasks, including identifying and correcting errors in the proofs and logical structure, improving the results and the organization of the manuscript based on the author's ideas, and proofreading.

The author has reviewed the entire manuscript and takes full responsibility for its content.
\bibliographystyle{abbrv}
\bibliography{bib/references}

\end{document}